\documentclass[11pt,a4paper]{article} 

\usepackage{adjustbox}
\usepackage{algorithm}
\usepackage[noend]{algpseudocode}
\usepackage{amsfonts,amsmath,amssymb,amsthm}
\usepackage{array}
\usepackage{authblk}
\usepackage{bbm}
\usepackage{booktabs}
\usepackage{dsfont}
\usepackage{enumerate}
\usepackage{enumitem}
\usepackage{float}
\usepackage{graphicx}
\usepackage{hyperref}
\usepackage[utf8]{inputenc}
\usepackage{latexsym}
\usepackage{lineno}
\usepackage{lmodern}
\usepackage{makecell}
\usepackage{multicol}
\usepackage{multirow}
\usepackage{scalerel}
\usepackage{subfigure}
\usepackage[skins]{tcolorbox}
\usepackage{textcomp}
\usepackage{tikz}
\usepackage{xcolor}
\usepackage{xspace}

\theoremstyle{plain}
\newtheorem{theorem}{Theorem}[section]

\newtheorem{remark}[theorem]{Remark}

\theoremstyle{definition}
\newtheorem{definition}[theorem]{Definition}

\newcommand{\MATLAB}{\textsc{Matlab}\xspace}
\def\NoNumber#1{{\def\alglinenumber##1{}\State #1}\addtocounter{ALG@line}{-1}}

\title{Data-driven quasi-interpolant spline surfaces for point cloud approximation}
\author[1,2]{Andrea Raffo}
\author[3]{Silvia Biasotti}
\affil[1]{Department of Mathematics and Cybernetics, SINTEF, Oslo, Norway.}
\affil[2]{Department of Mathematics, University of Oslo, Oslo, Norway.}
\affil[3]{Istituto di Matematica Applicata e Tecnologie Informatiche ``E. Magenes'' CNR, Genova, Italy.}
\date{}                     
\begin{document}
\maketitle

\begin{abstract}
    In this paper we investigate a local surface approximation, the \emph{Weighted Quasi Interpolant Spline Approximation} (wQISA), specifically designed for large and noisy point clouds. We briefly describe the properties of the wQISA representation and introduce a novel data-driven implementation, which combines prediction capability and complexity efficiency. We provide an extended comparative analysis with other continuous approximations on real data, including different types of surfaces and levels of noise, such as 3D models, terrain data and digital environmental data.
\end{abstract}

\section{Introduction}
Due to the recent progress in the acquisition by laser scanners, photogrammetry and diagnostic devices and the popularity of remote sensing technologies, 
there has been an exponential growth in the availability of data and the need of efficient representations. A good representation model must be, at the same time, \emph{efficient}, that is based on the minimum amount of data, yet \emph{effective}, that is able to support conservative conclusions and keep as much information as possible. To handle large data volumes it is often necessary to adopt approximation strategies so that a single point becomes representative of a region or a set of properties \cite{Daehlen1993,Ackermann2015}. 
Given a set of measurements (and their spatial locations), its model approximation has to permit deducing information about the process that generated those data, even at locations different from those at which the measurements were obtained. Surface reconstruction, terrain elevation estimation and the approximation of rainfall and pollution fields over Digital Elevation Models (DEMs) are all examples of spatial data approximation. Moreover, an efficient approximation allows the recovery of the digital representation of a physical shape that commonly contains a variety of properties and defects, such as: geometric features; noise and outliers; perturbations introduced  when the acquisition conditions are not optimal (low resolution of instruments, motions, etc.) or different acquisition techniques collect data sets of different resolution (laser scans or photogrammetry); incomplete data (e.g, in the acquisition of broken artefacts or scans of objects partially occluded \cite{Georgopoulos2017}).

Quasi-interpolation schemes \cite{deboor73,Sablonniere05} are popular for data approximation because, unlike traditional least squares approximations, they do not require solving a linear system. We use B-splines as basis functions because they are computationally convenient, for instance with respect to the common radial basis functions, as (piecewise) polynomial bases require low-order integration schemes to be exactly computed. The use of piecewise algebraic approximations makes our method suitable also to CAD applications, where B-splines and NURBS are \emph{de facto} the standard tools. Moreover, the treatment of essential boundary conditions is more natural for structured approaches like the spline-based ones than in meshless methods, as it relies on the number of repetition of each knot value. 

The wQISA implementation proposed in this paper is specifically designed to define powerful prediction methods from low quality points. Starting from the theoretical wQISA definition given in \cite{Raffo2019} for a single tensor-product mesh, in this work we derive a multi-level approximation algorithm and provide an extended comparative analysis with other methods in the literature. 
The main contributions of the paper include: a detailed description of the wQISA method when used for surface approximation; a multi-level approximation algorithm based on a data-driven definition of the weight functions; an extensive comparison of the wQISA outcome with other well-known continuous approximation methods.

The remainder of the paper is organized as follows. Section \ref{previouswork} overviews the literature on data approximation focusing on methods related to continuous surface approximation. Section \ref{wQISA} overviews the wQISA method and its properties, together with a multi-level implementation of the method based on a data-driven mesh refinement strategy. Section \ref{simulation} extensively compares the wQISA outcome on different real use cases, with respect to a number of well-known approximation methods and a set of performance indicators that are detailed in Section \ref{sec:settings}. Discussions and concluding remarks are provided in Section \ref{sec:conclusions}.

\section{Previous Work\label{previouswork}}

The literature on data approximation is vast and we cannot do justice to all contributions. We limit our review to methods whose output satisfies some smoothness requirement, either local or global. For a more complete list of methods related to our use cases,  we refer to \cite{Berger:2017} for a recent survey on surface reconstruction methods, to \cite{garnero2013comparisons} for an overview on methods for modelling terrain data and to \cite{Patane2017} for a comparative analysis of methods for approximating rainfall data.

\textbf{Meshless methods}.
Kriging is a geo-statistical approach largely used to approximate remote sensing measurements \cite{Oliver1990}. In particular, (ordinary) Kriging incorporates correlation information into data approximation through a variogram model. The main limitation of ordinary Kriging is the limited scalability; indeed, Kriging's computation scales quadratically with respect to the number of observations. The Moving Least Square (MLS) approach is largely adopted for surface reconstruction \cite{Alexa2001,Shen2004,Fleishman2005}. The basic idea behind MLS is to approximate the surface in the neighbour of a point with the tangent plane in that point. Several variations have been introduced, for instance Feng et al. \cite{Feng2014} devised a model that uses multiple curves/surfaces approximation, that allows separating mixed scanning points received from a thin-wall object and consists of a second-order extension of the method in \cite{Levin2003}.  Other implicit approximation techniques express the data as linear combination of basis elements. For instance, radial basis functions (RBFs) are central for scattered data approximation \cite{Franke}. The quality of the interpolation depends on the choice of the basis function. Originally introduced as a global method \cite{Savchenko1995,Carr2001}, RBFs have been adapted to compact support \cite{FLOATER1996}. A combination of MLS and RBF is presented in \cite{PATANE2012387}. Wavelets and wavelet transforms \cite{chui2016introduction} have been applied for data compression and noise reduction by truncating the wavelet decomposition. Among other techniques, we mention Poisson based methods, which have been applied for example to surface reconstruction from point sets \cite{Kazhdan2006}.

\textbf{Mesh-based methods}. Tensor-product spline surfaces (piecewise polynomial or NURBS) are a well established representation for modelling smooth shapes (see, for example, \cite{Schumaker2007,Farin1996}). Several local refinement methods have been proposed to overcome the limited mesh adaptivity to the data, such as T-splines \cite{Sederberg2003}, locally refined (LR) B-splines \cite{Dokken2013, Johannessen2014}  and (truncated) hierarchical B-splines \cite{GIANNELLI2016}. We are mainly interested on methods that preserves the tensor product structure, because its simplicity. Among the others, we mention Forsey and Bartels, who introduced hierarchical B-splines for interpolation and least square approximation of gridded data  \cite{Forsey1995};
Lee et al. \cite{Lee1997} defined a multi-level quasi-interpolant, the Multilevel B-spline Approximation (MBA), where the coefficients depends on values of tensor product B-splines defined on lattices. Similarly, Greiner and Horman \cite{Greiner1997} addressed approximation and interpolation of data by global least squares over hierarchical tensor-product spline spaces.  Kiss et al. \cite{Kiss2014} adapted this approach to handle THB-splines.

\section{Weighted Quasi Interpolant Spline Approximations\label{wQISA}}

Quasi-interpolation is a low-cost and accurate procedure in function approximation theory. The term quasi-interpolation has been interpreted differently according to the context. We follow the Cheney \cite{Cheney1995} definition of a quasi interpolant as any linear operator $L$ of the form
\begin{equation}
\label{equation:quasiinterpolant}
Lf:=\sum_{i=1}^{n_x}\sum_{j=1}^{n_y}f(x_i,y_j)g_{i,j},
\end{equation}
where $f:\Omega\subset\mathbb{R}^2\to\mathbb{R}$ is a function being approximated, $n_x\in\mathbb{N}\cup\{+\infty\}\ni n_y$, $(x_i,y_j)$ are given \emph{nodes} and $g_{i,j}:\Omega\subset\mathbb{R}^2\to\mathbb{R}$ are functions at our disposal. Differently from classic Quasi Interpolation methods that focus to function approximation, the Weighted Quasi Intepolant Spline Approximation (wQISA, \cite{Raffo2019}) aims at point cloud approximation, where the data are assumed to be affected by noise, outliers or partially missing.
In this Section we summarise the wQISA concept and list its theoretical properties. Then, we sketch a data-driven implementation of the wQISA algorithm, which allows to deal with both approximation and prediction problems. Finally, we discuss its computational complexity. 

\subsection{Preliminary concepts\label{overview}}
In a general pipeline for approximation, we envisage that a complex surface is decomposed in multi-charts, e.g., \cite{Ohtake2003,Sorgente2018} and then the approximated patches are stitched together. In this paper we focus on the quality of local approximations, reserving to future investigations the gluing of multiple patches. From the mathematical point of view, every surface can be locally projected onto a plane by using an injective map and, if locally regular, it can be expressed in local coordinates as $(x,y,z(x,y))$. 
Given a degree $p$, a knot vector is said to be $(p+1)$-\emph{regular} if no knot occurs more than $p+1$ times and each boundary knot occurs exactly $p+1$ times.

\begin{definition}
\label{definition:wQISA}
Let $\mathcal{P}\subset\mathbb{R}^3$ be a point cloud and $\mathbf{p}=(p_x,p_y)\in\mathbb{N}^\ast\times\mathbb{N}^\ast$ be a bi-degree. Let $\mathbf{x}$ be a $(p_x+1)$-regular knot vector with boundary knots $x_{p_x}=a_1$ and $x_{n_x}=b_1$ and let $\mathbf{y}$ be a $(p_y+1)$-regular knot vector with boundary knots $y_{p_y}=a_2$ and $x_{n_y}=b_2$. The \emph{Weighted Quasi Interpolant Spline Approximation} of bi-degree $\mathbf{p}$ to the point cloud $\mathcal{P}$ over the knot vectors $\mathbf{x}$ and $\mathbf{y}$ is defined by
\begin{equation}
\label{equation:dkNNVDSA}
f_w(x,y):=\sum_{i=1}^{n_x}\sum_{j=1}^{n_y}\hat{z}_w(x_i^\ast,y_j^\ast)\cdot B[\mathbf{x}_i,\mathbf{y}_j](x,y),
\end{equation}
where $x_i^{\ast}:=(x_i+\ldots+x_{i+p_x})/p_x$ and $y_j^{\ast}:=(y_j+\ldots+y_{j+p_y})/p_y$ are the \emph{knot averages}, the expression
\begin{equation}
    \label{equation:control_points_estimator}
    \hat{z}_w(u,v):=\dfrac{\sum\limits_{(x,y,z)\in\mathcal{P}}z\cdot w(x,y,u,v)}{\sum\limits_{(x,y,z)\in\mathcal{P}}w(x,y,u,v)}
\end{equation}
is the \emph{control points estimator} with weight function $w:\mathbb{R}^2\times\mathbb{R}^2\to[0,+\infty)$ and $B[\mathbf{x}_i,\mathbf{y}_j]$ denotes the tensor product B-spline of bi-degree $\mathbf{p}$ which is uniquely determined by the \emph{local knot vectors} $\mathbf{x}_i=[x_i,\ldots,x_{i+p_x+1}]$ and $\mathbf{y}_j=[y_j,\ldots,y_{j+p_y+1}]$.
\end{definition}

\begin{remark}
    \label{rmk:wQISA_input_data}
    In Definition \ref{definition:wQISA}, a wQISA depends on the following inputs: a point cloud, a tensor mesh (uniquely defined by a bi-degree and two knot vectors) and a weight function.
\end{remark}

\begin{remark}
    \label{rmk:weight_function}
    The weight function $w$ defines a \emph{window} around each point $(u,v)$ and the literature is rich of possible choices, including functions with global and local support. Very popular examples for the function $w$ are:
    \begin{itemize}[noitemsep,wide=0pt, leftmargin=\dimexpr\labelwidth + 2\labelsep\relax]
        \item Indicator of radius $r>0$:
            \begin{equation}
            w(x,y,u,v):=\mathbbm{1}_{||(x-u,y-v)||_2\le{}r},
            \label{equation:indicator_weight}
        \end{equation}
        \item Gaussian of standard deviation $\sigma$:
            \begin{equation}
            w(x,y,u,v):=e^{-||(x-u,y-v)||_2/2\sigma^2},
            \label{equation:gaussian_weight}
        \end{equation}
        \item $k$-Nearest Neighbors ($k$-NN):
            \begin{equation}
            w(x,y,u,v):=
            \begin{cases}
                1/k, 	& \text{ if } (x,y)\in N_k(u,v) \\
                0,	& \text{ otherwise }
            \end{cases},
            \label{equation:kNN_weight}
        \end{equation}
        where $k\in\mathbb{N}^\ast$ and $N_k(u,v)$ denotes the neighborhood of $(u,v)$ defined by the $k$ closest points of the point cloud.
        \item Inverse Distance Weight (IDW):
            
            \begingroup
            \footnotesize
            \noindent
            \thinmuskip=\muexpr\thinmuskip*4/8\relax
            \medmuskip=\muexpr\medmuskip*4/8\relax 
            \begin{equation}
            w(x,y,u,v):=
            \begin{cases}
                \dfrac{1}{||(x,y)-(u,v)||_2}, & \text{\hspace{-.5em}\scriptsize if } |C_{(u,v)}|=0 \\
                \begin{cases}
                    \dfrac{1}{|C_{(u,v)}|}, &   \hspace{-.5em}\forall\, (x,y)=(u,v) \\
                    0,                      &   \text{\hspace{-.5em}\scriptsize else }    
                \end{cases}
            , & \text{\hspace{-.5em}\scriptsize if } |C_{(u,v)}|\ne 0
            \end{cases}
            \label{equation:IDW_weight}
            \end{equation}
            \endgroup
            where 
            \begin{equation*}
            |C_{(u,v)}|:=\{(x,y,z)\in\mathcal{P} \text{ s.t. } (x,y)=(u,v)\}.
            \end{equation*}
    \end{itemize}
    Details on how to fix the free parameters, if any, are provided in Section \ref{parameter_settings}.
\end{remark}

\begin{figure}[!h]
    \begin{center}
        \fbox{
        \includegraphics[scale=0.20, trim={2cm 1cm 2cm 1cm},clip]{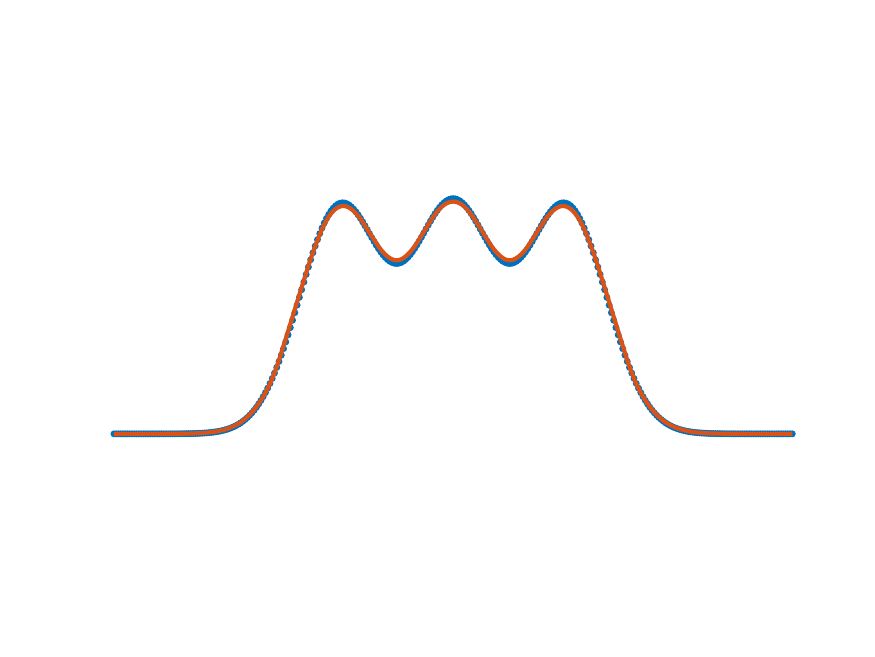}
        }
        \fbox{
        \includegraphics[scale=0.20, trim={2cm 1cm 2cm 1cm},clip]{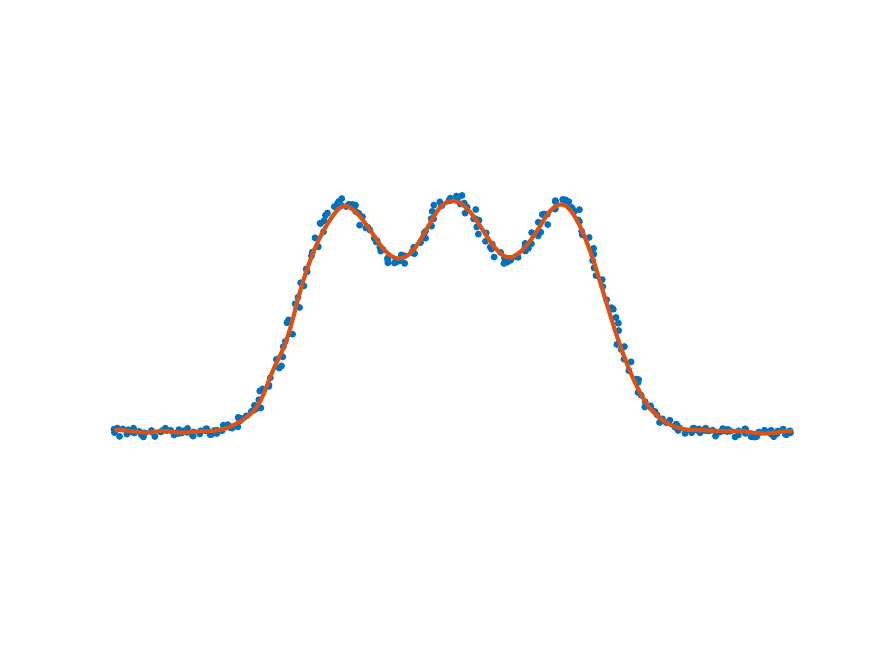}
        }
        \fbox{
        \includegraphics[scale=0.20, trim={2cm 1cm 2cm 1cm},clip]{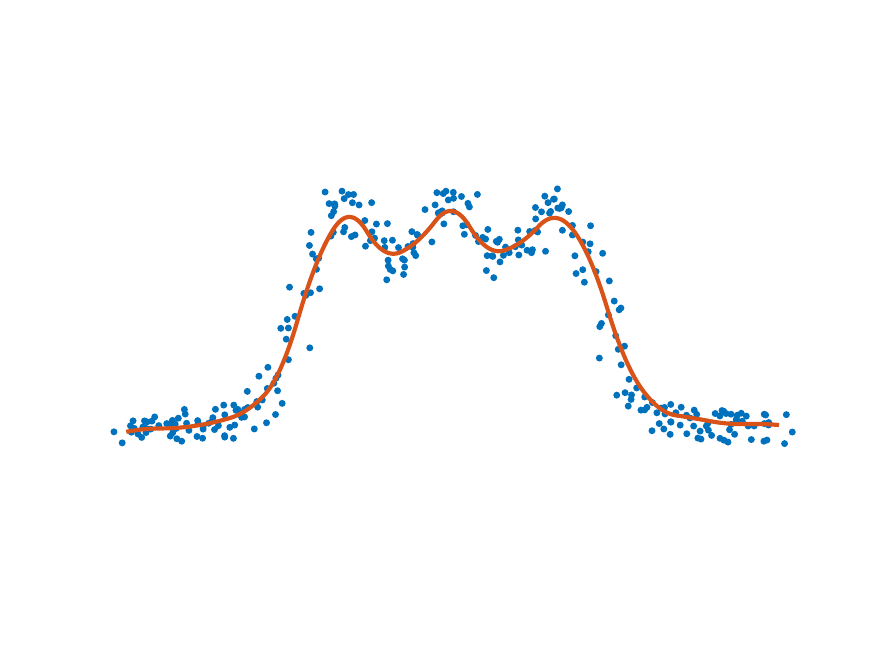}
        } \\
        \fbox{
        \includegraphics[scale=0.20, trim={2cm 1cm 2cm 1cm},clip]{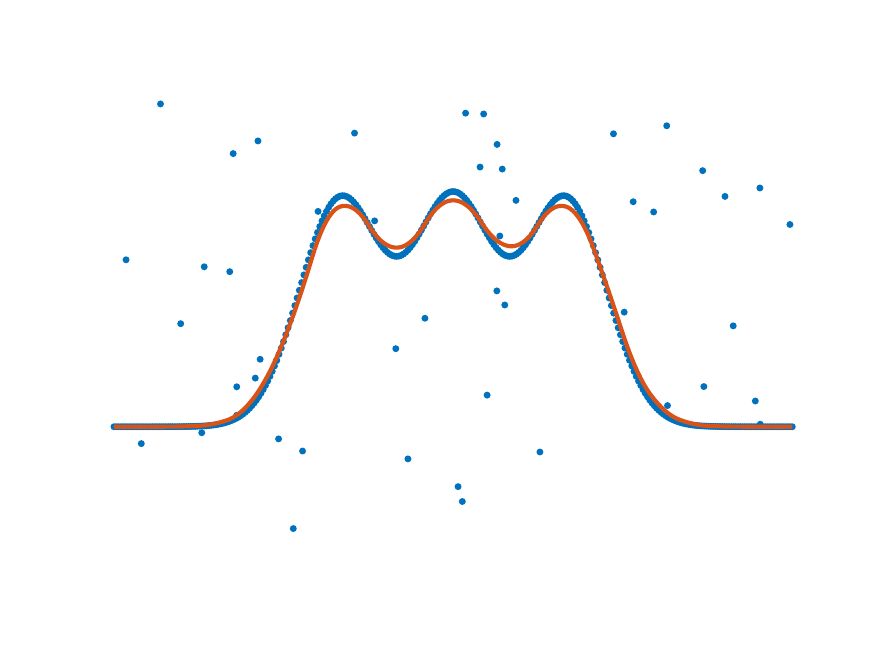}
        }
        \fbox{
        \includegraphics[scale=0.20, trim={2cm 1cm 2cm 1cm},clip]{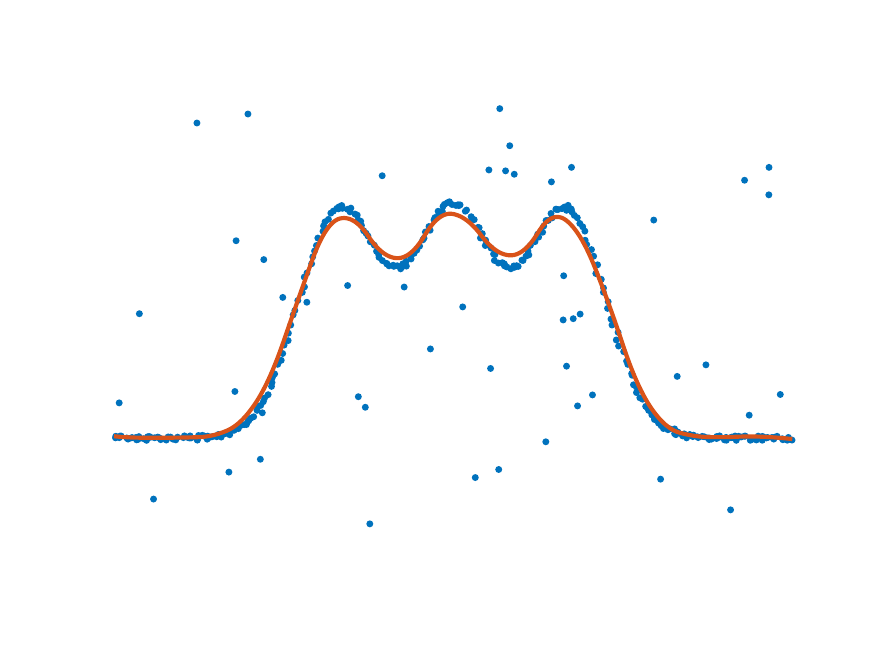}
        }
        \fbox{
        \includegraphics[scale=0.20, trim={2cm 1cm 2cm 1cm},clip]{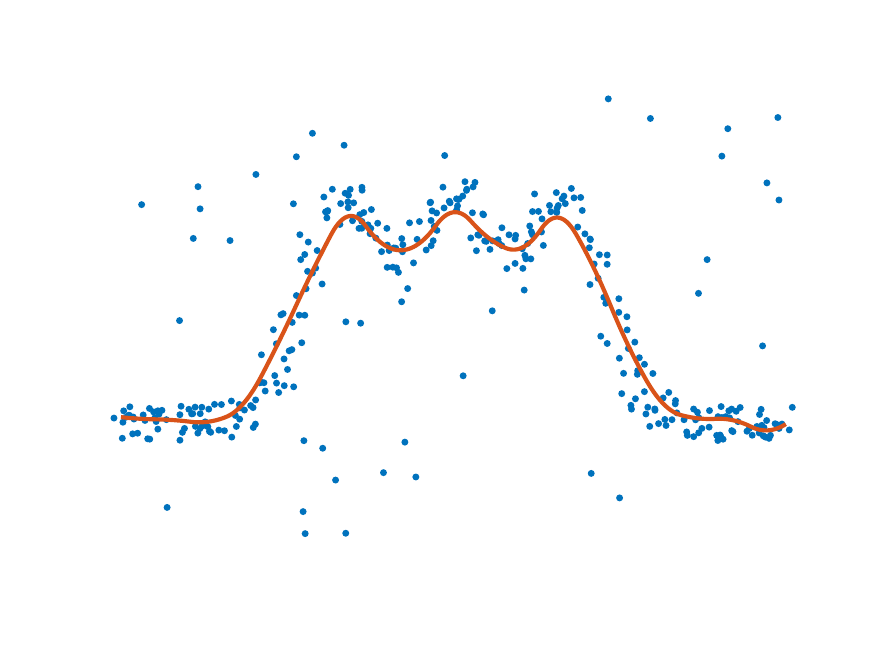}
        }
        \\
        \fbox{
        \includegraphics[scale=0.20, trim={2cm 1cm 2cm 1cm},clip]{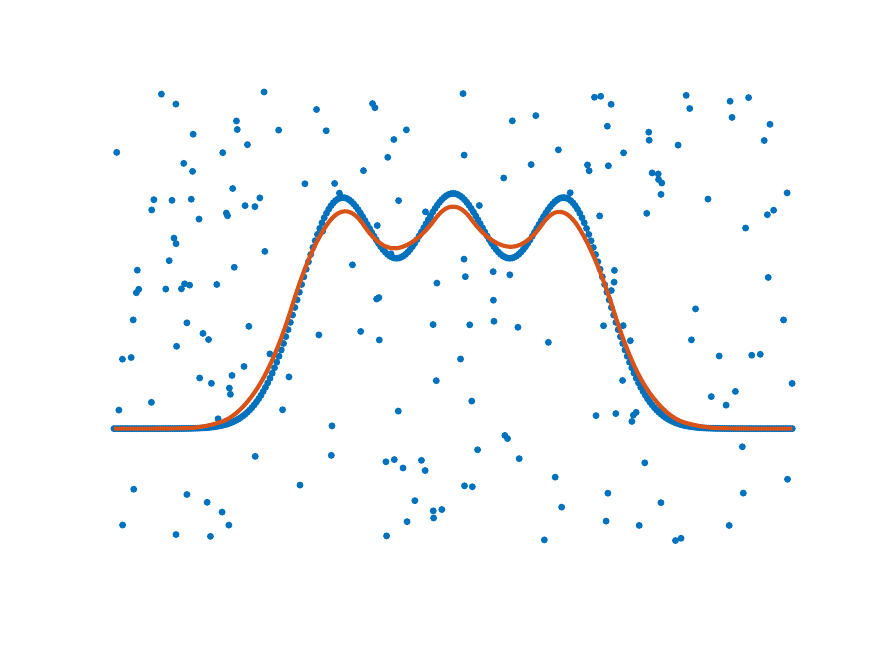}
        }
        \fbox{
        \includegraphics[scale=0.20, trim={2cm 1cm 2cm 1cm},clip]{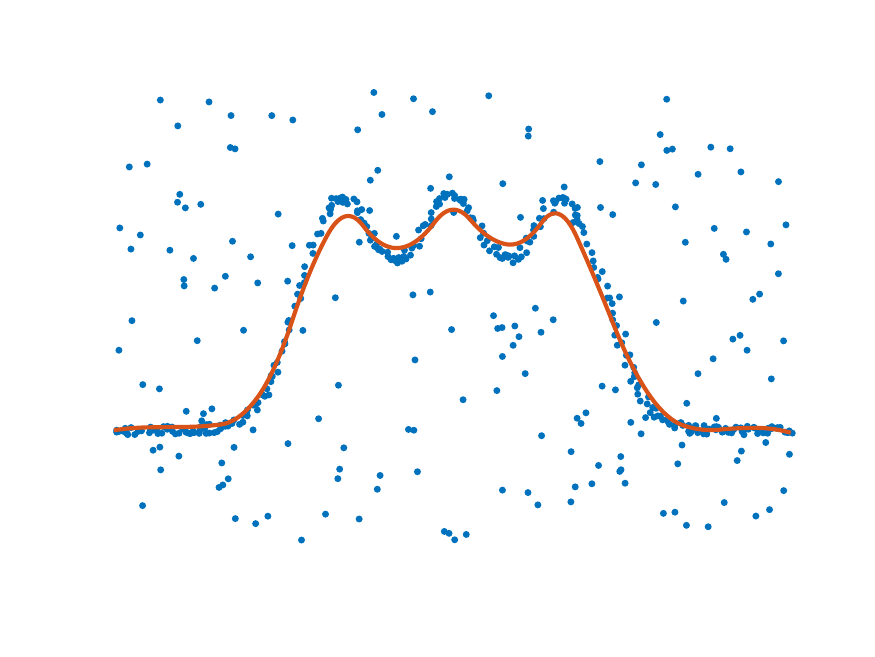}
        }
        \fbox{
        \includegraphics[scale=0.20, trim={2cm 1cm 2cm 1cm},clip]{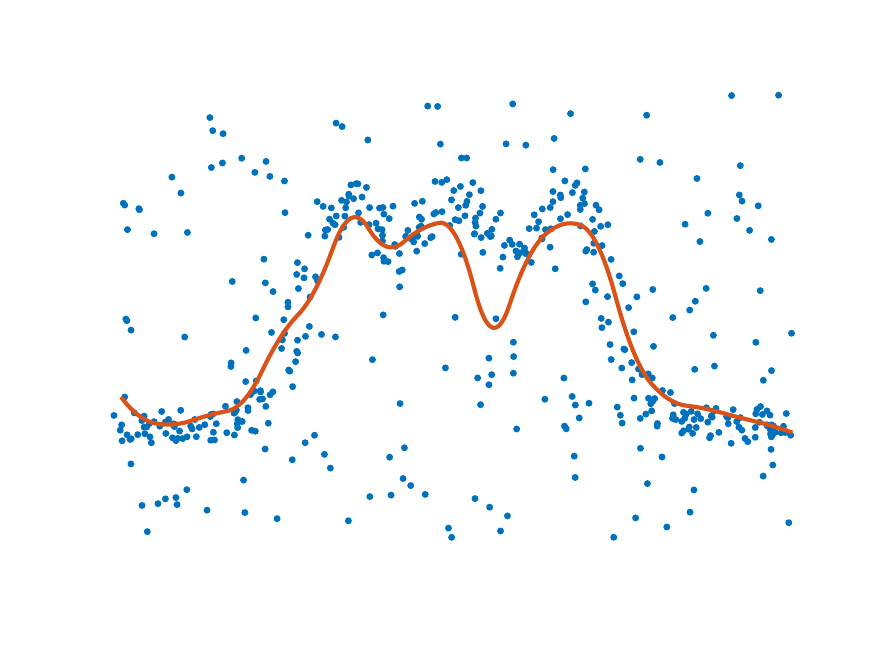}
        }
    \end{center}
    \caption{Noise and outliers robustness. The input point cloud is represented in blue, the output approximation is displayed in red. The noise increases from the left to the right. The outliers increase from top to bottom.}
    \label{figure:noiseandoutliers}
\end{figure}

Figure \ref{figure:noiseandoutliers} shows the stability of wQISA (until failure) against increasing amount of noise and outliers. In these simulations, a $k$-NN weight is filtered from outliers by using quartiles. 

\subsection{Properties}
Similarly to the classical quasi-interpolant schemes \cite{deboor73}, the wQISA method satisfies a number of desirable regularity properties, here briefly recalled. For more details and for a probabilistic interpretation of the method, we refer the reader to  \cite{Raffo2019}.

\paragraph*{Global bounds}
Let $\mathcal{P}\subset\mathbb{R}^3$ be a point cloud and $z_{\text{min}}, z_{\text{max}}\in\mathbb{R}$  that satisfy
\[
z_{\text{min}}\le z\le z_{\text{max}}, \quad \text{ for all } (x,y,z)\in\mathcal{P}.
\]
Then the weighted quasi interpolant spline approximation to $\mathcal{P}$ from some spline space $\mathbb{S}_{p,[\mathbf{x},\mathbf{y}]}$ and some weight function $w$ has the same (global) bounds
\begin{equation}
\label{equation:global_bounds}
z_{min}\le f_w(x,y)\le z_{max}, \quad \text{ for all } (x,y)\in\mathbb{R}^2.
\end{equation}
\paragraph*{Local bounds}
Let $x\in[x_\mu,x_{\mu+1})$ for some $\mu$ in the range $p_x+1\le\mu\le n_x$ and $y\in[y_\nu,y_{\nu+1})$ for some $\nu$ in the range $p_y+1\le\nu\le n_y$. The global bounds of Equation \ref{equation:global_bounds} can be refined, by using the B-splines' property of local representation, to
\begin{equation*}
\min_{\substack{i=\mu-p_x,\dots,\mu \\ j=\nu-p_y,\dots,\nu}}\hat{z}_w(x_i^\ast,y_j^\ast)\le f_w(x,y)\le\max_{\substack{i=\mu-p_x,\dots,\mu \\ j=\nu-p_y,\dots,\nu}}\hat{z}_w(x_i^\ast,y_j^\ast).
\end{equation*}

We can simplify these bounds to:
\begin{equation}
\label{equation:local_bounds}
\min_{(x,y,z)\in\mathcal{P}_{\mu,\nu}}z
\le f_w(x)\le
\max_{(x,y,z)\in\mathcal{P}_{\mu,\nu}}z,
\end{equation}
where
\[
\mathcal{P}_{\mu,\nu}:={\textstyle\bigcup\limits_{\substack{i=\mu-p_x,\ldots,\mu \\ j=\nu-p_y,\ldots,\nu}}}
\left\{\text{supp}\left(w(\cdot,\cdot,x_i^\ast,y_j^\ast)\right)
\right\}\cap\mathcal{P}.
\]

Note that the wQISA bounds depends on the functions $w$. By choosing weight functions having compact support, one can control these bounds while reducing the number of points required for each control point estimation.

\paragraph*{Special configurations, shape preservation and rate of convergence}
From basic spline theory, it is straightforward that if $x_{i+1}=\ldots=x_{i+p_x}<x_{i+p_x+1}$ and $y_{j+1}=\ldots=y_{j+p_y}<y_{j+p_y+1}$ then wQISA interpolates the control point estimate $\hat{z}(x_i^\ast,y_j^\ast)$. By assuming (1) any knot in $\mathbf{x}$ and $\mathbf{y}$ to have maximum multiplicity,
(2) $w$ to be a $1$-NN weight function and
(3) the point cloud to consist of samplings of a continuous function $f:[a_1,b_1]\times[a_2,b_2]\to\mathbb{R}$ at the knot averages, we have that $f_w$ corresponds to the \emph{Variation Diminishing Spline Approximation} (VDSA) of $f$ of bi-degree $\mathbf{p}$ to the point cloud $\mathcal{P}$ over the knot vectors $\mathbf{x}$ and $\mathbf{y}$ (see for example \cite{Lyche2011}). In this perspective, wQISA can be seen as a generalization of VDSA to perturbed data through a wider family of weight functions. 

There is also another point of view: wQISA corresponds to applying VDSA to $\hat{z}_w:\mathbb{R}^2\to\mathbb{R}$.  As pointed out in \cite{Lyche2011}, VDSA preserves certain shape properties - such as monotonicity and convexity - of the function being approximated. In our case, wQISA will preserve the monotonicity and convexity of $\hat{z}_w$, i.e., of the average trend of the point cloud, with respect to the chosen weight function. In case of points clouds with defects, the average trend is indeed more important than the position of a point with respect to the others.
By considering this analogy between VDSA and wQISA, we conclude that the latter is characterized by a linear convergence (see for example \cite{Lyche2011} for the rate of convergence of VDSA). 

\subsection{Data-driven implementation\label{implementation}}
When dealing with approximation, the quality of a method relies on its prediction capability over independent samples, i.e., on data that has not been used to ``train" the model. Whilst for function approximation the accuracy of a method can be quantified by directly comparing the approximation with the original function, in case of point clouds approximation the estimation of the approximation error on new samples is not trivial (for instance, it is often impossible to re-sample the data). Moreover, only a few benchmarks are available. One way to overcome this problem comes from statistical learning, and is here adopted in the form of a learning-based implementation of wQISA. Our implementation consists of four main steps:
\begin{description}
    \item[Step 1] \emph{Data pre-processing}. The input point cloud is used to generate a \emph{training}, a \emph{validation} and a \emph{test} set, each of which has a specific task in the approximation of the input data. Different strategies are proposed, depending on the size of the input (see Section \ref{data_pre_processing}).
    \item[Step 2] \emph{wQISA formulation}. Given a tensor mesh and a weight function, the training set is used to define the control point estimator $\hat{z}_w$. Notice that $\hat{z}_w$ may depend on free parameters to be tuned (see Section \ref{wQISA_computation}).
    \item[Step 3] \emph{Parameter settings}. A prediction error is defined on the validation set, in order to fix the possible free parameters of the weight functions \ref{parameter_settings}).
    \item[Step 4] \emph{Model assessment}. Having obtained an approximation model on the basis of the previous steps, the test set is used for estimating the generalization error on new data (see Section \ref{model_assessment}).
\end{description}
While the first and the last steps are done only once, the second and the third ones are run in a while-loop, by refining the mesh at each iteration (see Section \ref{termination_criteria}). 

\begin{remark}
    \label{rmk:data_overfitting}
    The use of different point clouds for Steps 1-3 in our data-driven implementation is necessary to avoid data overfitting. This is particularly relevant when looking for a model with a good prediction capability, rather than just aiming at minimizing the error on a finite set of points (see, for instance, \cite{Hastie2009, Reitermanova2010}).
\end{remark}

\begin{figure}[h!]
\centering
\includegraphics[scale= 0.5]{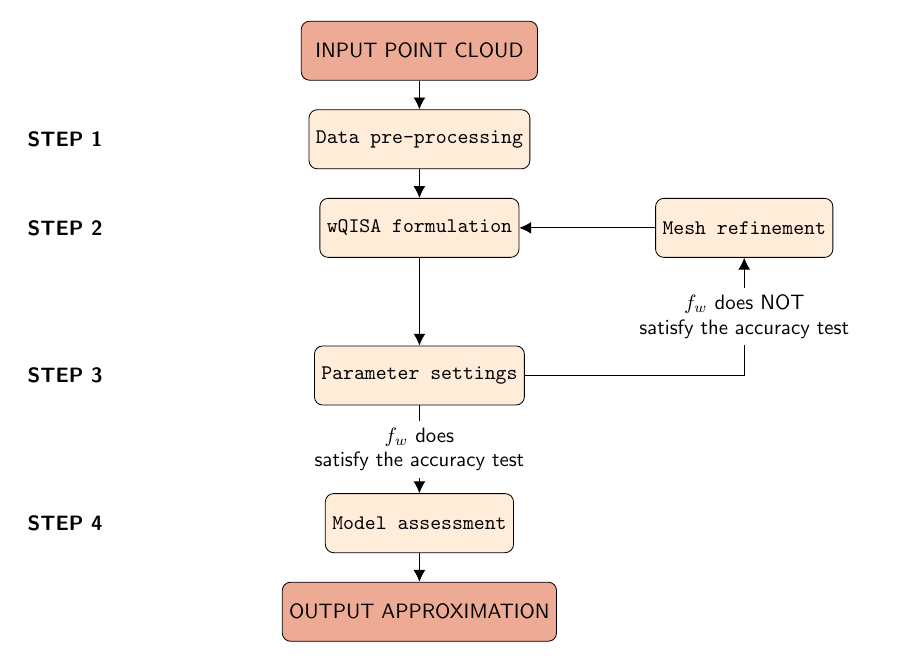}
\caption{Flowchart of the data-driven wQISA algorithm.}
\label{figure:main_algorithm}
\end{figure}

We provide a simplified flowchart of the algorithm in Figure \ref{figure:main_algorithm} and a pseudo-algorithm in Algorithm \ref{algorithm:pseudo_algorithm}, while referring to Sections \ref{data_pre_processing}-\ref{model_assessment} for a detailed description of the implemented procedures.

\subsubsection{Data pre-processing\label{data_pre_processing}}
\begin{figure}[!h]
    \centering
    \begin{tabular}{cccc}
         \includegraphics[height=1.4cm]{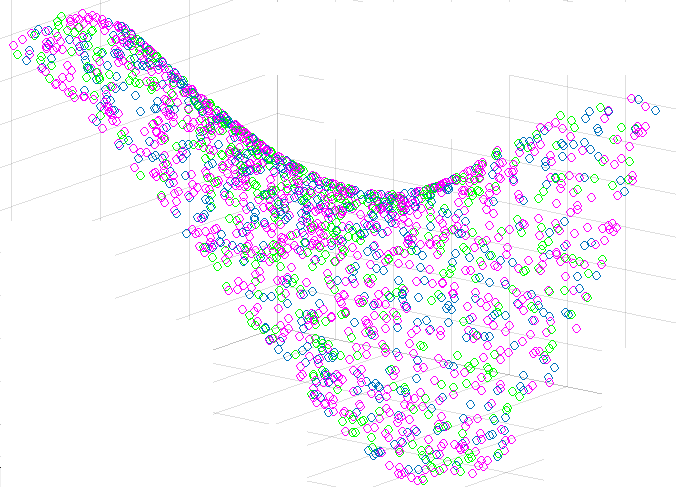} & \includegraphics[height=1.4cm]{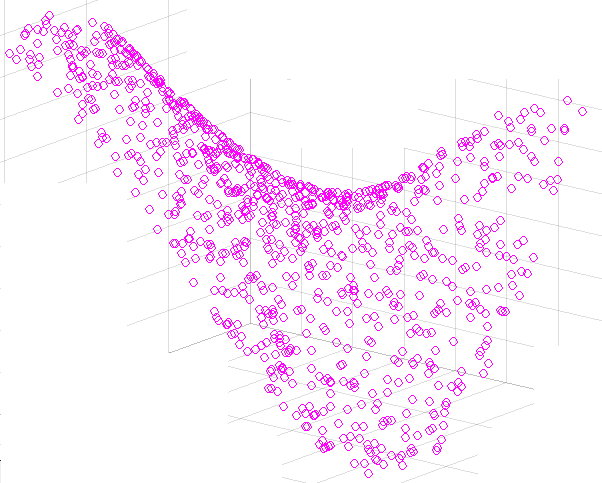} &
         \includegraphics[height=1.4cm]{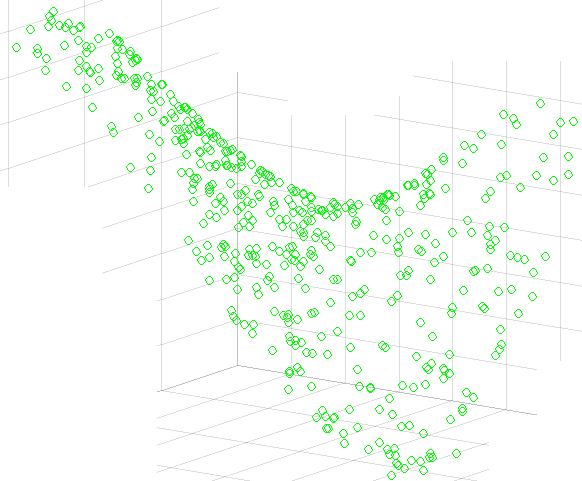} & \includegraphics[height=1.4cm]{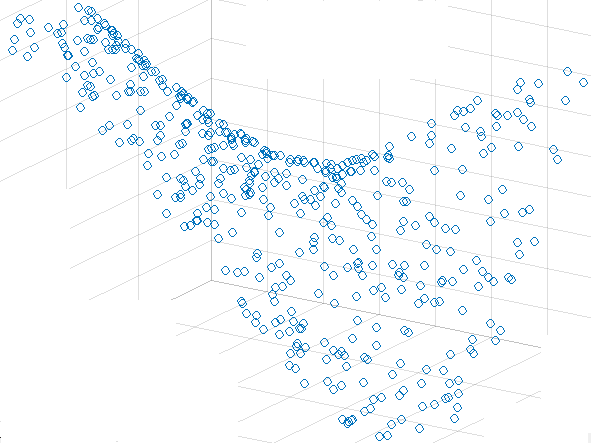} \\
         (a) & (b) & (c) & (d) \\
    \end{tabular}
    \caption{Point cloud partitioning. (a) The original point cloud $\mathcal{P}$; (b) the training set $\mathcal{T}$; (c) the validation set $\mathcal{V}$; (d) the test set $\mathcal{U}$.}
    \label{fig:sampling}
\end{figure}

The input point cloud $\mathcal{P}$ is used to generate three new sets: a training set $\mathcal{T}$, a validation set $\mathcal{V}$ and a test set $\mathcal{U}$. There is no general rule on how to choose the number of observations in each of the three sets.
As suggested in \cite{Hastie2009}, a typical split for large point clouds might be 50\% for training set, and 25\% each for validation and testing sets. The three sets must have similar data distribution properties, i.e., the noise distribution and the point sampling of the validation and test sets must be comparable to the one of the training set. For instance, the training set $\mathcal{T}$ can be easily obtained in \MATLAB{}  from $\mathcal{P}$ by applying the function \texttt{pcdownsample} with $50\%$ of sampling rate and then, $\mathcal{V}$ and $\mathcal{U}$ can be obtained with the sample procedure applied to the complement of $\mathcal{T}$. Figure \ref{fig:sampling} visually shows a possible data partitioning. In order to achieve even more robust results, one can consider several splits of the data into training, validation, and test sets.

{}

In many real-world applications there might be insufficient data to split the set of measurements into three parts. In these cases, we use, the $K$-fold cross-validation: the input point cloud is split into $K$ roughly equal-sized parts and $K-1$ parts are used to train  the model, while the remaining one is used to test it. This process is repeated by changing the part to be used for accuracy evaluation, and then combining the $K$ estimates of the prediction error. In this case, the validation and test sets are the same.

\subsubsection{wQISA formulation\label{wQISA_computation}}
The control point estimator is defined by Equation (\ref{equation:control_points_estimator}). As pointed out in Remark \ref{rmk:wQISA_input_data}, the wQISA input consists of: 
\begin{enumerate}
    \item \emph{A point cloud}. We use the training point cloud computed in Section \ref{data_pre_processing}.
    \item \emph{A tensor product spline space defined by regular knot vectors}. At the first iteration, one can consider a tensor mesh consisting of a single element with maximum knot multiplicities, with this element being the bounding box of the input data. The mesh is then (globally) refined, iteration-by-iteration,  where the approximation shows a lower precision (see Section \ref{termination_criteria}). Although the spline bi-degree can be considered as an additional parameter, we here choose to focus on $C^1$ bi-quadratic spline surfaces as they are smooth enough to represent data in a good way \cite{Skytt2015}.
    \item \emph{A weight function}. We here decide the weight function to be used on a case-by-case basis. One could nevertheless consider dictionaries of weight functions and then choose the more suitable via an opportune validation error.
\end{enumerate}
Having these inputs fixed, the approximation defined in Equation (\ref{equation:dkNNVDSA}) may still depend on additional parameters to be tuned (see Section \ref{parameter_settings}). 

\subsubsection{Parameter settings \label{parameter_settings}}
The control point estimator depends on the parameters defining a weight function (e.g., $k$ in a $k$-NN weight), if any. Statistical learning fixes the free parameters by minimizing some prediction error on an independent data set, here represented by the validation set. The reason for not using the training set for parameter settings is the need to avoid overfitting, i.e., a model that has good performances on the input set but has a low prediction capability on new data. We here fix the parameters by minimizing the \emph{Global Mean Squared Error} (GMSE), i.e., the Mean Squared Error (MSE) over the whole validation point cloud. The minimization of the chosen loss function can be undertaken via iterative methods (e.g., stochastic optimization), or by directly evaluating the loss function at a finite number of parameters (e.g., in case of $k$-NN weight).

\subsubsection{Termination criterion \label{termination_criteria}}
The steps described in the Sections \ref{wQISA_computation} and \ref{parameter_settings} are run in a while loop. The refinement level of a mesh is fixed by considering the last iteration before the GMSE starts increasing (for instance, this occurs when the model starts overfitting the data), with a maximum number of iterations (here we set it to $15$). At each iteration of the while loop, we compute a local validation error and use it, together with a user-defined threshold, to decide which elements should be split by knot insertion. We here insert the knot at the mid-value, in each of the two coordinate directions.

We here consider the \emph{Local Mean Squared Error} (LMSE). Given an element of the mesh, the LMSE on that element is the Mean Squared Error over any validation point whose projection onto the plane $(x,y)$ falls into that element. If no projection lies inside that element, the LMSE for that element is set to zero. 

\subsubsection{Model assessment \label{model_assessment}}
Once the model has been selected, the test set is used to estimate the generalization error. In a data-rich situation, validation and test sets are distinct, and so are in general the validation and test errors. In case the data set was not large enough to be split into three parts, these validation and test sets are equal. This leads the validation and test errors to be the same, and results, in general, in an underestimation of the generalization error (see once more \cite{Hastie2009}).

\begin{algorithm}[t]
    \caption{wQISA pseudocode for $k$-NN weight\label{algorithm:pseudo_algorithm}}
    \begin{algorithmic}[1]
    \State Input: $\mathcal{P}$, $\mathbf{x}$, $\mathbf{y}$, $w$, $\varepsilon$;
    \State Output: $f_w(\mathbf{x}, \mathbf{y})$ (spline approximation), MSE (generalization error)
    \State \Return $\mathcal{T}$, $\mathcal{V}$, $\mathcal{U}$;
    
    \Function{wQISA-Definition}{$\mathcal{T}$, $\mathbf{x}$, $\mathbf{y}$, $w$}
    \State Apply Equation \ref{equation:control_points_estimator} and get $f_w(\mathbf{x}, \mathbf{y},k)$;
    \State \Return $f_w(\mathbf{x}, \mathbf{y},k)$;
    \EndFunction
    \Function{wParameters}{$\mathcal{V}$, $f_w(x,y,k)$}
    \State Find the parameter $k$ that minimizes GMSE;
    \State Get $f_w(x,y)$ by fixing the optimal parameter $k$;
    \State \Return [$f_w(x,y)$,GMSE];
    \EndFunction
    \Function{RefineMesh}{$\mathcal{V}$,$\mathbf{x}$, $\mathbf{y}$,$f_w(\mathbf{x}, \mathbf{y})$,$\varepsilon$}
    \State Compute LMSE from $\mathcal{V}$ and $f_w(\mathbf{x}, \mathbf{y})$;
    \State Check elements where LMSE is above $\varepsilon$;
    \State Split those elements by knot refining both $\mathbf{x}$ and $\mathbf{y}$;
    \State \Return $\mathbf{x}$ and $\mathbf{y}$;
    \EndFunction
    \Function{GeneralizationError}{$\mathcal{U}$,$f_w(\mathbf{x}, \mathbf{y})$}
    \State \Return MSE of $f_w(\mathbf{x}, \mathbf{y})$ over $\mathcal{U}$;
    \EndFunction
    
   \Function{main}{($\mathcal{P}$, $\mathbf{x}$, $\mathbf{y}$, $w$, $\varepsilon$)}
    \NoNumber{\texttt{/*Split the data*/}}
    \State [$\mathcal{T}$, $\mathcal{V}$, $\mathcal{U}$]=DataPre-process($\mathcal{P}$); 
    \NoNumber{\texttt{/*Define the family of approximations*/}}
    \State $f_w(x,y,k)$=wQISA-Definition($\mathcal{T}$, $\mathbf{x}$, $\mathbf{y}$, $w$);
    \NoNumber{\texttt{/*Fix the parameters (here, $k$)*/}}
    \State [$f_w(x,y)$,GMSE]=wParameters({$\mathcal{V}$, $f_w(x,y,k)$});
    \NoNumber{\texttt{/*Repeat while refining the mesh until the stopping criteria is satisfied*/}}
    \While{stopping criteria not satisfied}
    \State $[\mathbf{x},\mathbf{y}]$ = RefineMesh($\mathcal{V}$,$\mathbf{x}$, $\mathbf{y}$,$f_w(\mathbf{x}, \mathbf{y})$,$\varepsilon$);
    \State $f_w(x,y,k)$=wQISA-Definition($\mathcal{T}$, $\mathbf{x}$, $\mathbf{y}$, $w$);
    \State [$f_w(x,y)$,GMSE]=wParameters({$\mathcal{V}$, $f_w(x,y,k)$});
    \EndWhile
    \State MSE=GeneralizationError($\mathcal{U}$,$f_w(\mathbf{x}, \mathbf{y})$);
    \State {\Return  {$[f_w(\mathbf{x}, \mathbf{y}), MSE]$}}
    \EndFunction
    \end{algorithmic}
\end{algorithm}

\subsection{Computational complexity}
The computational complexity of a single wQISA iteration depends on the chosen weight function and tensor product spline space. Given the linear space $\mathbb{S}_{\mathbf{p},[\mathbf{x},\mathbf{y}]}$ of tensor product B-splines of bi-degree $\mathbf{p}$ over the regular knot vectors $\mathbf{x}\in\mathbb{R}^{n_x+p_x+1}$ and $\mathbf{y}\in\mathbb{R}^{n_y+p_y+1}$, there are exactly 
$\text{dim}(\mathbb{S}_{\mathbf{p},[\mathbf{x},\mathbf{y}]})=n_{x}\cdot{}n_{y}$ control points to be estimated.

The Gaussian weight is global and thus computes, for a single control point estimation, the linear combination of $N$ addends. Since the weight of any point is at most as expensive as the exponential of an Euclidean norm, the computational complexity of a control point estimate is $O(N)$. To considerably reduce the computational complexity, one can recover local support by composing global weight functions with a $k$-NN tree (see Equation \ref{equation:kNN_weight}) or with a weight functions with local support, such as the indicator weight function (see Equation \ref{equation:indicator_weight}). 
Indeed, in case of a $k$-NN weight function the time needed to compute all the coefficients through k-d trees is proportional to $O(N\log(N))$, where $N$ is the number of points of the cloud \cite{Friedman:1977}. For its efficiency, k-d trees are already adopted for noise point clouds reconstruction, see for instance \cite{Giraudot:2013}.

\begin{figure}[!h]
\begin{center}
\includegraphics[scale=0.3, trim={1.0cm 7.5cm 1.0cm 7.5cm},clip]{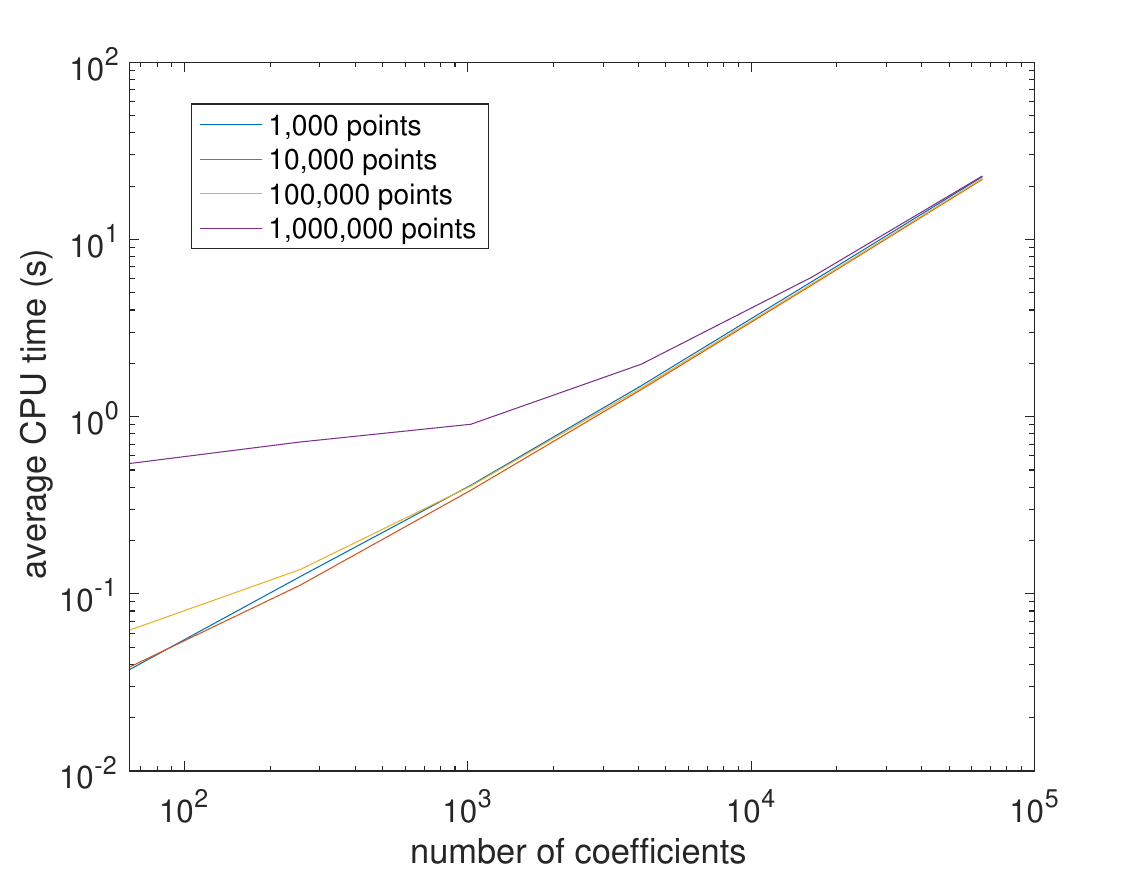}
\caption{CPU times. Log-log plot of the CPU times when increasing the point cloud cardinality and the number of coefficients.}
\label{figure:loglog_times}
\end{center}
\end{figure}

The number of iterations to reach the local minimum of the loss function depends on the input point cloud, the chosen approximation method, and the loss function itself (here: the GMSE over the validation point cloud). In these settings, the use of a method with a high convergence rate can cause an overshooting of the minimum. Experimental results on the number of iterations of our method are reported in Section \ref{simulation} and show that wQISA converges in a reasonably low number of iterations.

n the following, we analyze the CPU times to approximate a point cloud with increasing cardinality, randomly sampled from the mathematical function
\begin{equation*}
    z=\sqrt{64-81\bigl((x-0.5)^2+(y-0.5)^2\bigr)}/(9-0.5)
\end{equation*}
over a mesh with an increasing number of knots. We here approximate an IDW weight function by considering, for each knot average, only the $K$ closest points ($K$ is set to $500$). We have used, for this test, a prototype implementation developed in Python. These tests run on a 2019 MacBook pro with 2.4 GHz 8-cores Intel\textregistered Core\textsuperscript{TM} i9-processor, resulting in the CPU times shown in the log-log plot of Figure \ref{figure:loglog_times}.

\section{Experimental settings}
\label{sec:settings}
Despite the popularity of data approximation, we could not find a benchmark able to address all the use cases we are targeting (surface reconstruction, terrain modelling and spatial data measurement approximation); thus we propose a new comparative analysis. In this section we detail the approximation methods and the performance metrics we adopt in our work.

\subsection{Approximation methods}
We here overview the three popular approximation schemes used as wQISA counterparts throughout this paper: namely, the implicit approximation with Radial Basis Functions (RBF), Kriging and the Multilevel B-spline Approximation (MBA).

\subsubsection{Implicit approximation with radial basis functions} 
This implicit approximation has the form
\begin{equation}
\label{equation:RBF}
    f(x,y):=\sum_{i=1}^Nw_i\phi(||(x-x_i,y-y_i)||_2),
\end{equation}
where the approximating function is represented as a linear combination of \emph{radial basis functions}, each associated with a different center $(x_i,y_i)$ and weighted by an unknown coefficient $w_i$. The weights $w_i$ can be computed by imposing interpolatory constraints of the data, where the trivial null solution is avoided by adding normal constraints. Depending on the properties of $\phi$, RBFs can be locally- or globally-supported. The computational complexity of global and local approximations is respectively $O(N^3)$ and $O(N\log{N})$. We here consider the following kernels:
\begin{subequations}
\begin{align}
    & \phi(r):=e^{(r/\varepsilon)^2}, 
    & (\text{Gaussian})			   
    \label{equation:gaussian_kernel_RBF} 	
    \\
	& \phi(r):=\sqrt{1+(r/\varepsilon)^2}, 			& (\text{Multiquadric})			\label{equation:multiquadric_kernel_RBF}	\\
	&\phi(r):=\dfrac{1}{\sqrt{1+(r/\varepsilon)^2}}   
	& (\text{Inverse multiquadric})		\label{equation:inverse_multiquadric_kernel_RBF} 
	\\
	&\phi(r):=e^{-\sqrt{r}},
	& (\text{Modified Gaussian})		\label{equation:modified_gaussian_kernel_RBF}
\end{align}
\end{subequations}
where $\varepsilon$ is a shape parameter that here approximates the average distance between nodes. An additional parameter $\alpha\ge0$ is introduced by performing an $L^2$-regularization on the interpolation matrix, in order to increase the smoothness of the approximation ($\alpha=0$ is for interpolation). We consider the implementation from the Python class \texttt{interpolate.Rbf}, contained in the \texttt{Scipy} library.

\subsubsection{Kriging}
Kriging owes its importance to its capability to take into account the correlation among input data, which may strongly affect the approximation, e.g., when unevenly distributed. Kriging is defined by the weighted average
\begin{equation}
    \label{equation:kriging}
    f(x,y):=\mathbf{w}^t\mathbf{z}=\sum\limits_{i=1}^Nw_iz_i,
\end{equation}
where the weights $\mathbf{w}\in\mathbb{R}^N$ are the solution to the linear system $\mathbf{C}\mathbf{w}=\mathbf{c}$, where $\mathbf{C}$ is the covariance matrix of the input point cloud and $\mathbf{c}$ is the vector of covariances between the input points and the prediction point $(x,y)\in\mathbb{R}^2$. Despite its popularity, Kriging suffers several issues when applied to large data sets. For each sample, the ordinary kriging estimator needs to solve a linear system and thus the computational complexity scales quadratically with the number of input points. We use here the Python package \texttt{PyKrige}.

\subsubsection{Multilevel B-spline approximation}
Multilevel B-spline Approximation (MBA) was originally introduced in \cite{Lee1997} for approximating scattered data by tensor product B-spline surfaces. A peculiarity of MBA, that makes it analogue to wQISA, is its explicit formulation, which allows it to compute an approximation without solving any equation system. In this respect, MBA can be considered as an example of quasi-interpolation method. At the first step, the tensor product mesh consists of only one element. At each successive iteration, each element is halved in the two coordinate directions. The coefficient $c_k$ of a B-spline $B_k$ is given by
\begin{equation*}
    c_k=\dfrac{\sum_iB_k(x_{i,k},y_{i,k})^2\phi_i}{\sum_iB_k(x_{i,k},y_{i,k})^2},
\end{equation*}
where ${(x_{i,k},y_{i,k})}_i$ is the set of points within the support of the B-spline $B_k$ and 
\begin{equation*}
    \phi_i:=\dfrac{B_k(x_{i,k},y_{i,k})z_i}{\sum_lB_l(x_{i,k},y_{i,k})^2},
\end{equation*}
where the sum in the denominator is taken over all B-splines which contains $(x_{i,k},y_{i,k})$ in their support. The computational complexity scales linearly with the number of input points. We here use the implementation of the Geometry Group at SINTEF ICT, available at \url{https://github.com/orochi663/MBA}. 

\subsection{Evaluation measures\label{evaluation_measures}}
To evaluate the quality of an approximation against the input point cloud, we need to define some error measures. We here present a comparison with a number of measures, each one able to highlight different approximation aspects. 

\subsubsection{Punctual error and its statistics}
Given an approximation defined on the training point cloud, the \emph{absolute punctual error} is computed on the validation point cloud. The absolute punctual error is studied via measures of central tendency (mean) and of statistical dispersion (standard deviation and mean square error). 

\subsubsection{Hausdorff distance and $L^\infty$ norm}
The \emph{Hausdorff distance} (or \emph{two-sided Hausdorff distance}) $d_{\text{Haus}}$ between two non-empty subsets $A$ and $B$ in $(\mathbb{R}^d,||\cdot||_2)$ is given by
\begin{equation*}
    d_{\text{Haus}}(A,B):=\max
    \left\{
    \sup_{a\in A}\inf_{b\in B} ||a-b||_2,
    \sup_{b\in B}\inf_{a\in A} ||a-b||_2
    \right\},
\end{equation*}
where $\sup$ represents the supremum and $\inf$ the infimum. In our case, $A$ is chosen to be the validation point cloud, while $B$ is the image of an approximation $f$ over its domain. In our settings, the $\sup$ and $\inf$ can be replaced by $\max$ and $\min$ respectively. We will use this distance for all sets of points embedded in an Euclidean space. For data measurements that associate a scalar value to a point of a grid we will use instead the $L^\infty$ norm.

\section{Experimental simulations\label{simulation}}
In our experiments we consider real data coming from different use cases.
We classify the data as affected by different levels of noise: low (e.g., 3D point clouds from high-quality laser scans), average (e.g., terrains from Lidar or sonar acquisition) and high (e.g., air pollution and rainfall measurements from sensors and radars). The data size varies from few dozens to hundreds of thousands of points. All data used in our experiments are provided as additional material. To provide a fair comparative analysis - and avoid any method to overfit - a data-driven implementation is used for all the considered techniques. 
Hence, the approximation parameters and the model validation are kept distinct and performed on the same training and validation datasets for all methods. In addition, the optimal shaping parameter $\alpha$ for the RBF approximations is found through a constrained trust-region method, while for Kriging we always use the default package settings (linear variogram model).

In all Tables, we use the following naming convention: Gaussian kernel (RBF\textsuperscript{1}); multiquadric kernel (RBF\textsuperscript{2}); inverse multiquadric kernel (RBF\textsuperscript{3}); modified Gaussian kernel (RBF\textsuperscript{4}).

\subsection{Surface reconstruction}
\label{sec:3Dmodels}
\begin{figure}[!h]
    \begin{center}
        \begin{tabular}{|cc|}
            \hline
            \includegraphics[width=2.8cm]{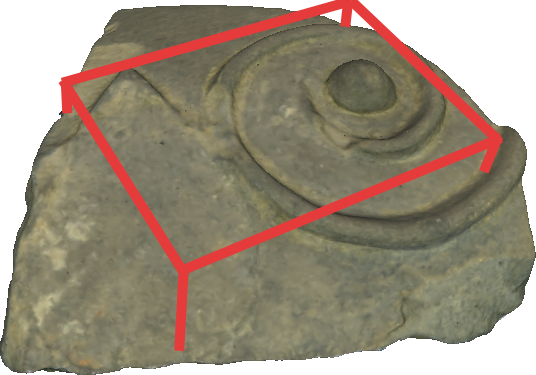} &
            \includegraphics[width=3.3cm,height=4.6cm,keepaspectratio,trim={14cm 4.4cm 13cm 5.3cm}, clip]{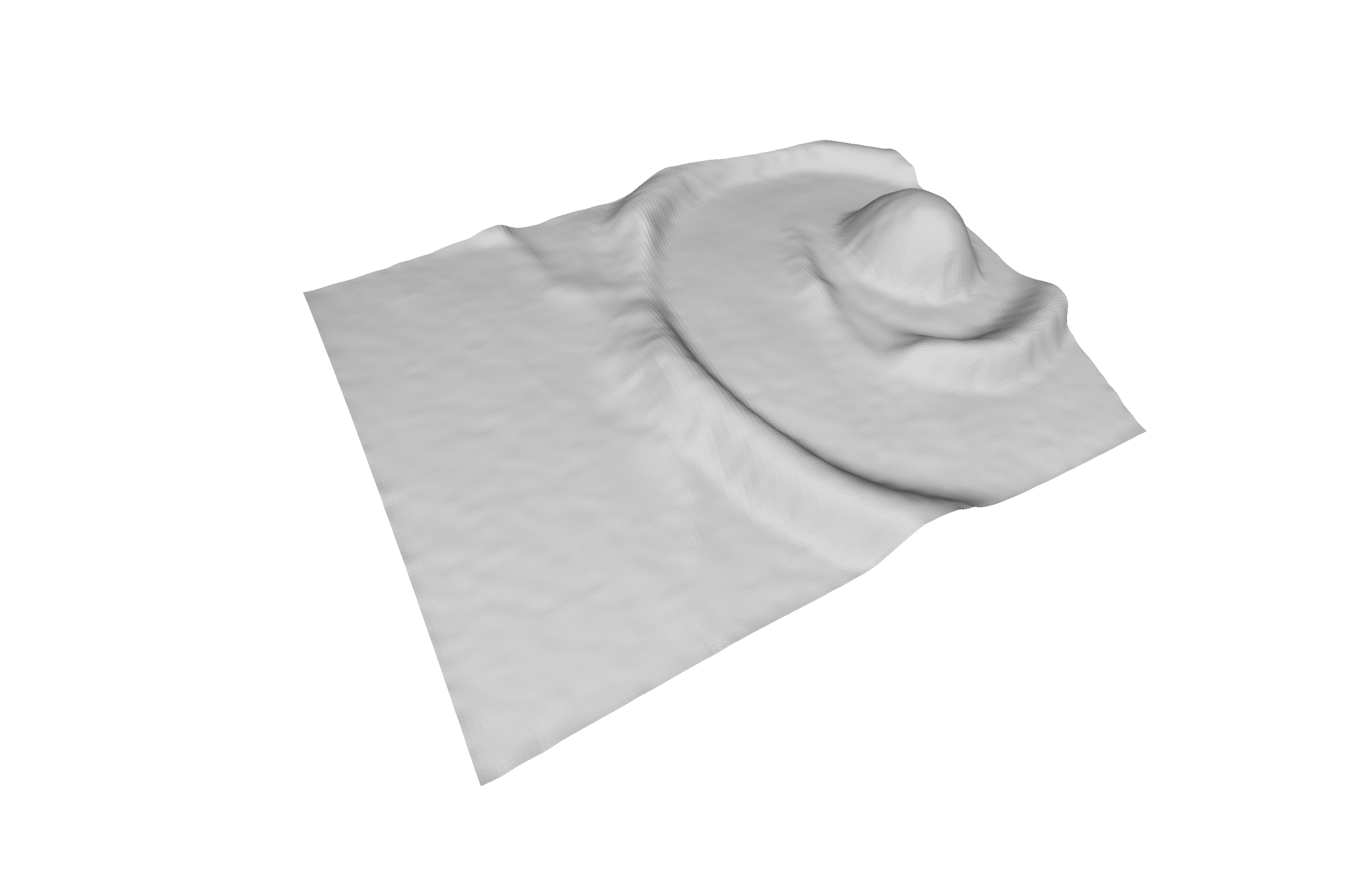} \\
            \multicolumn{2}{|c|}{A.1} \\
            \hline
            \includegraphics[width=2.6cm]{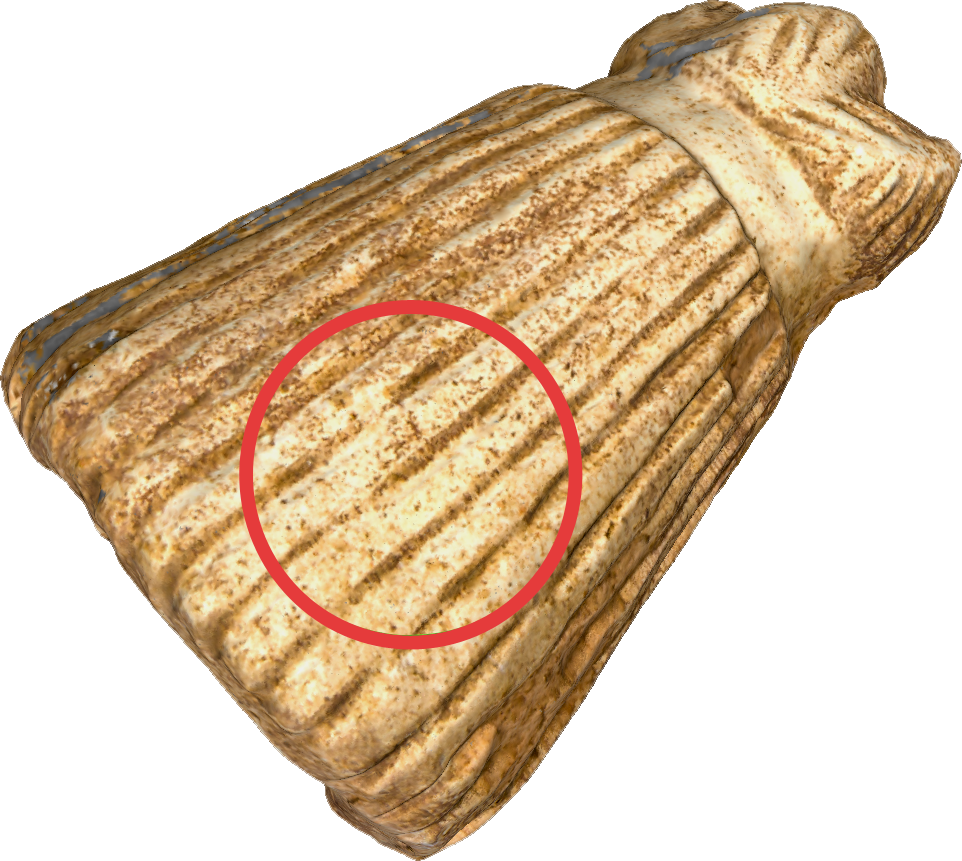} &
            \includegraphics[width=3.3cm,height=4.6cm,keepaspectratio,trim={12.5cm 4.4cm 12.5cm 4.3cm}, clip]{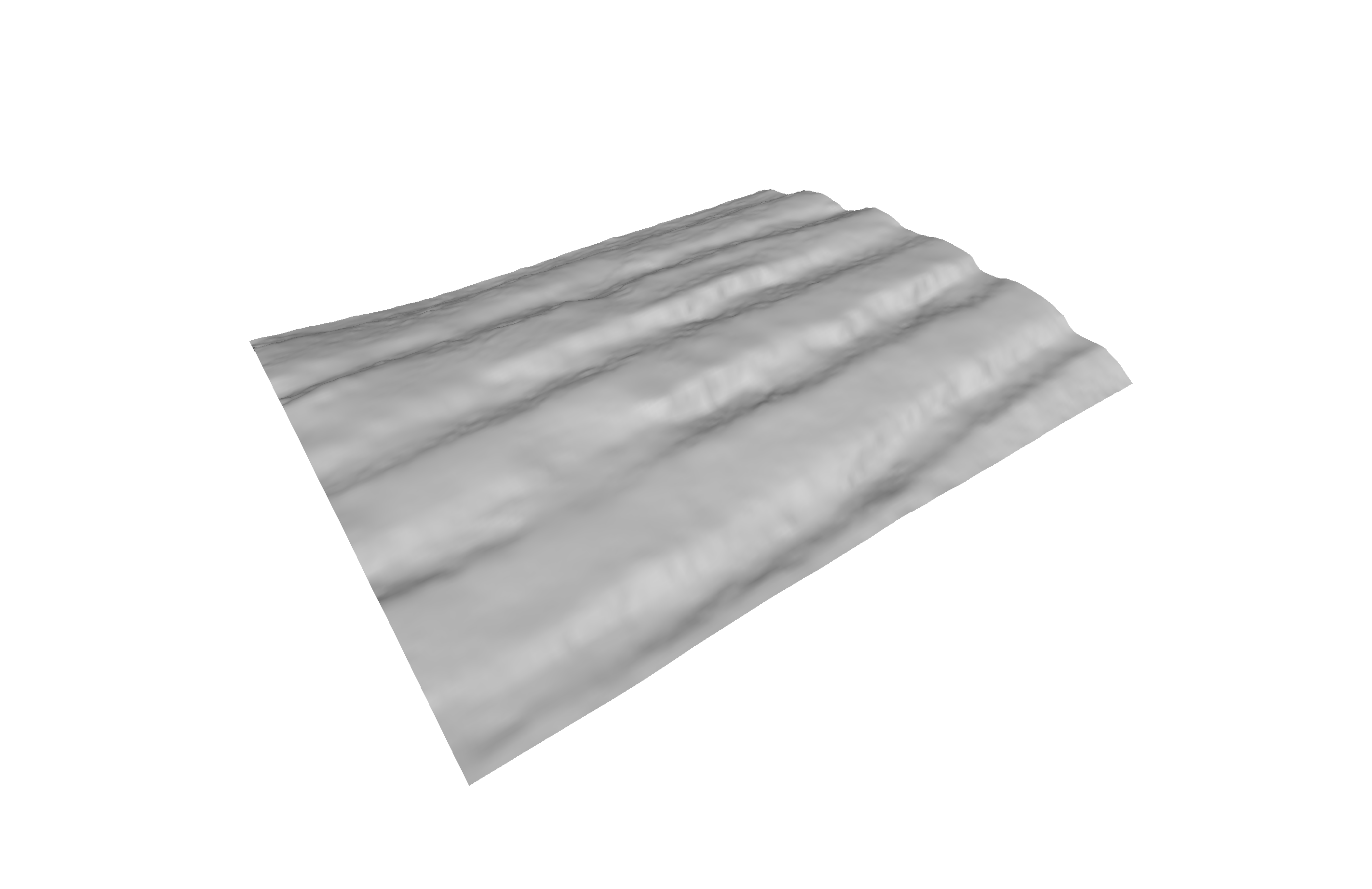} \\
            \multicolumn{2}{|c|}{A.2} \\
            \hline
            \includegraphics[width=2.8cm]{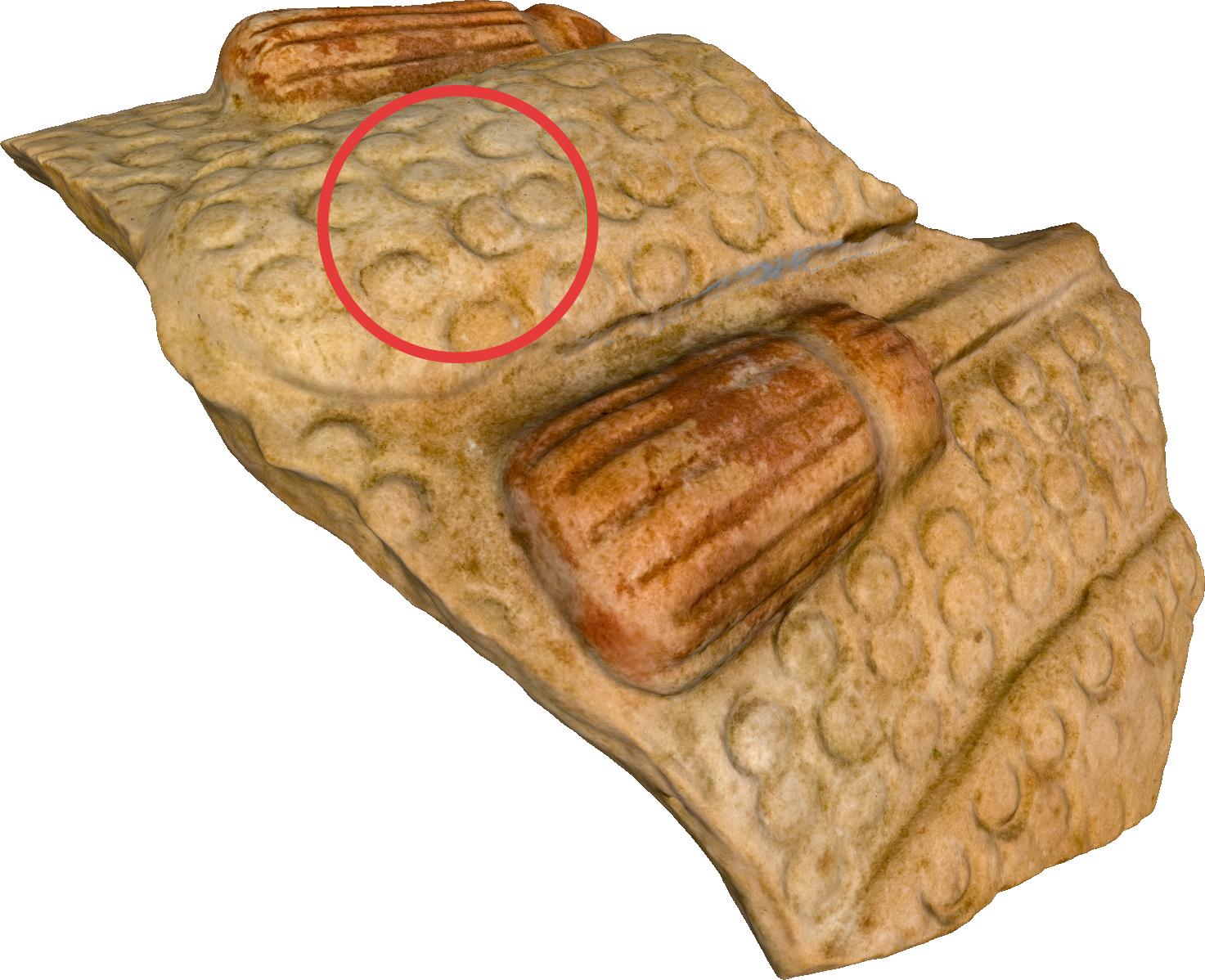} 
             &  
             \includegraphics[scale=0.040,trim={12.5cm 4.1cm 12.5cm 4.1cm}, clip]{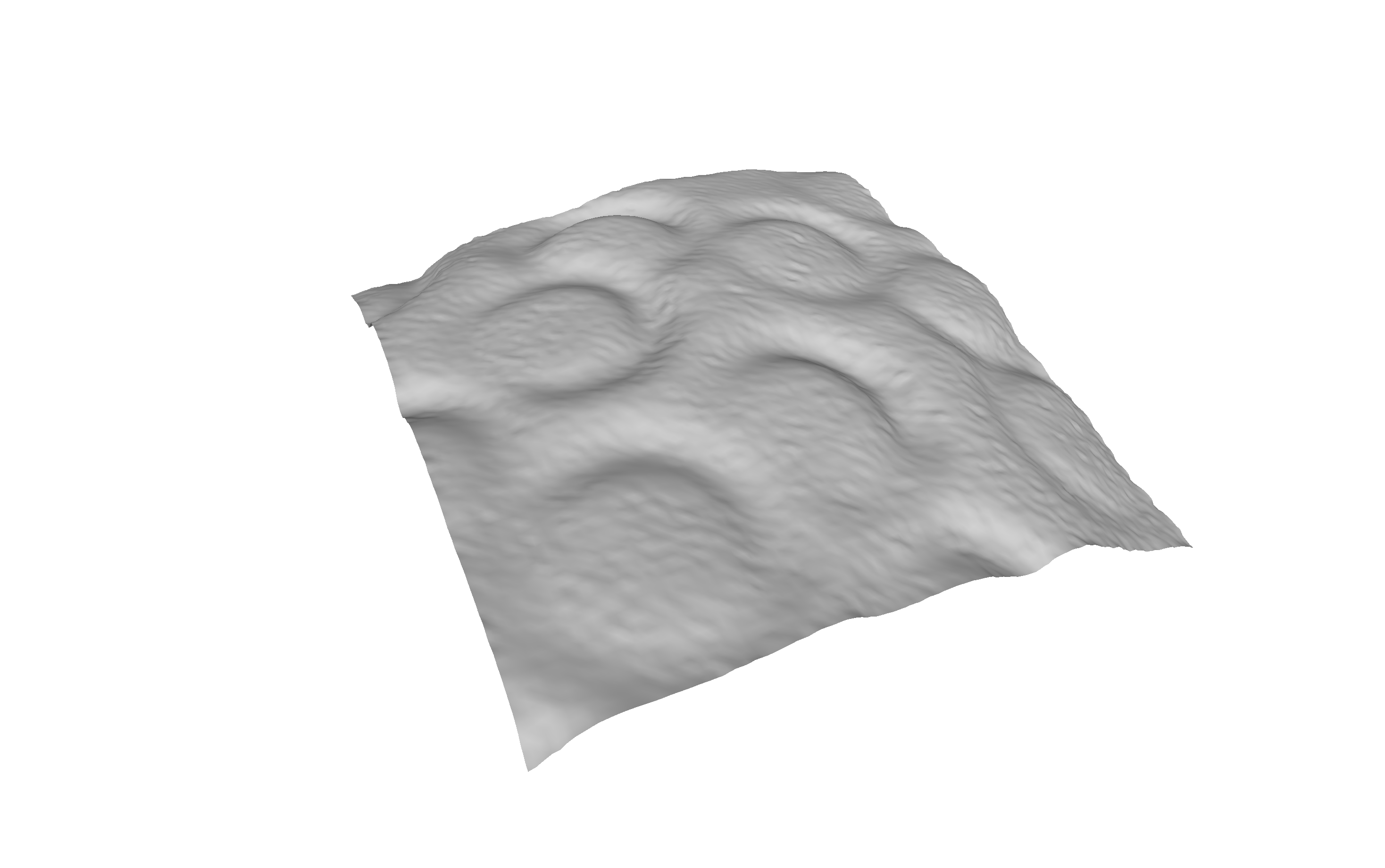}
             \\
             \multicolumn{2}{|c|}{A.3} \\
             \hline
        \end{tabular}
    \caption{Left sides: the original models (from STARC repository \cite{STARC}). Right sides: outcomes of the wQISA method via Gaussian weights for the areas highlighted in red. \label{figure:3Dmodels}}
    \end{center}
\end{figure}

As examples of scattered 3D point clouds, we consider three models from the Science and Technology in Archaeology Research Center (STARC, \cite{STARC}) of The Cyprus Institute, which consist of a digitisation of archaeological fragments. The complete models are shown in Figure \ref{figure:3Dmodels} (left column): the model labelled as A.1 is part of a vessel; models A.2 and A.3 are parts of votive statues. The models are acquired via high-quality laser scans. For our experiments, we select a region of interest in every model (highlighted in red), which is then up-sampled via a uniform Montecarlo approach; namely, A.1 is sampled with $30,000$ points, A.2 with $47,642$ points and A.3 with $29,843$ points. 
We here test wQISA with the Gaussian weight functions of Equation \ref{equation:gaussian_weight}; the standard deviation is $\sigma=10^{-4}$ for the model A.1 and $\sigma=10^{-3}$ for A.2 and A.3.

\begin{figure}[!h]
\begin{center}
\tcbox{
\begin{tikzpicture}[scale=0.80]
\node[inner sep=0pt] (b11) at (0,3.5)
    {\includegraphics[scale=0.0375,trim={12.5cm 4cm 12.5cm 4cm}, clip]{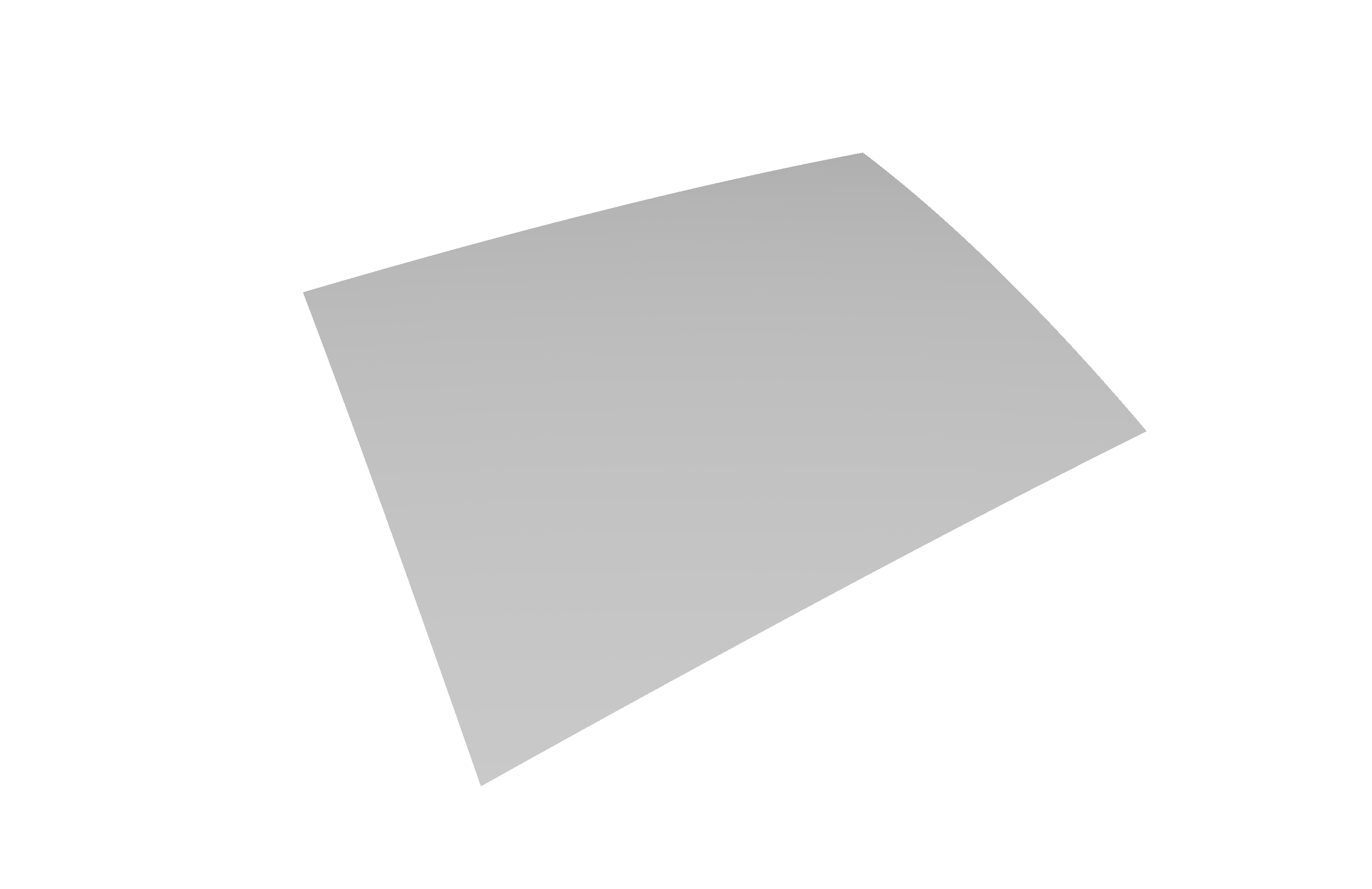}};
    
\node[inner sep=0pt] (b12) at (4.5,3.5)
    {\includegraphics[scale=0.0375,trim={12.5cm 4cm 12.5cm 4cm}, clip]{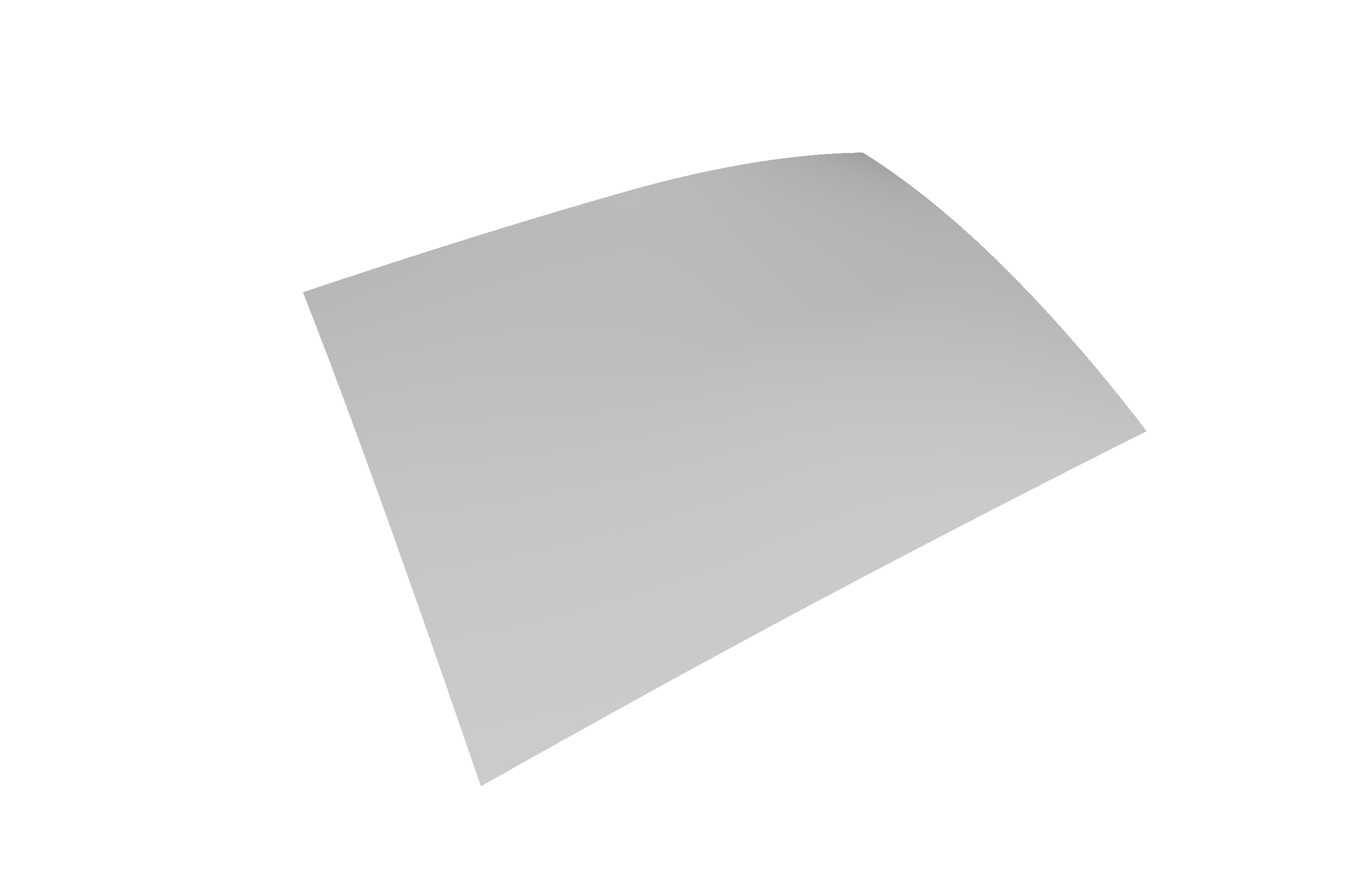}};

\node[inner sep=0pt] (b13) at (9,3.5)
    {\includegraphics[scale=0.0375,trim={12.5cm 4cm 12.5cm 4cm}, clip]{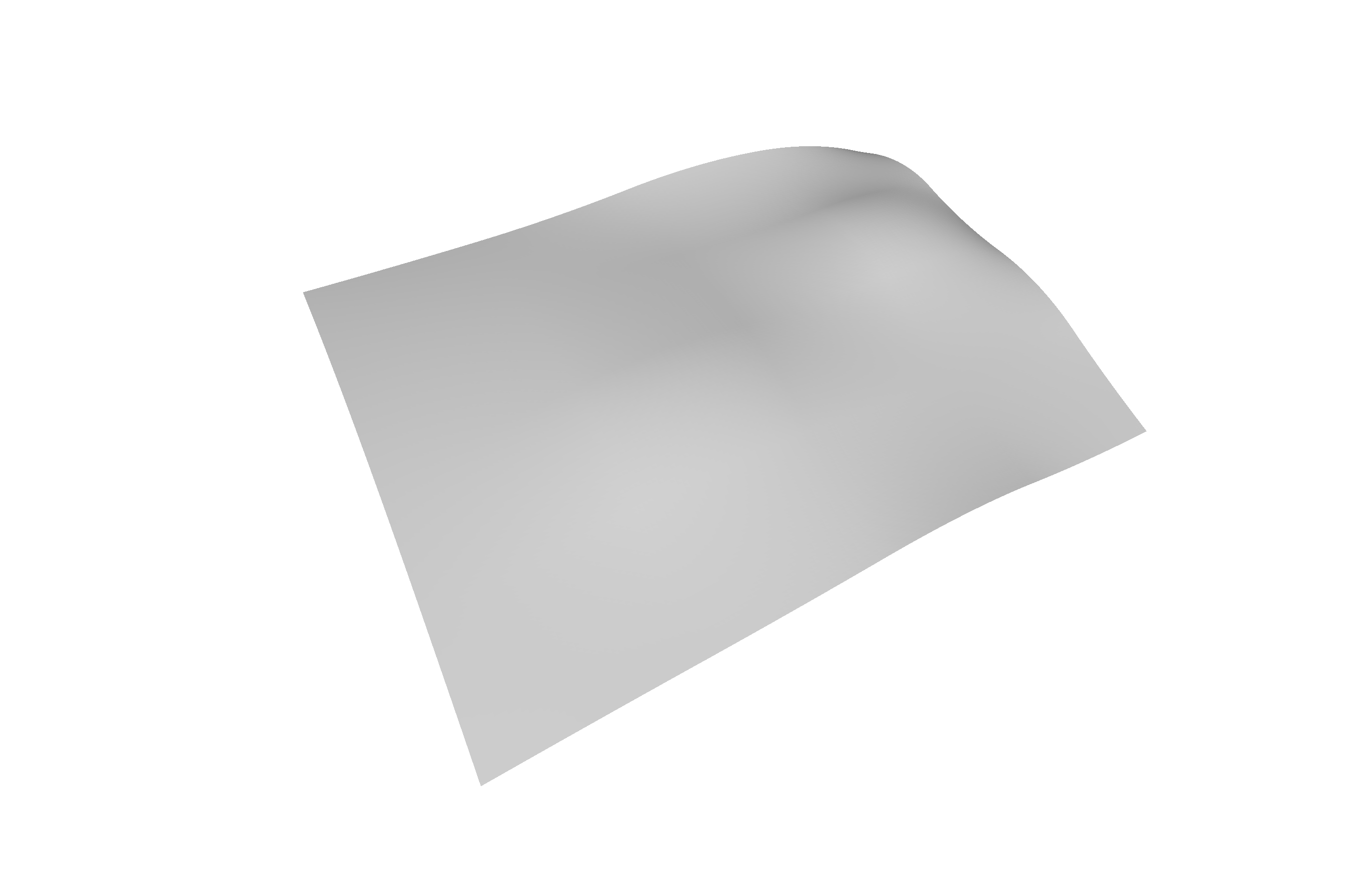}};

\node[inner sep=0pt] (b14) at (13.5,3.5)
    {\includegraphics[scale=0.0375,trim={12.5cm 4cm 12.5cm 4cm}, clip]{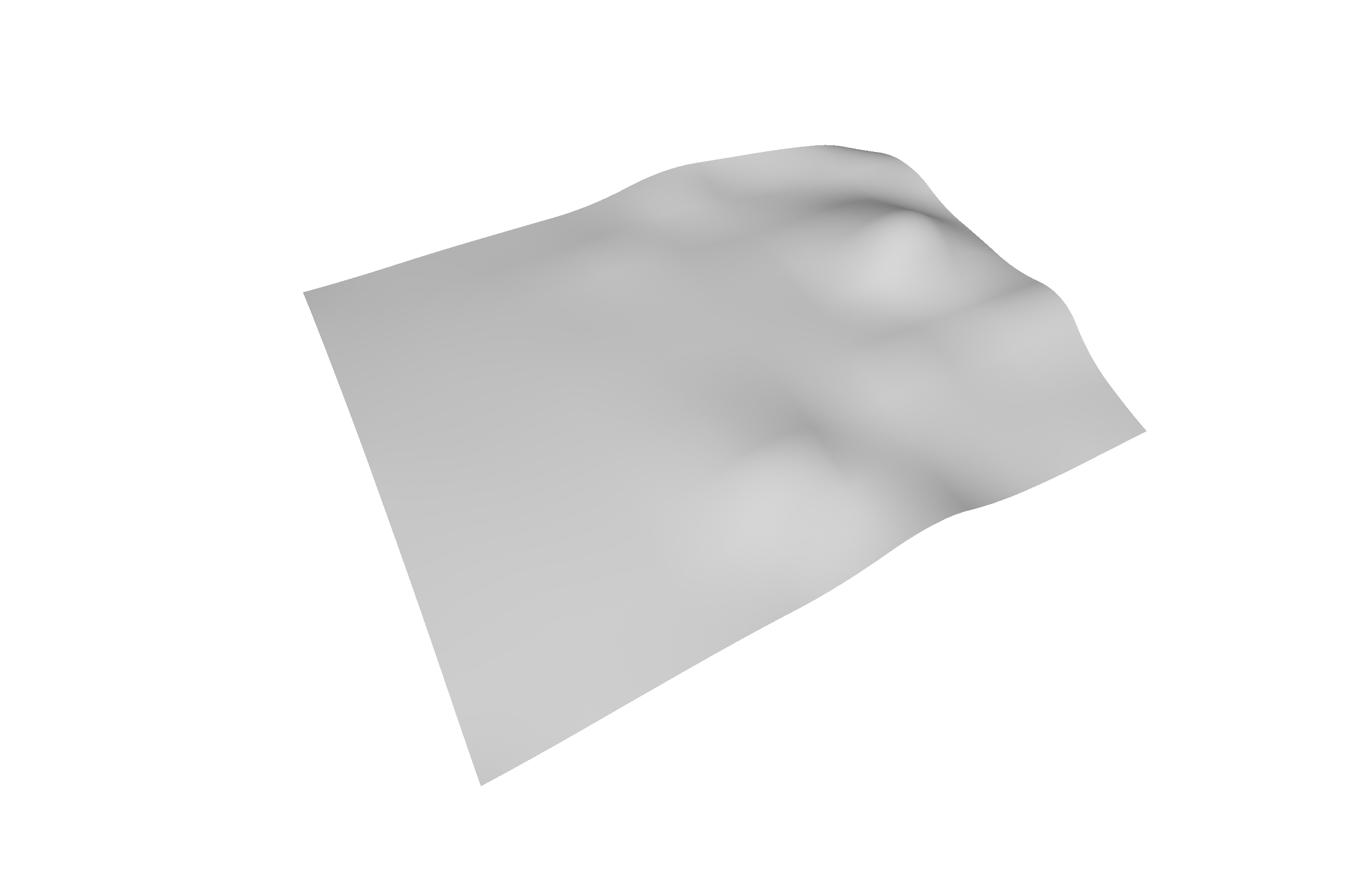}};
    
\node[inner sep=0pt] (b15) at (13.5,0)
    {\includegraphics[scale=0.0375,trim={12.5cm 4cm 12.5cm 4cm}, clip]{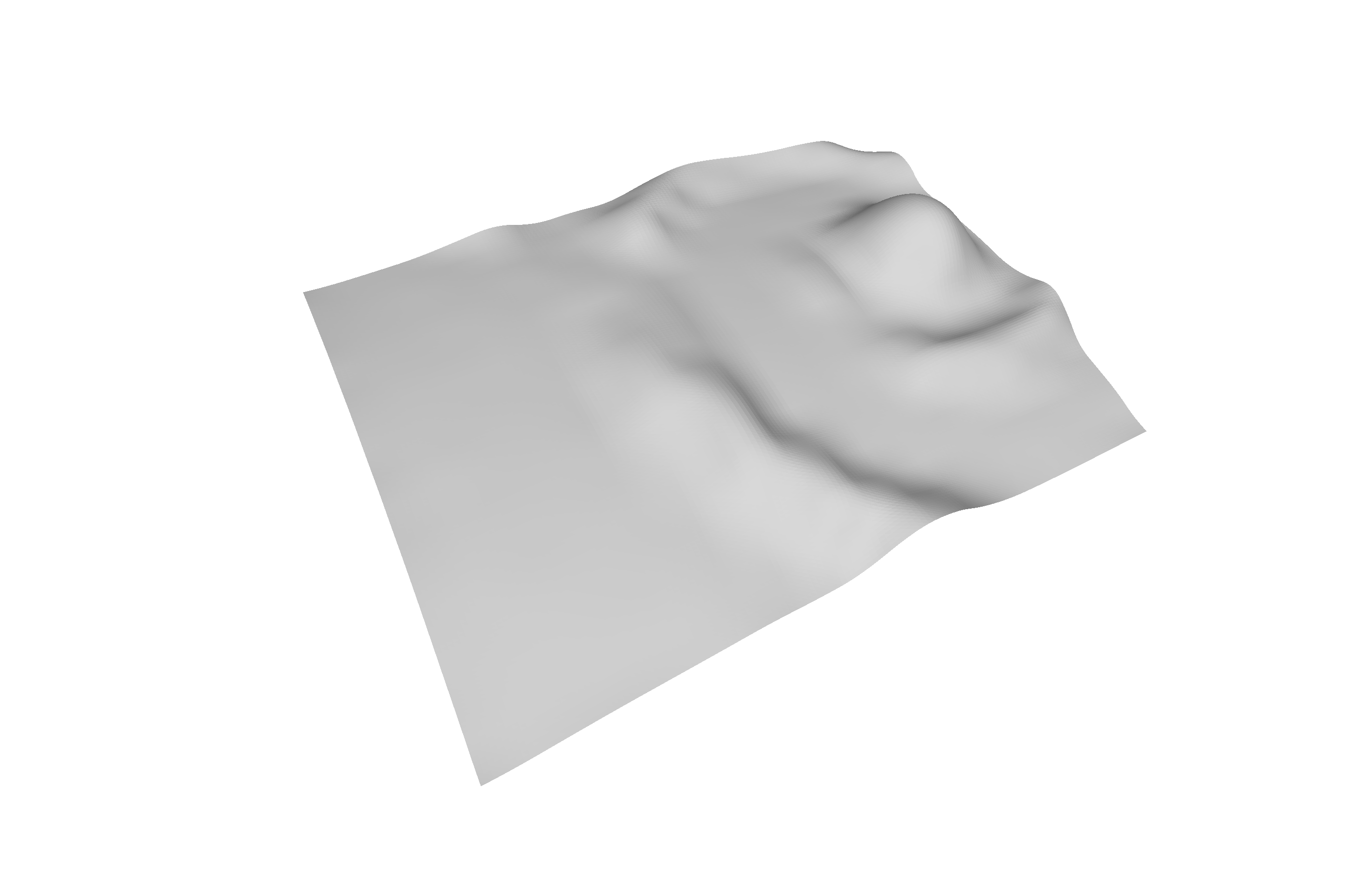}};
    
\node[inner sep=0pt] (b16) at (9,0)
    {\includegraphics[scale=0.0375,trim={12.5cm 4cm 12.5cm 4cm}, clip]{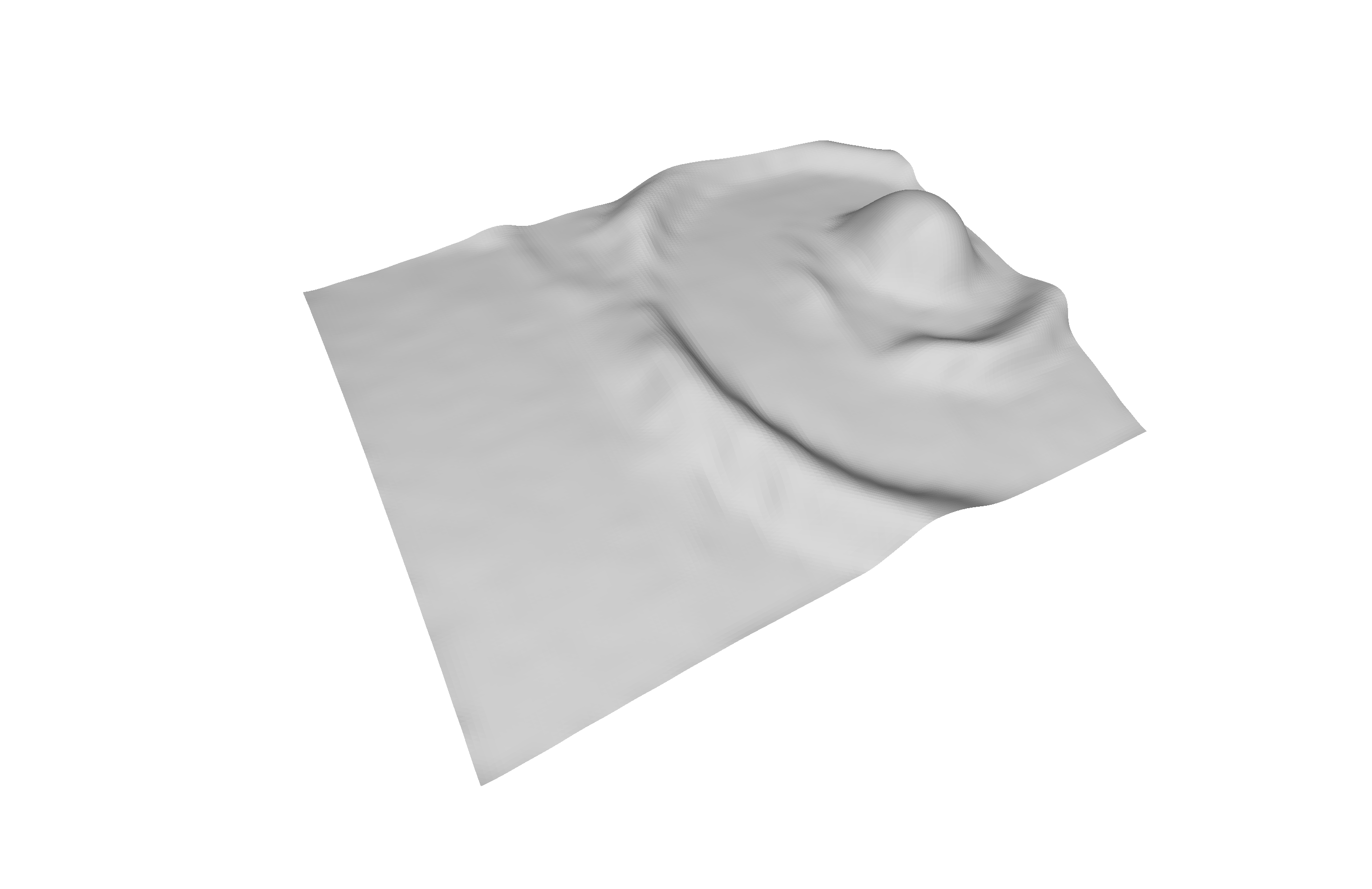}};

\node[inner sep=0pt] (b17) at (4.5,0)
    {\includegraphics[scale=0.0375,trim={12.5cm 4cm 12.5cm 4cm}, clip]{spiral00_L06.png}};
    
\draw [->,thick] (b11) -- (b12);
\draw [->,thick] (b12) -- (b13);
\draw [->,thick] (b13) -- (b14);
\draw [->,thick] (b14) -- (b15);
\draw [->,thick] (b15) -- (b16);
\draw [->,thick] (b16) -- (b17);

\end{tikzpicture}
}
\textbf{A.1} The MSE over the test set is $1.47\cdot{}10^{-6}$.
\tcbox{
\begin{tikzpicture}[scale=0.80]
\node[inner sep=0pt] (c11) at (0,3.5)
    {\includegraphics[scale=0.0375,trim={12.5cm 4cm 12.5cm 4cm}, clip]{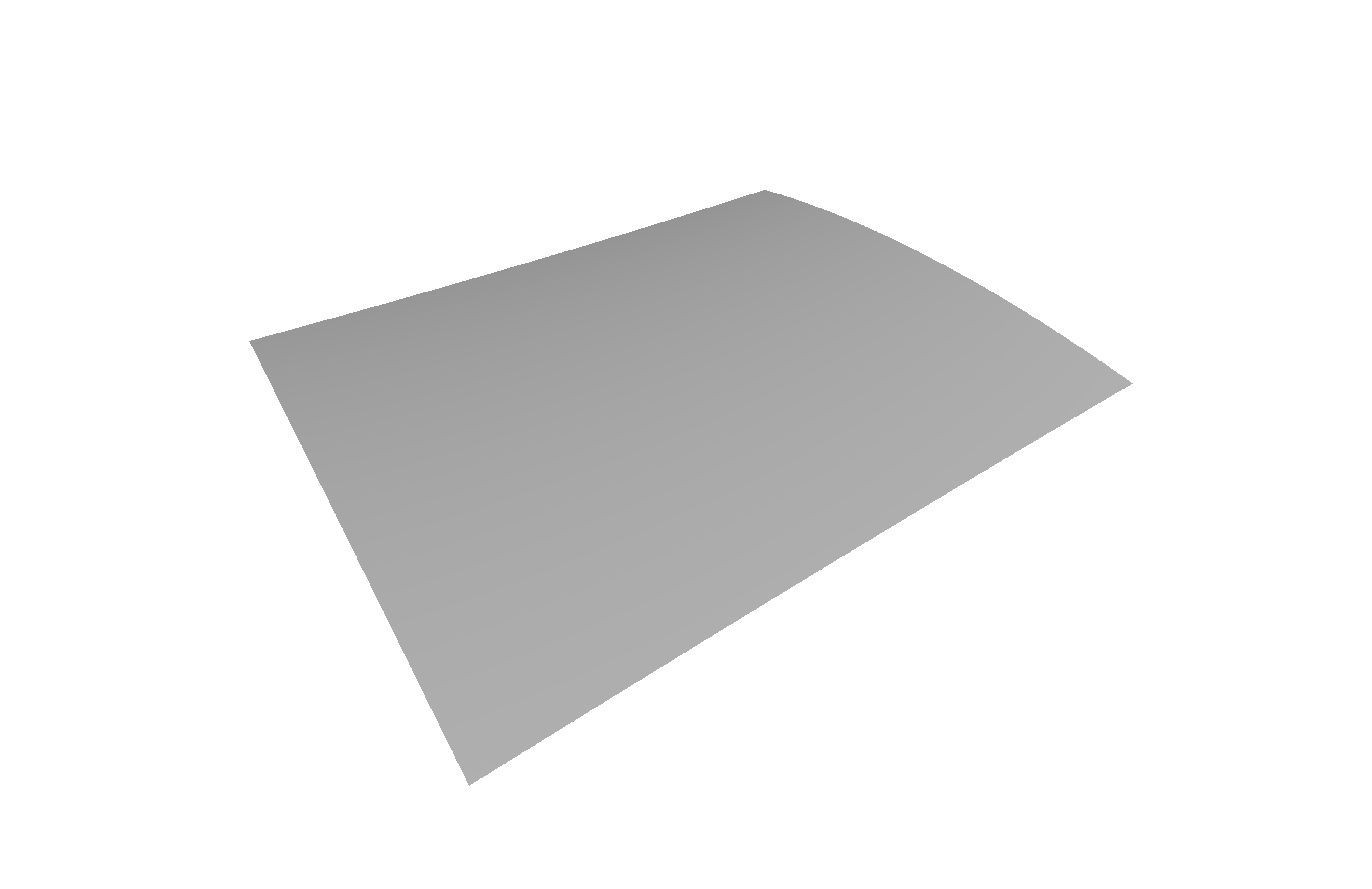}};
\node[inner sep=0pt] (c12) at (4.5,3.5)
    {\includegraphics[scale=0.0375,trim={12.5cm 4cm 12.5cm 4cm}, clip]{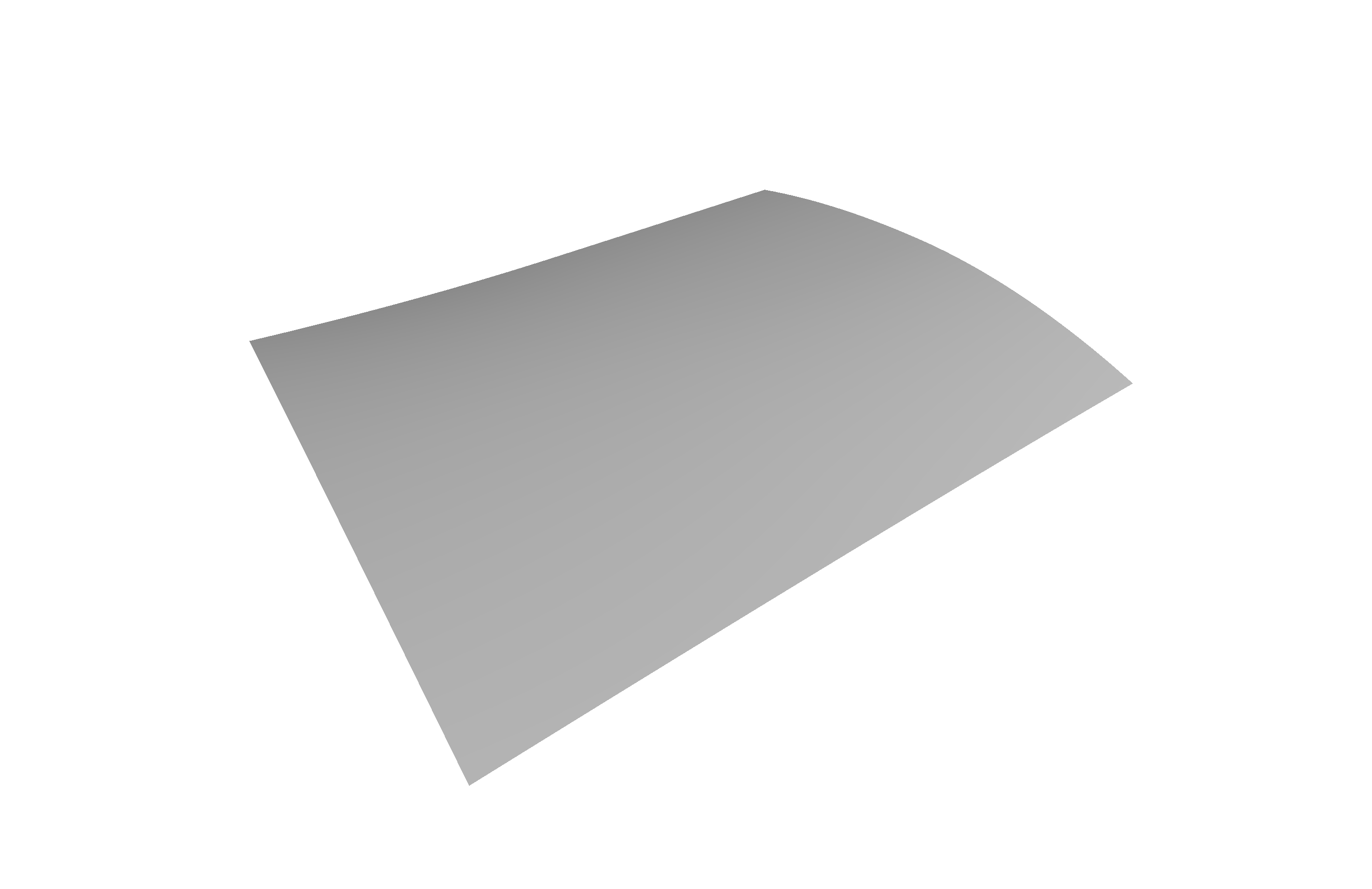}};

\node[inner sep=0pt] (c13) at (9,3.5)
    {\includegraphics[scale=0.0375,trim={12.5cm 4cm 12.5cm 4cm}, clip]{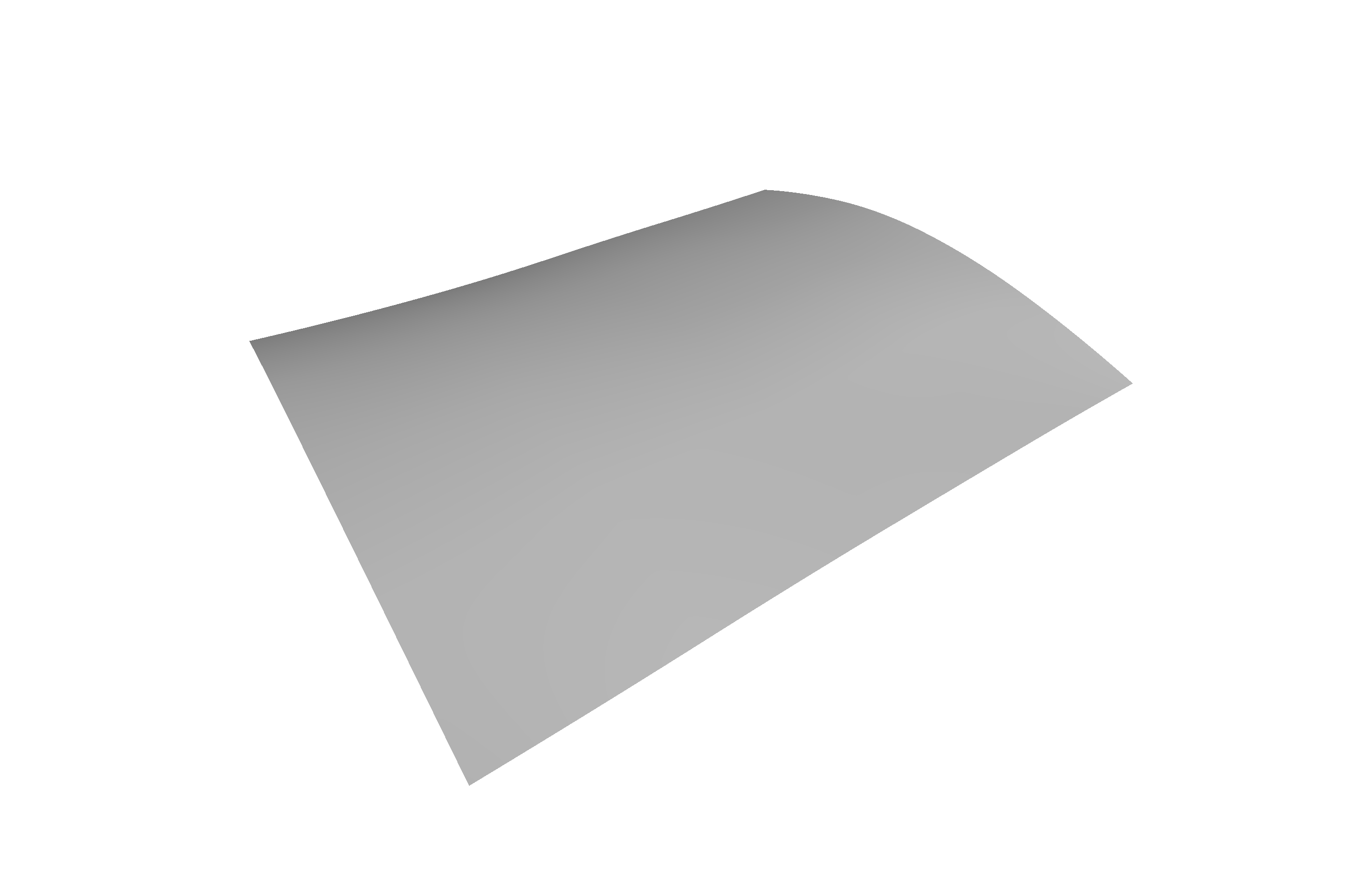}};

\node[inner sep=0pt] (c14) at (13.5,3.5)
    {\includegraphics[scale=0.0375,trim={12.5cm 4cm 12.5cm 4cm}, clip]{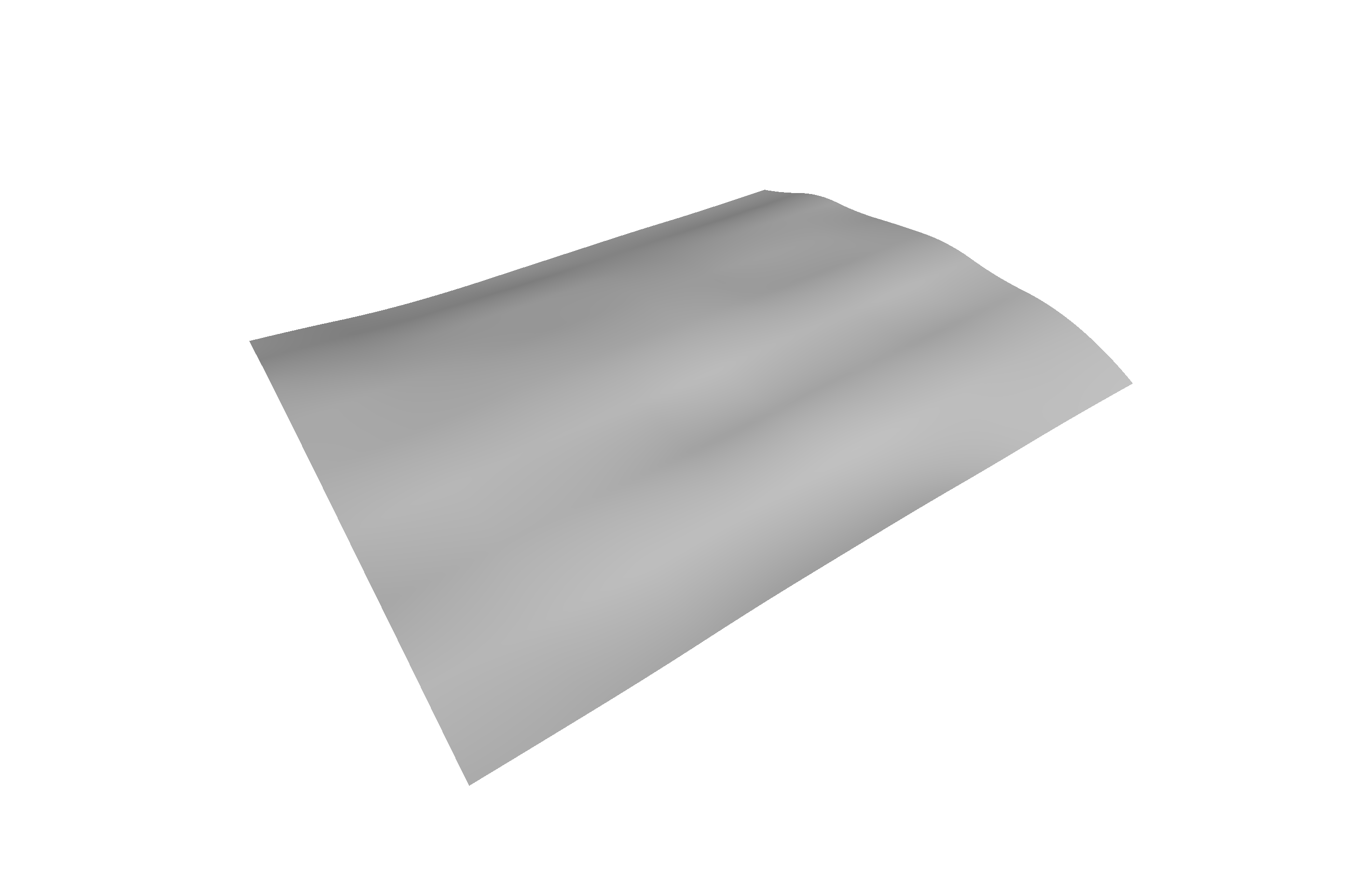}};
    
\node[inner sep=0pt] (c15) at (13.5,0)
    {\includegraphics[scale=0.0375,trim={12.5cm 4cm 12.5cm 4cm}, clip]{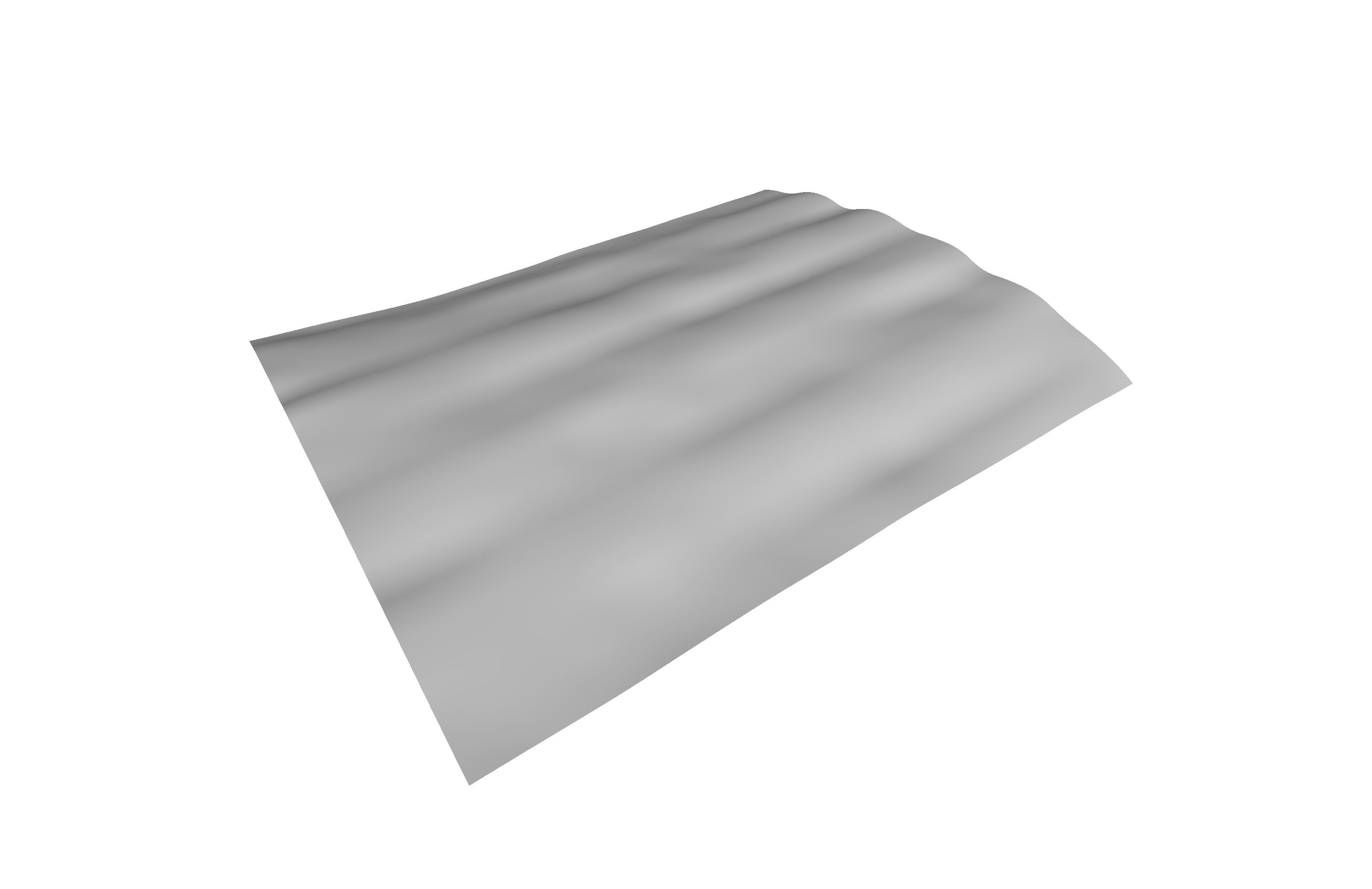}};
    
\node[inner sep=0pt] (c16) at (9,0)
    {\includegraphics[scale=0.0375,trim={12.5cm 4cm 12.5cm 4cm}, clip]{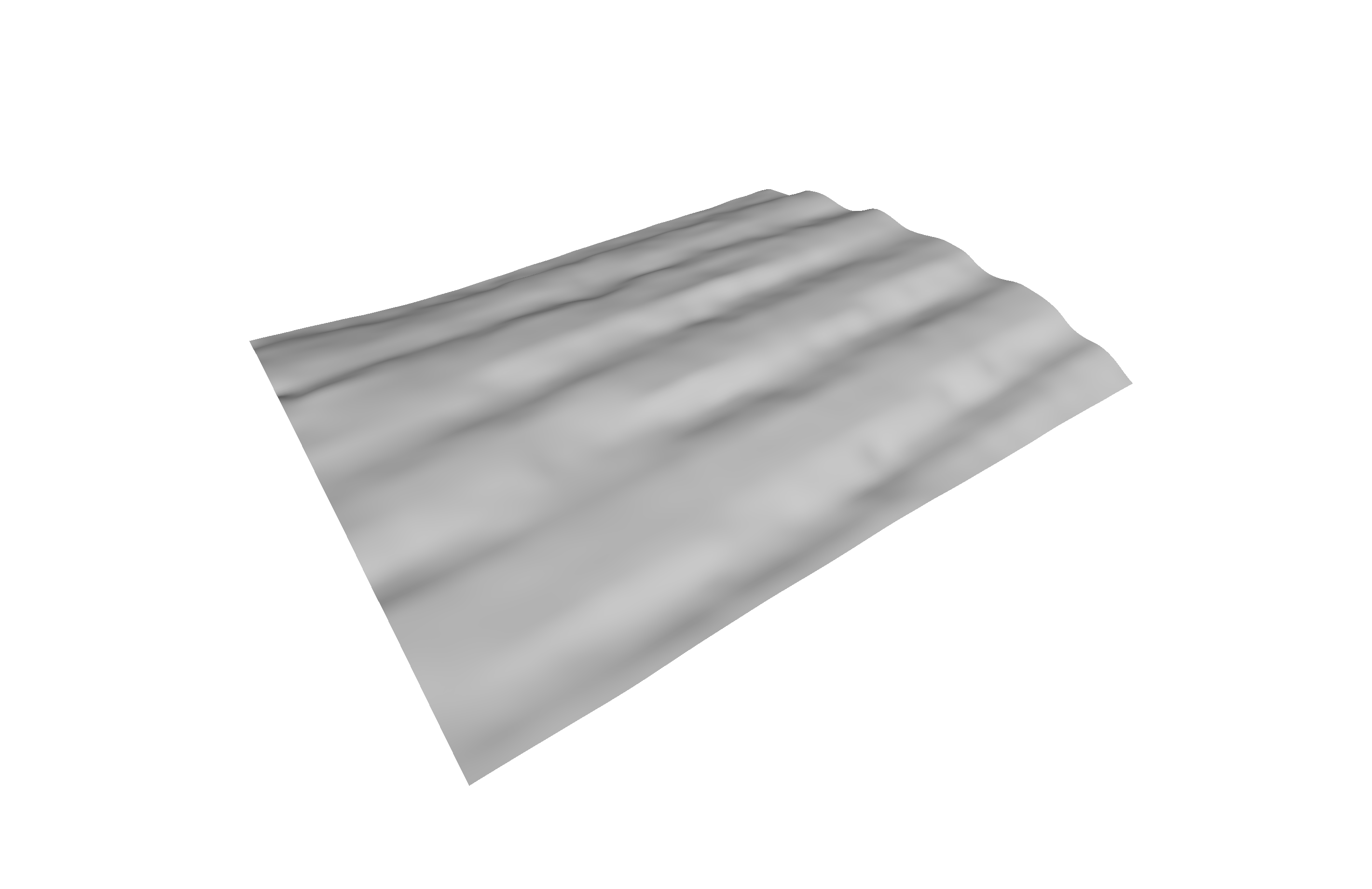}};

\node[inner sep=0pt] (c17) at (4.5,0)
    {\includegraphics[scale=0.0375,trim={12.5cm 4cm 12.5cm 4cm}, clip]{fringe00_L05.png}};

\node[inner sep=0pt] (c18) at (0,0)
    {\includegraphics[scale=0.0375,trim={12.5cm 4cm 12.5cm 4cm}, clip]{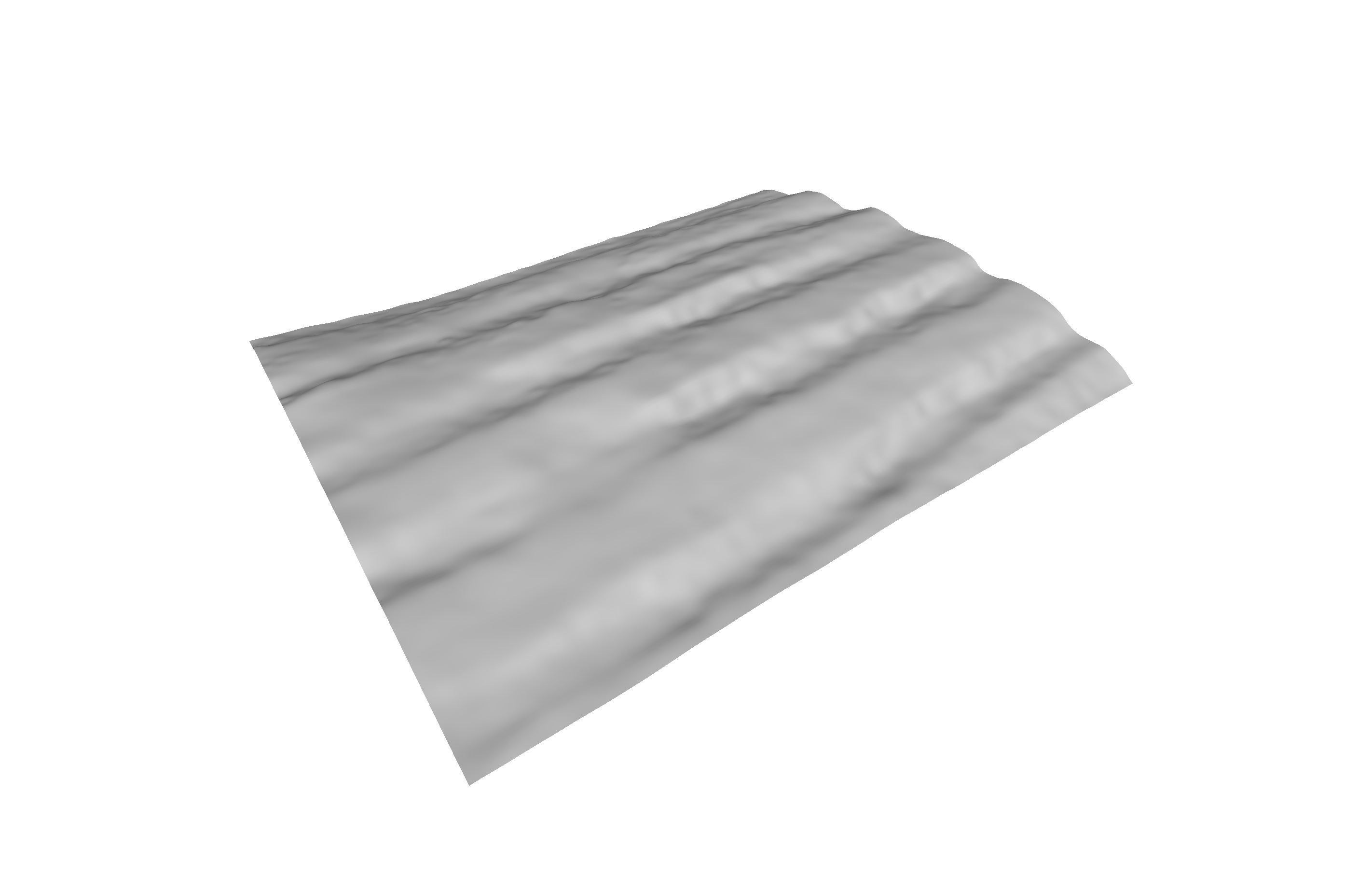}};
    
\draw [->,thick] (c11) -- (c12);
\draw [->,thick] (c12) -- (c13);
\draw [->,thick] (c13) -- (c14);
\draw [->,thick] (c14) -- (c15);
\draw [->,thick] (c15) -- (c16);
\draw [->,thick] (c16) -- (c17);
\draw [->,thick] (c17) -- (c18);

\end{tikzpicture}
}
\textbf{A.2} The MSE over the test set is $2.16\cdot{}10^{-6}$.

\tcbox{
\begin{tikzpicture}[scale=0.80]
\node[inner sep=0pt] (a11) at (0,3.5)
    {\includegraphics[scale=0.0375,trim={12.5cm 4cm 12.5cm 4cm}, clip]{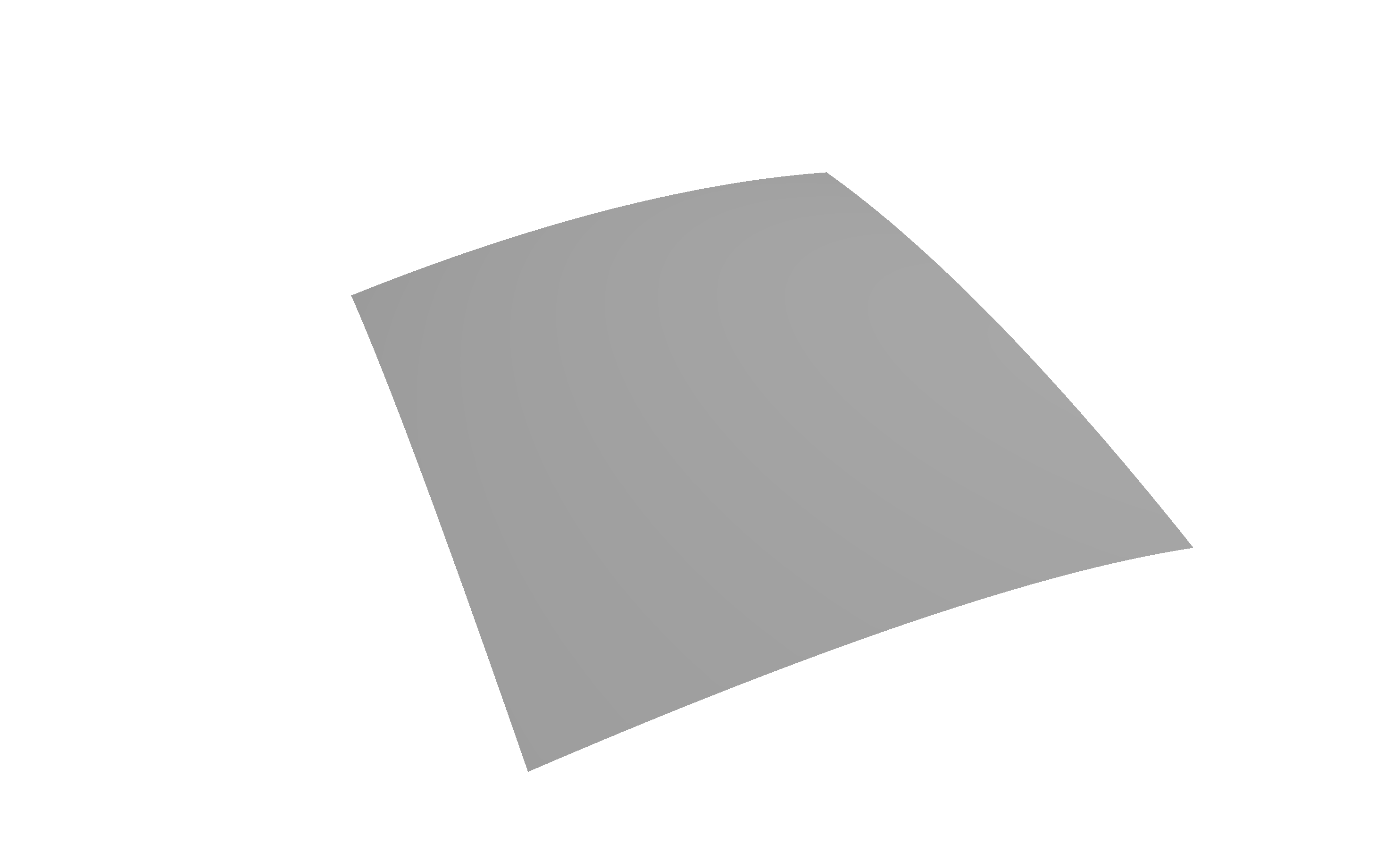}};
    
\node[inner sep=0pt] (a12) at (4.5,3.5)
    {\includegraphics[scale=0.0375,trim={12.5cm 4cm 12.5cm 4cm}, clip]{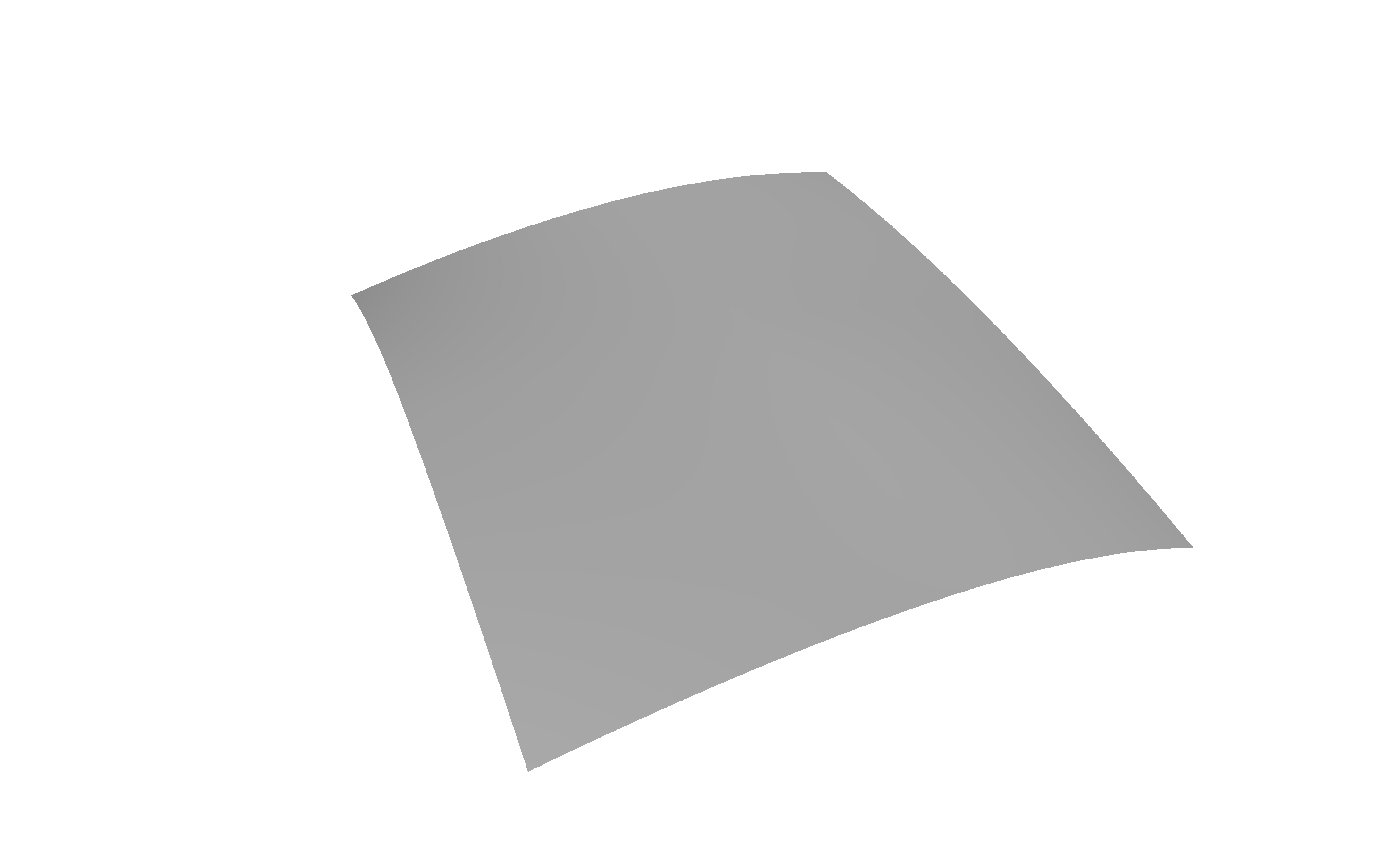}};

\node[inner sep=0pt] (a13) at (9,3.5)
    {\includegraphics[scale=0.0375,trim={12.5cm 4cm 12.5cm 4cm}, clip]{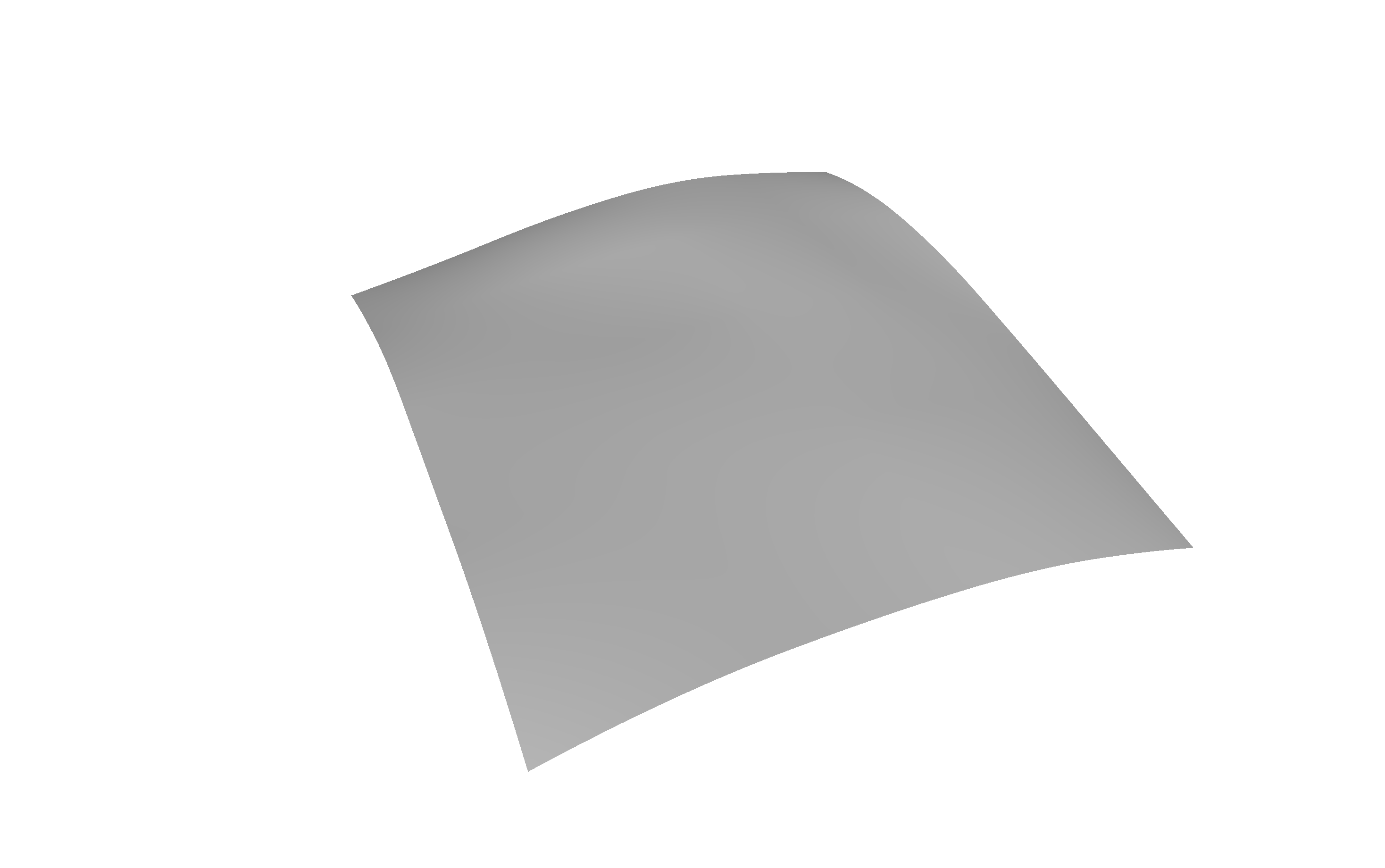}};

\node[inner sep=0pt] (a14) at (13.5,3.5)
    {\includegraphics[scale=0.0375,trim={12.5cm 4cm 12.5cm 4cm}, clip]{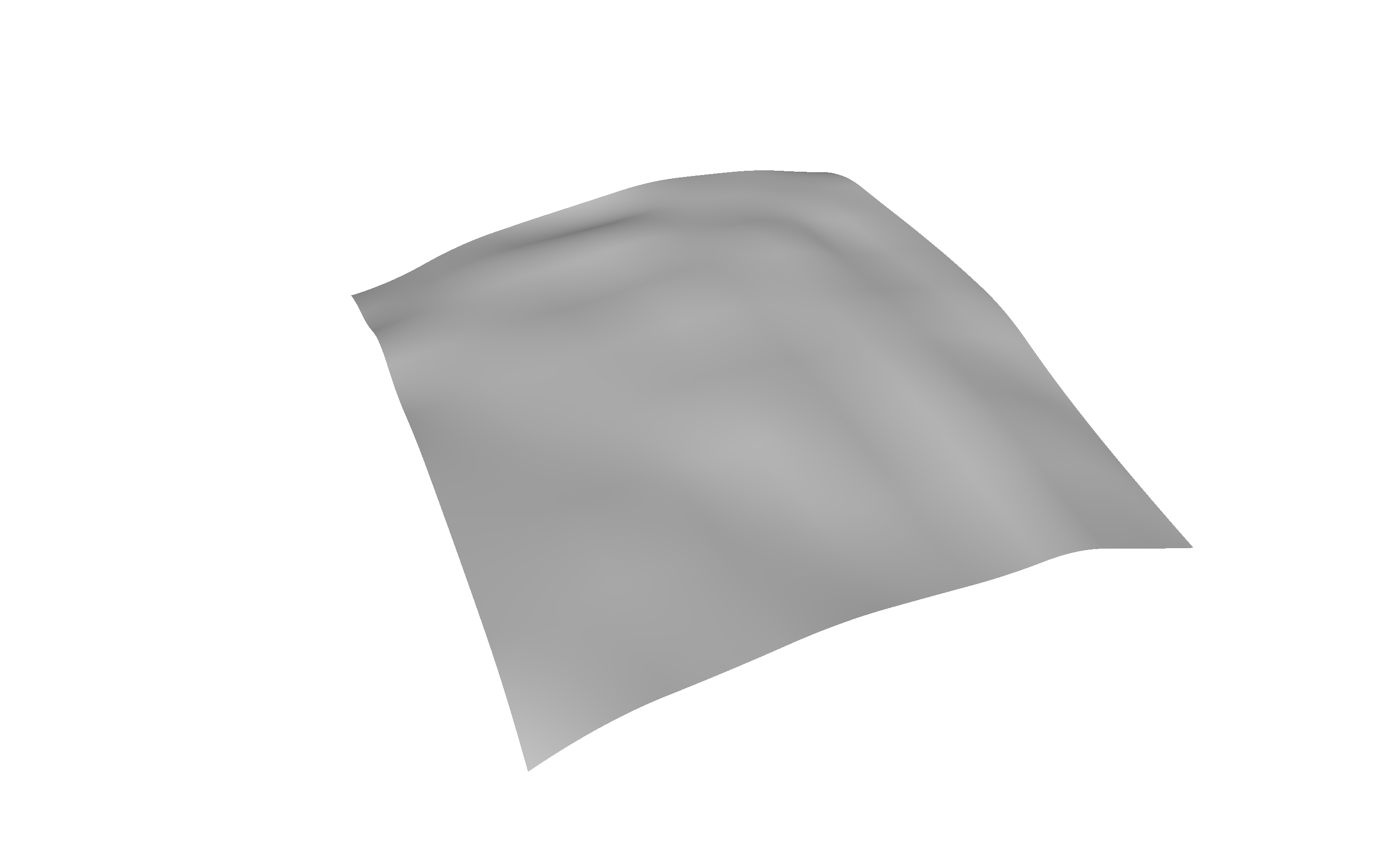}};
    
\node[inner sep=0pt] (a15) at (13.5,0)
    {\includegraphics[scale=0.0375,trim={12.5cm 4cm 12.5cm 4cm}, clip]{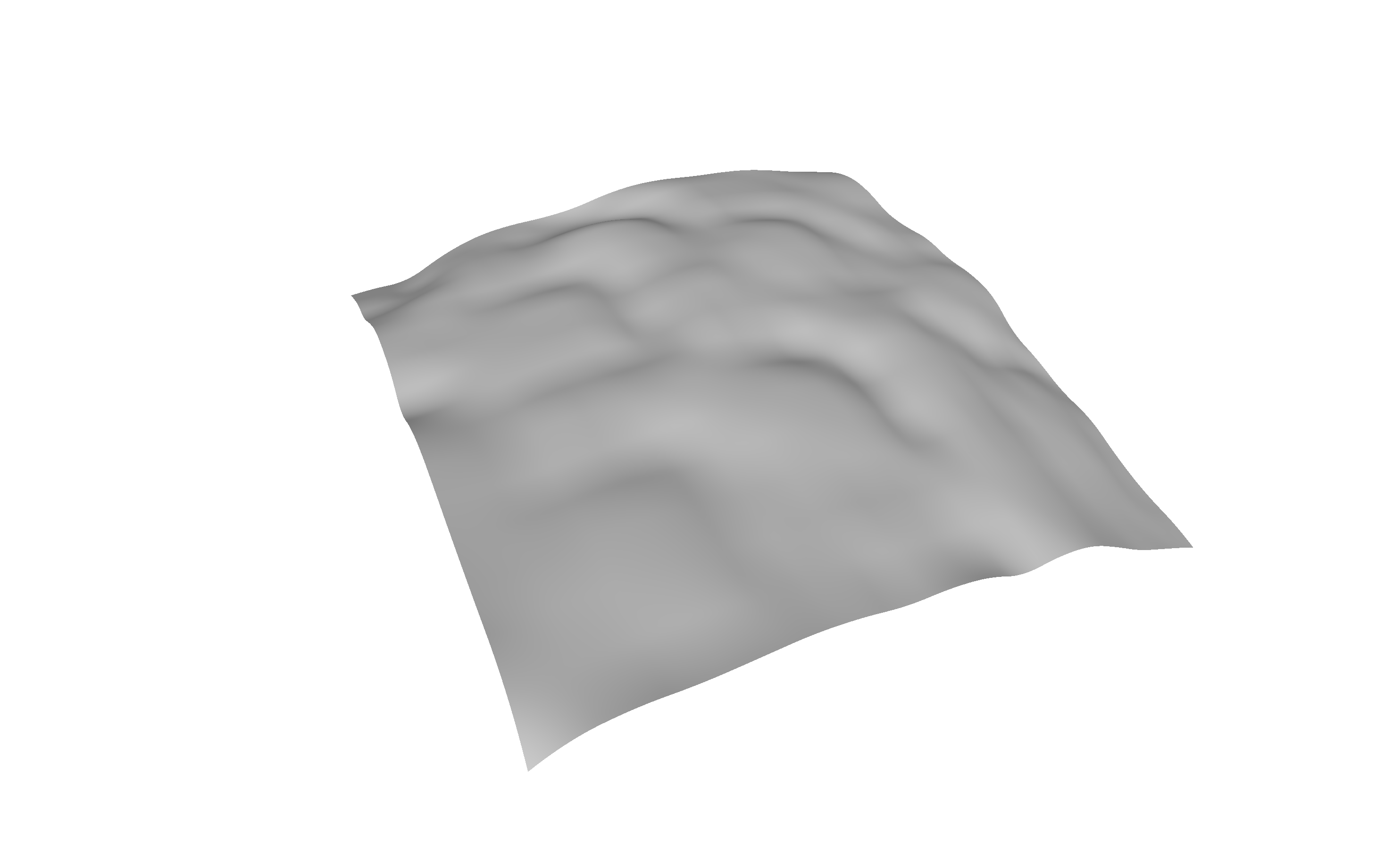}};
    
\node[inner sep=0pt] (a16) at (9,0)
    {\includegraphics[scale=0.0375,trim={12.5cm 4cm 12.5cm 4cm}, clip]{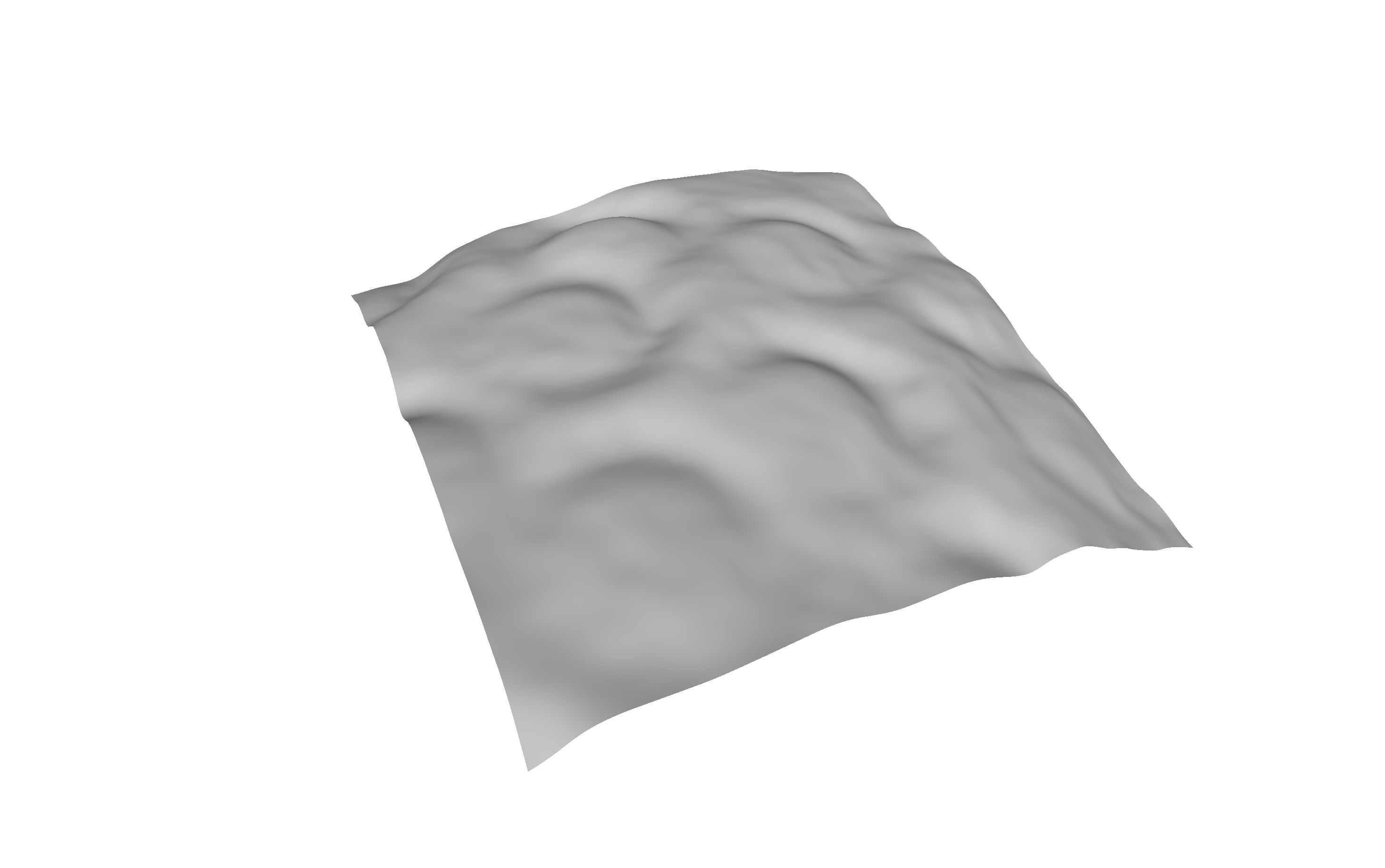}};

\node[inner sep=0pt] (a17) at (4.5,0)
    {\includegraphics[scale=0.0375,trim={12.5cm 4cm 12.5cm 4cm}, clip]{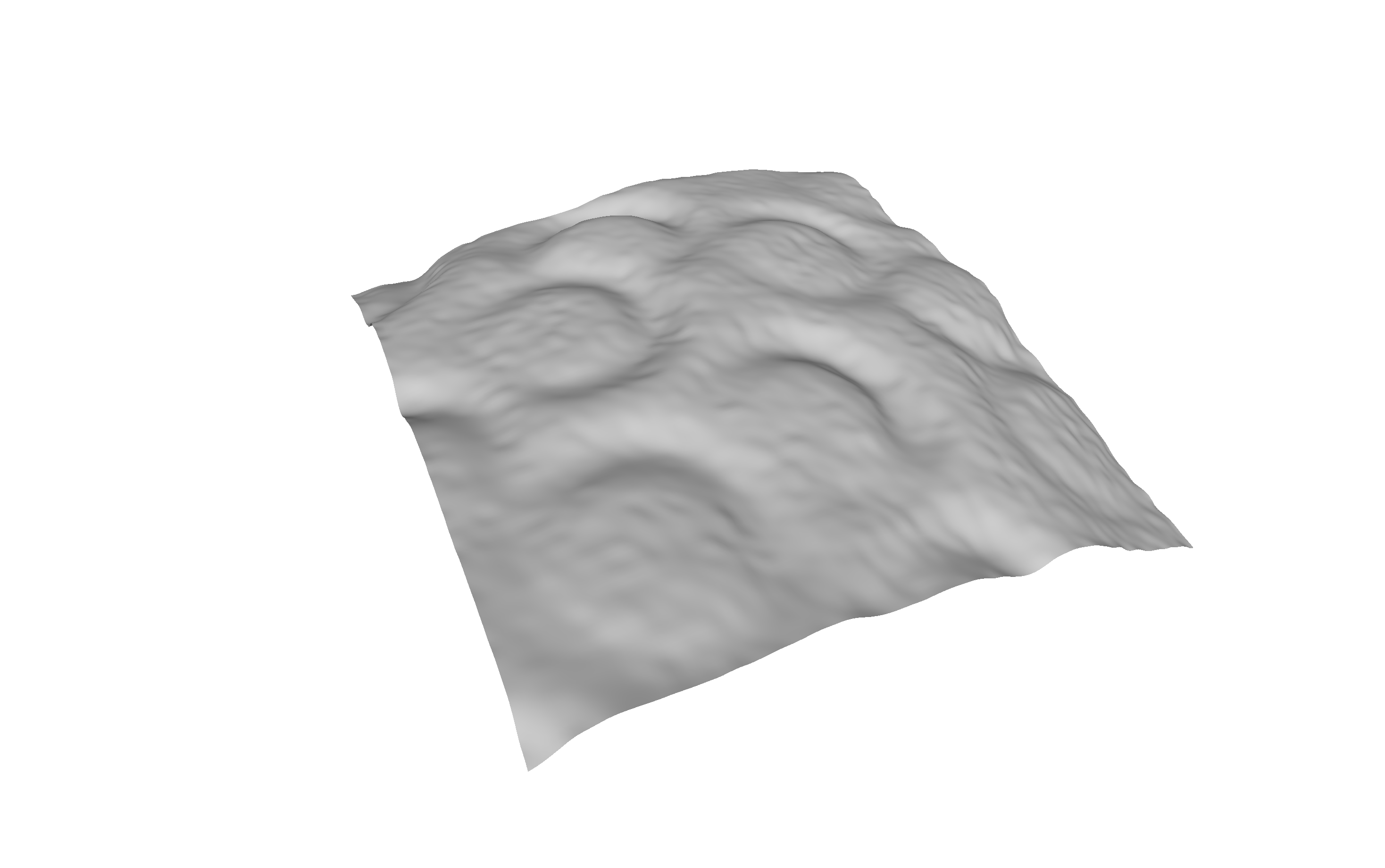}};

\node[inner sep=0pt] (a18) at (0,0)
    {\includegraphics[scale=0.0375,trim={12.5cm 4cm 12.5cm 4cm}, clip]{fringe200_L07.png}};
    
\draw [->,thick] (a11) -- (a12);
\draw [->,thick] (a12) -- (a13);
\draw [->,thick] (a13) -- (a14);
\draw [->,thick] (a14) -- (a15);
\draw [->,thick] (a15) -- (a16);
\draw [->,thick] (a16) -- (a17);
\draw [->,thick] (a17) -- (a18);

\end{tikzpicture} 
}
\textbf{A.3} The MSE over the test set is $6.92\cdot{}10^{-9}$.

\caption{Representation step by step of the approximations of the artifacts in Fig. \ref{figure:3Dmodels} via wQISA with Gaussian weight function.}
\label{fig:iterative_3Dmodels}
\end{center}
\end{figure}

Figure \ref{figure:3Dmodels} (right column) depicts the local reconstructions obtained via wQISA. The approximations show a correct recovery of the main details of the artefacts, with an accurate shape preservation: in A.1, the symmetry of the bulge at the center of the spiral is maintained; in A.2, 
we can notice a faithful reconstruction even in presence of small geometric variations, such as chiselings or local reliefs; finally, A.3 shows a reconstruction characterized by circular patterns.
Figure \ref{fig:iterative_3Dmodels} depicts the iterative process that leads to the final wQISA approximations and reports the MSE values computed in Step 4 of Algorithm \ref{algorithm:pseudo_algorithm}.

We refer to Table \ref{table:3Dmodels} for a detailed quantitative comparison of the approximation performances of wQISA, RBFs, Kriging and MBA of the models in Figure \ref{figure:3Dmodels}. We highlight in bold the values with the lowest order of magnitude, which  correspond to the best performances. This Table also lists the required number of iterations. From these experiments we can conclude that wQISA performs comparatively well with respect to the other methods when the standard deviation and the Hausdorff distance are considered. This is not surprising because the idea behind wQISA is to fit the overall data trend rather that converging to single values. Tables \ref{table:mesh} reports the size of the optimal grid for wQISAs and MBAs. The optimal smoothing parameters for the RBF approximations can be found in Table \ref{table:RBF}.

\begin{table}[!h]
\centering
\small
\caption{Local approximation of 3D models. The labels A.1, A.2 and A.3 refer to the models in Figure \ref{figure:3Dmodels}. In bold, we highlight the errors with the lowest order of magnitude. The single asterisk highlights those cases where the maximum number of iterations, here set to $15$, is reached. Kriging runs in a single step and no number is displayed. \label{table:3Dmodels}}
\scalebox{0.9}{
\begin{tabular}{l l c c c c c}
\toprule
& & \multicolumn{3}{c}{Punctual error} & \multirow{2}{*}{Hausdorff} & \multirow{2}{*}{\makecell{\small Number \\ of \\ Iterations}}  \\ 
\cmidrule(l){3-5} & & mean &  st.d. & MSE  & \\[1.5ex]

\midrule
\multirow{7}{*}{A.1} & \multicolumn{1}{|l}{wQISA} &   $4.445{\cdot}10^{-5}$ &  $\mathbf{6.881{\cdot}10^{-5}}$ & $\mathbf{6.723{\cdot}10^{-9}}$  & $\mathbf{9.229{\cdot}10^{-4}}$ &  \multicolumn{1}{|c}{7}   \\

 & \multicolumn{1}{|l}{RBF\textsuperscript{1}} & 
  $1.194{\cdot}10^{-5}$ &  $\mathbf{3.399{\cdot}10^{-5}}$ & $\mathbf{1.305{\cdot}10^{-9}}$  & $\mathbf{9.303{\cdot}10^{-4}}$ & \multicolumn{1}{|c}{$15^\ast$} \\

 & \multicolumn{1}{|l}{RBF\textsuperscript{2}} & 
  $3.94{\cdot}10^{-3}$ &  $3.052{\cdot}10^{-3}$ & $2.480{\cdot}10^{-5}$  & $1.968{\cdot}10^{-2}$ & \multicolumn{1}{|c}{13} \\
 
 & \multicolumn{1}{|l}{RBF\textsuperscript{3}} & 
   $4.496{\cdot}10^{-3}$ &  $3.635{\cdot}10^{-3}$ & $3.335{\cdot}10^{-5}$  & $2.436{\cdot}10^{-2}$ & \multicolumn{1}{|c}{$15^\ast$} \\
 
 & \multicolumn{1}{|l}{RBF\textsuperscript{4}} & 
   $2.760{\cdot}10^{-3}$ &  $2.241{\cdot}10^{-3}$ & $1.261{\cdot}10^{-5}$  & $1.425{\cdot}10^{-2}$ & \multicolumn{1}{|c}{$15^\ast$} \\
 
 & \multicolumn{1}{|l}{Kriging} &   $2.801{\cdot}10^{-5}$ &  $\mathbf{5.352{\cdot}10^{-5}}$ & $\mathbf{3.650{\cdot}10^{-9}}$  & $\mathbf{9.307{\cdot}10^{-4}}$ & \multicolumn{1}{|c}{-} \\
     
 & \multicolumn{1}{|l}{MBA} & 
  $\mathbf{9.451{\cdot}10^{-6}}$ &  $\mathbf{3.137{\cdot}10^{-5}}$ & $\mathbf{1.071{\cdot}10^{-9}}$  & $\mathbf{9.301{\cdot}10^{-4}}$ & \multicolumn{1}{|c}{7}  \\

\midrule
\multirow{7}{*}{A.2} & \multicolumn{1}{|l}{wQISA} &   $7.046{\cdot}10^{-4}$ &  $\mathbf{9.590{\cdot}10^{-4}}$ & $1.420{\cdot}10^{-6}$  & $\mathbf{6.670{\cdot}10^{-2}}$ & \multicolumn{1}{|c}{$8$} \\

 & \multicolumn{1}{|l}{RBF\textsuperscript{1}}
 & 9.112 
 & 22.504 & 589.495  & $1.351{\cdot}10^{3}$ & \multicolumn{1}{|c}{$15^\ast$} \\

 & \multicolumn{1}{|l}{RBF\textsuperscript{2}} & 
   $\mathbf{8.043{\cdot}10^{-5}}$ &  $\mathbf{2.044{\cdot}10^{-4}}$ & $\mathbf{4.826{\cdot}10^{-8}}$  & $\mathbf{6.412{\cdot}10^{-2}}$ & \multicolumn{1}{|c}{$12$}\\
 
 & \multicolumn{1}{|l}{RBF\textsuperscript{3}} & 
 $1.097{\cdot}10^{-2}$ &  $4.990{\cdot}10^{-2}$ & $2.61{\cdot}10^{-3}$  & $2.264$ & \multicolumn{1}{|c}{$15^\ast$}\\
 
 & \multicolumn{1}{|l}{RBF\textsuperscript{4}} & 
   $2.199{\cdot}10^{-2}$ &  $2.973{\cdot}10^{-2}$ & $1.361{\cdot}10^{-3}$  & $1.154$ & \multicolumn{1}{|c}{$15^\ast$} \\
 
 & \multicolumn{1}{|l}{Kriging} &   $1.153{\cdot}10^{-3}$ &  $1.301{\cdot}10^{-3}$ & $3.027{\cdot}10^{-6}$ & $\mathbf{6.695{\cdot}10^{-2}}$ & \multicolumn{1}{|c}{-} \\
     
 & \multicolumn{1}{|l}{MBA} & 
  $1.214{\cdot}10^{-4}$ &  $\mathbf{2.503{\cdot}10^{-4}}$ & $\mathbf{7.738{\cdot}10^{-8}}$  & $\mathbf{6.695{\cdot}10^{-2}}$ & \multicolumn{1}{|c}{$13$}  \\

\midrule
\multirow{7}{*}{A.3} & \multicolumn{1}{|l}{wQISA} &   $1.054{\cdot}10^{-3}$ &  $\mathbf{9.926{\cdot}10^{-4}}$ & $2.081{\cdot}10^{-6}$  & $\mathbf{5.640{\cdot}10^{-2}}$ & \multicolumn{1}{|c}{$8$}  \\

 & \multicolumn{1}{|l}{RBF\textsuperscript{1}} & 
   $23.964$ & $98.608$ & $1.030{\cdot}10^{4}$  & $4.342{\cdot}10^{3}$ & \multicolumn{1}{|c}{$15^\ast$} \\

 & \multicolumn{1}{|l}{RBF\textsuperscript{2}} & 
   $\mathbf{1.440{\cdot}10^{-4}}$ &  $\mathbf{2.212{\cdot}10^{-4}}$ & $\mathbf{6.958{\cdot}10^{-8}}$  & $\mathbf{5.552{\cdot}10^{-2}}$ & \multicolumn{1}{|c}{$12$} \\
 
 & \multicolumn{1}{|l}{RBF\textsuperscript{3}} & 
  $3.008{\cdot}10^{-2}$ &  $4.139{\cdot}10^{-2}$ & $2.617{\cdot}10^{-2}$  & $7.891$ & \multicolumn{1}{|c}{$15^\ast$} \\
 
 & \multicolumn{1}{|l}{RBF\textsuperscript{4}} & 
  $8.760{\cdot}10^{-2}$ &  $1.189{\cdot}10^{-1}$ & $2.179{\cdot}10^{-2}$  & $7.066$ & \multicolumn{1}{|c}{$15^\ast$} \\
 
 & \multicolumn{1}{|l}{Kriging} &  $\mathbf{1.510{\cdot}10^{-4}}$ &  $\mathbf{1.601{\cdot}10^{-4}}$ & $\mathbf{4.845{\cdot}10^{-8}}$  & $\mathbf{5.527{\cdot}10^{-2}}$ & \multicolumn{1}{|c}{-} \\
     
 & \multicolumn{1}{|l}{MBA} & 
  $\mathbf{1.635{\cdot}10^{-4}}$ &  $\mathbf{2.185{\cdot}10^{-4}}$ & $\mathbf{7.446{\cdot}10^{-8}}$  & $\mathbf{5.506{\cdot}10^{-2}}$ & \multicolumn{1}{|c}{$10$}  \\
 
\bottomrule
\end{tabular}
}
\end{table}

\begin{table}[!h]
\centering
\caption{Size of the tensor mesh for wQISA and MBA.}
\label{table:mesh}
\scalebox{0.9}{
\begin{tabular}{lccc}
\toprule
\textbf{Method} & \textbf{A.1} & \textbf{A.2} & \textbf{A.3}\\
\midrule
wQISA                  & $50\times{}64$ & $128\times{}128$ & $128\times{}128$ \\
MBA                    & $128\times{}128$ & $4096\times{}4096$ & $128\times{}128$ \\
\bottomrule
\end{tabular}
}
\end{table}

\begin{table}[!h]
\centering
\caption{Optimal smoothing parameter $\alpha$ for the RBFs implementations.}
\label{table:RBF}
\scalebox{0.9}{
\begin{tabular}{lccc}
\toprule
\textbf{Method} & \textbf{A.1} & \textbf{A.2} & \textbf{A.3}\\
\midrule
RBF\textsuperscript{1} & $1.27{\cdot}10^{-3}$ & $1.49{\cdot}10^{-1}$ & $1.49{\cdot}10^{-1}$\\
RBF\textsuperscript{2} & $2.37{\cdot}10^{-1}$ & $3.49{\cdot}10^{-2}$ & $3.43{\cdot}10^{-2}$ \\
RBF\textsuperscript{3} & $2.83{\cdot}10^{-1}$ & $4.54{\cdot}10^{-2}$ & $2.89{\cdot}10^{-1}$\\
RBF\textsuperscript{4} & $2.86{\cdot}10^{-1}$ & $1.42{\cdot}10^{-1}$ & $1.42{\cdot}10^{-1}$\\
\bottomrule
\end{tabular}
}
\end{table}

\begin{figure*}[h!]
\centering
    \begin{tabular}{ccc}
        \includegraphics[scale=0.225, trim={1cm 0cm 1cm 0cm},clip]{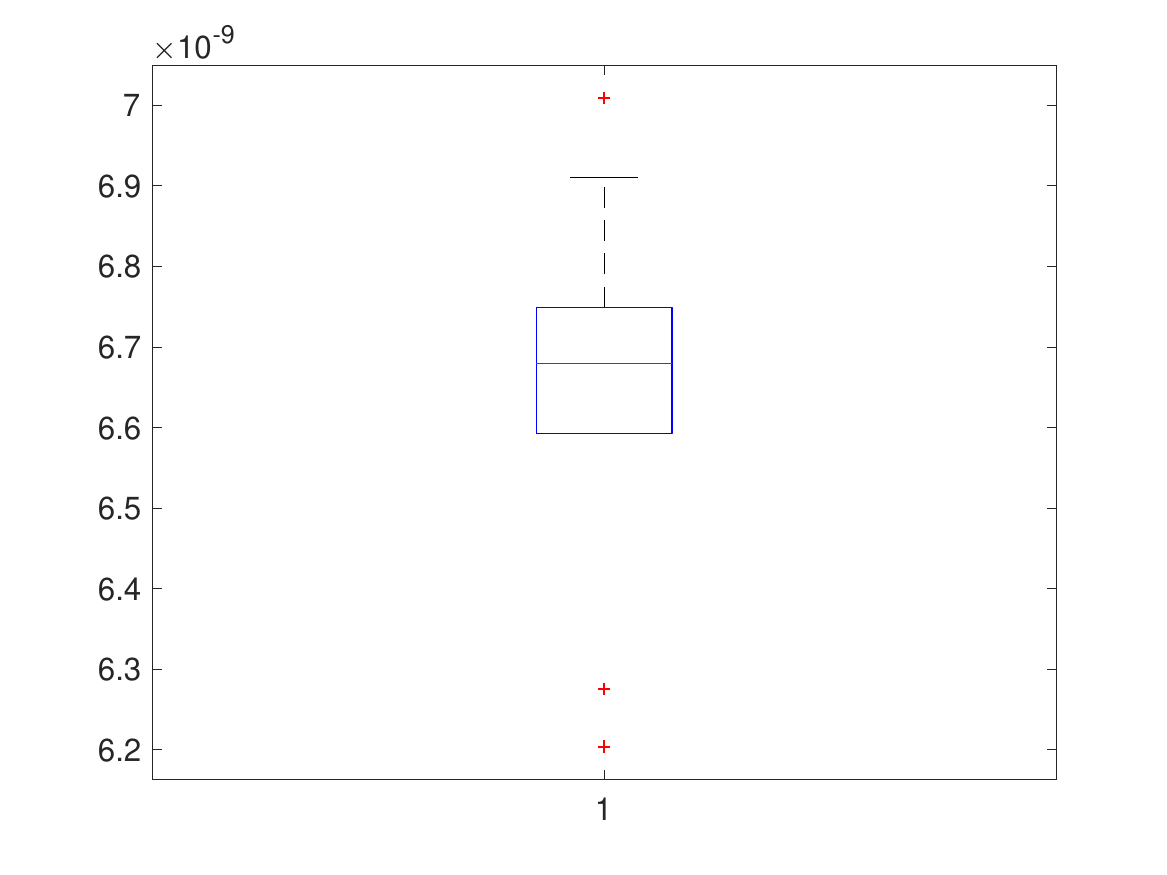}
        &
        \includegraphics[scale=0.225, trim={1cm 0cm 1cm 0cm},clip]{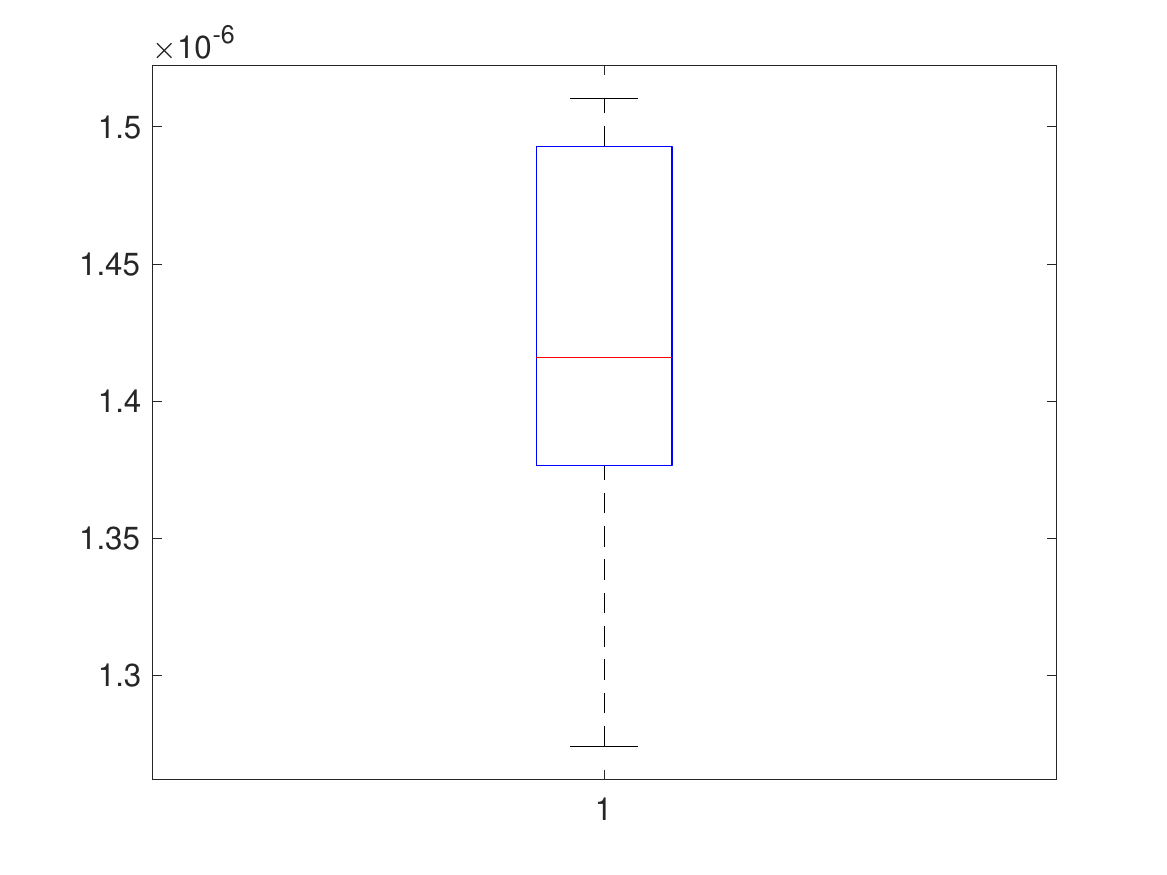}
        &
        \includegraphics[scale=0.225, trim={1cm 0cm 1cm 0cm},clip]{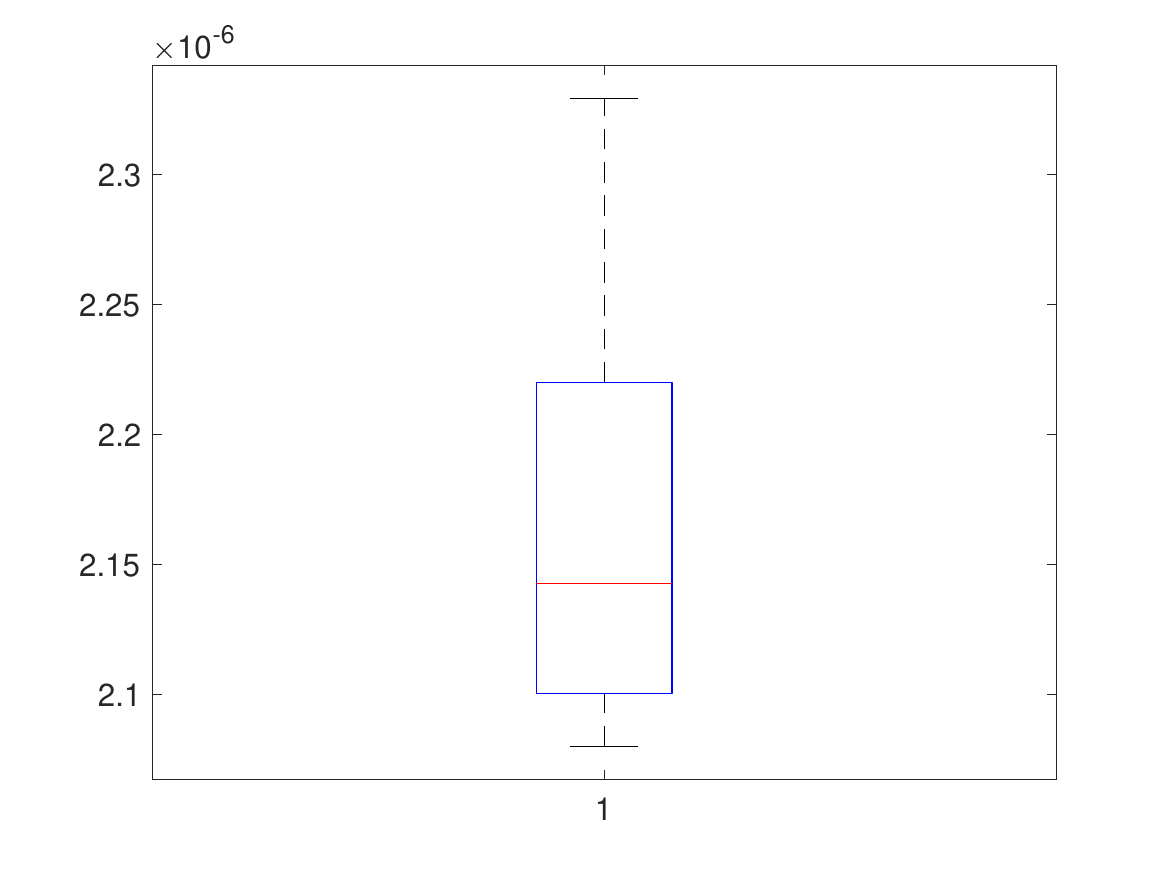} \\
       (A.1) & (A.2) & (A.3)
    \end{tabular}
    \caption{Robustness to different training, validation, and test sets for the models in Figure \ref{figure:3Dmodels}. We show with box-plots the MSE over $10$ validation sets.}
    \label{figure:boxplots_3Dmodels}
\end{figure*}

The robustness of the data-driven implementation to different choices of the training, validation and test sets (see Section \ref{data_pre_processing}) is here assessed by considering ten different data splits for each 3D model. The resulting MSEs over the respective validation point clouds are graphically represented in the form of boxplots (see Figure \ref{figure:boxplots_3Dmodels}).

\newpage
\subsection{Terrain data modelling}

Acquisition methods such as Lidar, photogrammetry and sonar produce huge amounts of data that may be perturbed. Methods to efficiently and effectively handle these point clouds are thus required. Since the Kriging and RBFs implementations we are considering are not tailored for big data and present some scalability limitations, here, we limit our attention to the two local methods: wQISA and MBA.

The first data set comes from the island municipality of Værøy and consists of many islands. The data are provided by the Norwegian Mapping Authority Kartverket, and are freely downloadable at \url{https://hoydedata.no/LaserInnsyn/}. The original point cloud is pre-processed to remove the data outside the mainland, and reduced to 44,529 points. The peculiarity of this terrain is that alternate regions with high variability and almost flat regions.

The second data set corresponds to submarine sand dunes off the coast of France. The data were obtained using bathymetric surveys and are detailed in \cite{FRANZETTI201317}. The selected part of the data set contains 759,952 points. The difference in sea depth is about 33.2 m. This terrain was already adopted in \cite{Skytt2015}, it is relatively large and presents slight ripples on a wide wavefront.

\begin{figure}[!h]
\begin{center}
\begin{tabular}{cc}
\includegraphics[scale=0.225]{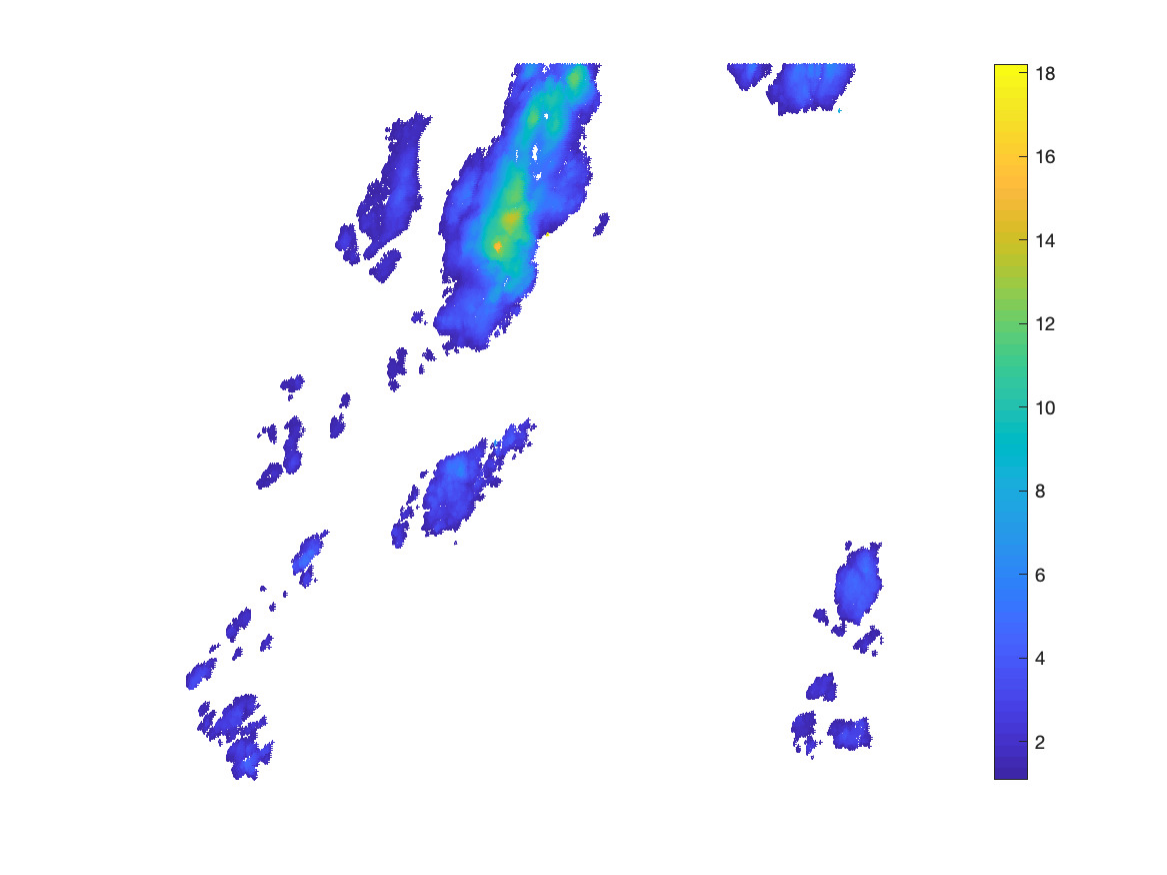}
&
\includegraphics[scale=0.225]{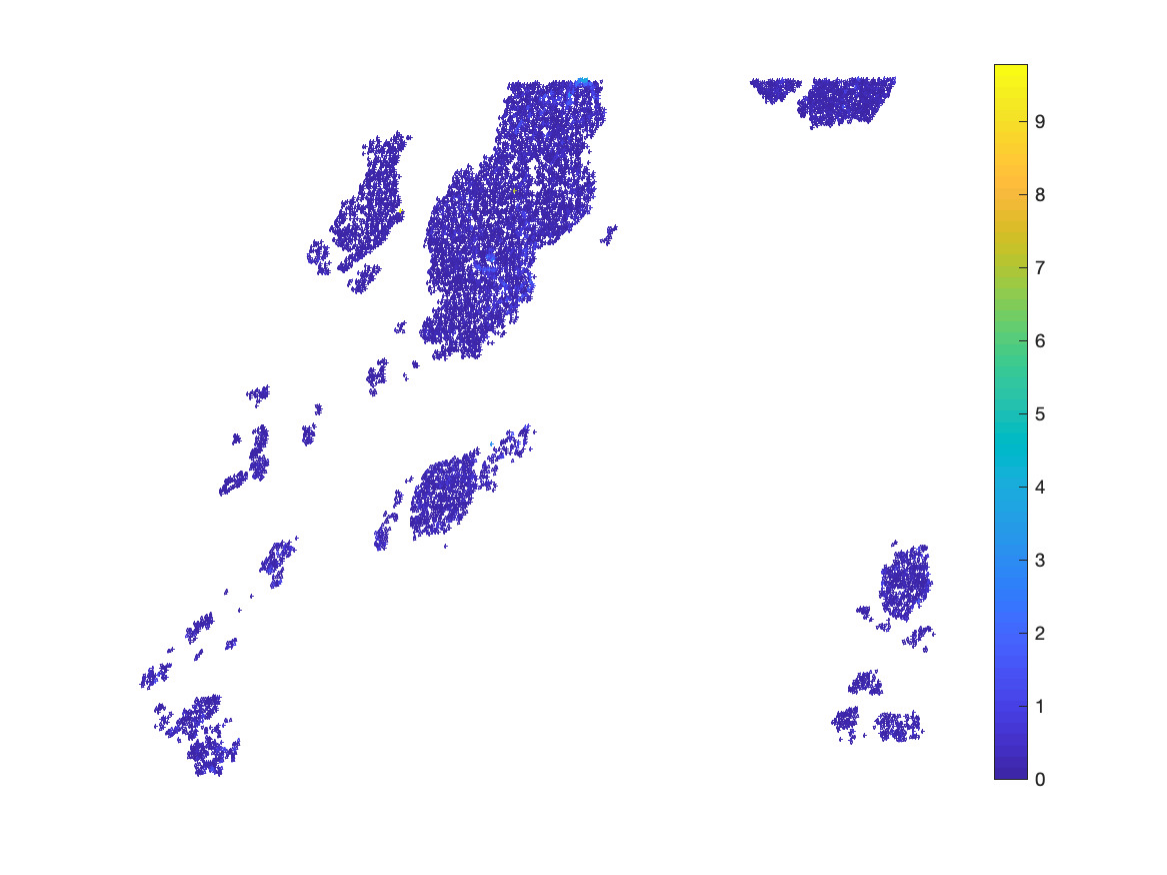} \\
(a) & (b) 
\end{tabular}
\end{center}
\caption{Værøy. Figure (a) shows a wQISA of the islands. Figure (b) shows the absolute punctual error.}
\label{figure:varoy}
\end{figure}

\begin{figure}[!h]
\begin{center}
\begin{tabular}{AA}
\includegraphics[scale=0.225]{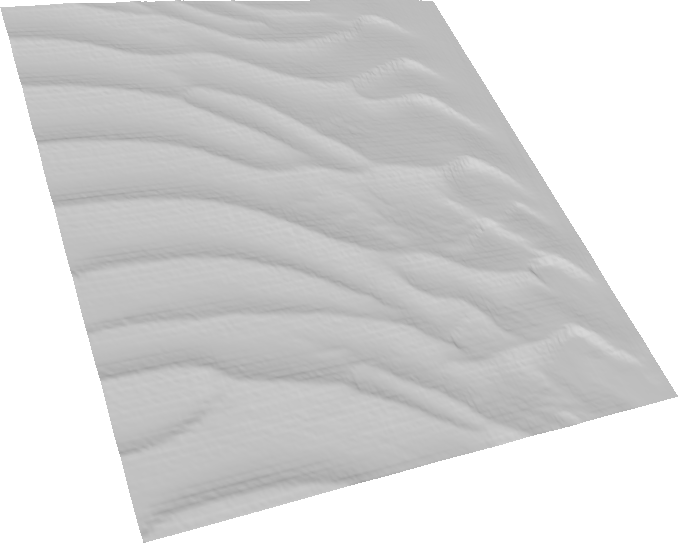}
&
\includegraphics[scale=0.625]{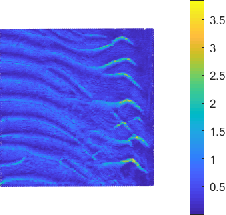} \\
(a) & (b) \\
\end{tabular}
\caption{Sand dunes. Figure (a) shows a wQISA of the sand dunes (north is pointing upwards to the left). Figure (b) shows the absolute punctual error.}
\label{figure:sand_dunes}
\end{center}
\end{figure}

For both models, we consider a wQISA representation with IDW weight (see Equation \ref{equation:IDW_weight}). Figures \ref{figure:varoy} shows the wQISA representation of the first terrain, together with the absolute punctual error. The optimal mesh is reached in 11 iteration ($1024\times1024$ elements) for wQISA  and in 10 iterations ($512\times{}512$ elements) for MBA. Similarly, the approximation and the absolute punctual error of the second terrain are shown in Figure \ref{figure:sand_dunes}. 
For this model, the optimal mesh for wQISA contains $6400\times{}6400$ elements while for MBA it contains $1024\times{}1024$ elements.

\newcolumntype{A}{ >{\centering\arraybackslash} m{4.0cm} }

A quantitative comparison of the wQISA and the MBA performances can be found in Table \ref{table:terrain_data}. MBA provides a better punctual approximation (highlighted by the MSE measure), but it is again less accurate than wQISA for the overall fit (as highlighted by the Hausdorff distance).

\begin{table}
\centering
\caption{Terrain data: Værøy (B.1) and sand dunes (B.2).\label{table:terrain_data}}
\scalebox{0.9}{
\begin{tabular}{l l c c c c}
\toprule
& & \multicolumn{3}{c}{Punctual error} &  \\ 
\cmidrule(l){3-5}
&  & mean &  st.d. & MSE & Hausdorff\\

\hline

\multirow{2}{*}{B.1} & \multicolumn{1}{|l}{wQISA} & 0.2158 & 0.2582 & 0.1132 & \textbf{9.6266}  \\
     
& \multicolumn{1}{|l}{MBA} & \textbf{0.1513} & \textbf{0.2257} & \textbf{0.0738} & 13.1604   \\

\hline

\multirow{2}{*}{B.2} & \multicolumn{1}{|l}{wQISA} & 0.0492 & 0.0641 & 0.0065 & \textbf{8.2569}  \\
     
& \multicolumn{1}{|l}{MBA} & \textbf{0.0330} &  \textbf{0.0543} & \textbf{0.0040} & 8.2579   \\

\bottomrule
\end{tabular}
}
\end{table}

\subsection{Data measurements}
We analyse the approximation of strongly perturbed data sets. We consider two precipitation events, with the aim of deepening the comparisons previously presented in \cite{Patane2017}.
Although focusing here on a specific environmental application, we highlight that these approximation methods do not limit to rainfall field simulations or to a specific sample size. We propose in the additional material the study of a rather small data set (just $11$ samples) measuring air pollution. In this case the B-spline approximations rather reduces to linear polynomial approximations, and thus we can consider it as a limiting case.

The first event occurred on September 29, 2013, and was characterized by light rain over Liguria with two different thunderstorms that caused local flooding and landslides. This data set counts only $143$ measurements, all acquired by the rain gauge network maintained by Regione Liguria, which is spread over the whole region. A quantitative comparison is provided in Table \ref{table:rainfall_fields} and is computed by performing a LOO cross-validation on each method. The optimal parameters are here listed:
\begin{itemize}
    \item \emph{wQISA with k-NN weight}. The optimal number of neighbors $k$ is $k=2$, where we have checked for any $k=1,\ldots,10$. The optimal mesh has $8\times8$ elements.
    \item \emph{RBF}. The optimal smoothing parameter $\alpha$ for the RBF approximations is computed by minimizing the mean squared error: $\alpha=1.00$ for the Gaussian kernel (RBF\textsuperscript{1}); $\alpha=1.22$ for the RBF with multiquadric kernel (RBF\textsuperscript{2}); $\alpha=1.00$ for the RBF with inverse kernel (RBF\textsuperscript{3}); $\alpha=34.31$ for the RBF with modified Gaussian kernel (RBF\textsuperscript{4}).
    \item \emph{MBA}. The mean squared error start increasing after the at iteration number $6$, with an optimal tensor mesh containing $16\times16$ elements.
\end{itemize}

wQISA has the best generalization performances over the data set with respect to standard deviation, mean squared error and $L^\infty$ norm. This means that wQISA offers more accurate predictions on independent data. 

\begin{table}[!h]
\begin{center}
\caption{Rainfall fields, statistics for the error distribution for the LOO cross-validation. \label{table:rainfall_fields}}
\scalebox{0.9}{
\begin{tabular}{c m{0.5in} m{0.5in} m{0.5in} m{0.5in} m{0.5in}}
\midrule
Method & \vtop{\hbox{\strut Mean}\hbox{\strut [mm]}} & \vtop{\hbox{\strut Std}\hbox{\strut [mm]}} & \vtop{\hbox{\strut MSE}\hbox{\strut [mm\textsuperscript{2}]}} &
\vtop{\hbox{\strut $L^\infty$}\hbox{\strut [mm\textsuperscript{2}]}} &
\vtop{\hbox{\strut Numb. of}\hbox{\strut iterations}}  \\[-2ex]
\midrule
wQISA                   & 0.405 & \textbf{0.565} & \textbf{0.483} & \textbf{4.952} & \multicolumn{1}{|c}{4}\\ 
RBF\textsuperscript{1}  & 4.058 & 4.720 & 38.747 & 39.957  & \multicolumn{1}{|c}{114} \\
RBF\textsuperscript{2}  & \textbf{0.390} & 0.587 & 0.497 & 5.321  & \multicolumn{1}{|c}{20} \\
RBF\textsuperscript{3}  & 5.472 & 4.326 & 48.659  & 42.285 & \multicolumn{1}{|c}{132}\\
RBF\textsuperscript{4}  & 0.790 & 0.936 & 1.500  & 6.400 & \multicolumn{1}{|c}{26} \\
Kriging                 & 0.474 & 0.613 & 0.600 & 5.668 &  \multicolumn{1}{|c}{-} \\
MBA                     & 0.412 & 0.629 & 0.562 & 4.983 & \multicolumn{1}{|c}{5} \\
\midrule
\end{tabular}
}
\end{center}
\end{table}

The second event occurred on  January, 2014. The data set contains $19,187$ measurements, gathered from different devices: rain gauges and weather radar. Besides the $143$ rain gauges by Regione Liguria, another ~$25$ measure stations by Genoa municipality are considered. The remaining points come from raw radar acquisitions, at first as reflectivity measurements with a range of $400$ km. The frequency of mountains over the whole Ligurian territory affects the quality of radar acquisitions and a pre-processing step is needed to remove ground clutter effects; processed data are then  combined with observations gathered from rain gauges, which are more reliable measurements but do not cover the whole region. The integration of radar data in the interpolation of the precipitation field makes it possible to extend rainfall fields also to areas surrounding Liguria, and therefore to have a clearer picture about the temporal evolution of precipitation events. Since the temporal interval is different for each acquisition device, rainfall measurements have been cumulated. In this study, a $30$ minutes cumulative step has been used (240 time samples).
This procedure of integrating gauge and radar data was made to alleviate the well-known error and uncertainty that characterize radar estimates. Spurious signals may be caused, for example, by radar failure or by shielding of the radar beam by mountain ranges \cite{Harrison2000}. However, various outliers still perturbs the data.
The simulation over five time intervals is given in Table \ref{table:richer_rainfall_fields}.

\begin{table}[!h]
\begin{center}
\caption{Rainfall fields, statistics for the error distribution. \label{table:richer_rainfall_fields}}
\scalebox{0.9}{
\begin{tabular}{c m{0.55in} m{0.55in} m{0.55in} m{0.55in} m{0.55in} m{0.55in}}
\midrule
& Method & \vtop{\hbox{\strut Mean}\hbox{\strut [mm]}} & \vtop{\hbox{\strut Std}\hbox{\strut [mm]}} & \vtop{\hbox{\strut MSE}\hbox{\strut [mm\textsuperscript{2}]}}  &
\vtop{\hbox{\strut $L^\infty$}\hbox{\strut [mm]}}  &
\vtop{\hbox{\strut Numb. of}\hbox{\strut iterations}}  \\[-2ex]
\midrule

\multirow{7}{*}{00:30} & \multicolumn{1}{|l}{wQISA}  & \textbf{0.321} & \textbf{0.597} & \textbf{0.460} & 14.578 & \multicolumn{1}{|c}{7}\\

& \multicolumn{1}{|l}{RBF\textsuperscript{1}}  & 6.193 & 89.241 & $8\cdot 10^3$ & $4.16\cdot{}10^3$ & \multicolumn{1}{|c}{121} \\

& \multicolumn{1}{|l}{RBF\textsuperscript{2}}  & 0.625 & 0.924 & 1.244 & 14.356 & \multicolumn{1}{|c}{80}\\
& \multicolumn{1}{|l}{RBF\textsuperscript{3}}  & 1.549 & 30.776 & 949.591 & $1.55\cdot{}10^3$ & \multicolumn{1}{|c}{6} \\
& \multicolumn{1}{|l}{RBF\textsuperscript{4}}  & 1.135 & 1.570 & 3.753 & 15.600 & \multicolumn{1}{|c}{23} \\
& \multicolumn{1}{|l}{Kriging}                 & 0.958 & 1.248 & 2.476 & 14.473 & \multicolumn{1}{|c}{-} \\
& \multicolumn{1}{|l}{MBA}                     & 0.675 & 1.053 & 1.563 & \textbf{14.238} & \multicolumn{1}{|c}{5} \\
\midrule

\multirow{7}{*}{01:00} & \multicolumn{1}{|l}{wQISA}  & \textbf{0.270} & \textbf{0.589} & \textbf{0.419} & 11.125 & \multicolumn{1}{|c}{10}\\ 
& \multicolumn{1}{|l}{RBF\textsuperscript{1}}  & 130.886 & 225.675 & $6.80\cdot 10^4$ & $6.94\cdot{}10^3$ & \multicolumn{1}{|c}{122} \\
& \multicolumn{1}{|l}{RBF\textsuperscript{2}}  & 0.362 & 2.03 &   4.985 & 91.068 & \multicolumn{1}{|c}{122} \\
& \multicolumn{1}{|l}{RBF\textsuperscript{3}}  & 144.857 & 124.710 & $3.65\cdot 10^4$ & 903.952 & \multicolumn{1}{|c}{122} \\
& \multicolumn{1}{|l}{RBF\textsuperscript{4}}  & 0.790 & 0.936 & 1.500 & \textbf{6.400} & \multicolumn{1}{|c}{728} \\
& \multicolumn{1}{|l}{Kriging}                 & 0.964 & 1.233 & 2.245 & 12.449 &\multicolumn{1}{|c}{-} \\
& \multicolumn{1}{|l}{MBA}                     & 0.595 & 1.009 & 1.373 & 12.815 & \multicolumn{1}{|c}{5} \\
\midrule

\multirow{7}{*}{01:30} & \multicolumn{1}{|l}{wQISA}  & 0.305 & 0.647 & 0.512 & 10.792 & \multicolumn{1}{|c}{8}\\ 
& \multicolumn{1}{|l}{RBF\textsuperscript{1}}  & 
1.297 & 3.062 & 11.055 & 64.532 & \multicolumn{1}{|c}{6} \\
& \multicolumn{1}{|l}{RBF\textsuperscript{2}}  & 0.302 & 0.645 & \textbf{0.509} & 10.822 & \multicolumn{1}{|c}{35} \\
& \multicolumn{1}{|l}{RBF\textsuperscript{3}}  & 0.503 & 1.254 & 1.826 & 20.700 & \multicolumn{1}{|c}{7} \\
& \multicolumn{1}{|l}{RBF\textsuperscript{4}}  & 1.109 & 1.641 & 3.363 & 13.334 & \multicolumn{1}{|c}{24} \\
& \multicolumn{1}{|l}{Kriging}                 & 0.474 & \textbf{0.613} & 0.600 & \textbf{5.668} & \multicolumn{1}{|c}{-} \\
& \multicolumn{1}{|l}{MBA}                     & \textbf{0.299} & 0.653 & 0.516 & 10.810 & \multicolumn{1}{|c}{8} \\
\midrule

\multirow{7}{*}{02:00} & \multicolumn{1}{|l}{wQISA}  & \textbf{0.457} & \textbf{1.078} & \textbf{1.370} & 28.484 & \multicolumn{1}{|c}{7}\\ 
& \multicolumn{1}{|l}{RBF\textsuperscript{1}}  & 266.258 & 246.256 & $1.32\cdot 10^5$  & $2.69\cdot{}10^3$ & \multicolumn{1}{|c}{122} \\
& \multicolumn{1}{|l}{RBF\textsuperscript{2}}  & 0.749 & 4.389 &   19.822  & 101.409 & \multicolumn{1}{|c}{122} \\
& \multicolumn{1}{|l}{RBF\textsuperscript{3}}  & 369.316 & 326.457 & $2.43\cdot 10^5$ & \textbf{18.775} & \multicolumn{1}{|c}{137} \\
& \multicolumn{1}{|l}{RBF\textsuperscript{4}}  & 1.283 & 2.222 & 6.582 & 34.145 & \multicolumn{1}{|c}{12} \\
& \multicolumn{1}{|l}{Kriging}                 & 1.197 & 1.892 & 5.014  & 33.140 & \multicolumn{1}{|c}{-} \\
& \multicolumn{1}{|l}{MBA}                     & 0.916 & 1.560 & 3.271 & 31.4356 & \multicolumn{1}{|c}{5} \\
\midrule

\multirow{7}{*}{02:30} & \multicolumn{1}{|l}{wQISA}  & 0.423 & \textbf{0.974} & \textbf{1.127} & \textbf{16.069} & \multicolumn{1}{|c}{9}\\ 
& \multicolumn{1}{|l}{RBF\textsuperscript{1}}  & 36.954 & 42.373 & $3.16\cdot 10^3$ & 470.072 & \multicolumn{1}{|c}{124} \\
& \multicolumn{1}{|l}{RBF\textsuperscript{2}}  & 0.432 & 0.991 &   1.169 & $3.22\cdot{}10^3$ &  \multicolumn{1}{|c}{90} \\
& \multicolumn{1}{|l}{RBF\textsuperscript{3}}  & 0.850 & 3.112 & 10.408 & 128.381 & \multicolumn{1}{|c}{11} \\
& \multicolumn{1}{|l}{RBF\textsuperscript{4}}  & 1.544 & 2.933 & 10.984 & 23.973 & \multicolumn{1}{|c}{12} \\
& \multicolumn{1}{|l}{Kriging}                 & 1.776 & 2.334 & 8.600 & 22.421 & \multicolumn{1}{|c}{-} \\
& \multicolumn{1}{|l}{MBA}                     & \textbf{0.391} & 1.058 & 1.271 & 31.9741 & \multicolumn{1}{|c}{8} \\
\midrule

\end{tabular}
}
\end{center}
\end{table}

Again, wQISA provides the best generalization performances with respect to standard deviations and mean squared error. The difficulty of handling these datasets is reflected by the variability of the $L^\infty$ norm; because of the presence of outliers, the maximum of the absolute differences is largely unstable and, for each dataset, there is a different winner, even if the wQISA is generally one of the best performing methods also in with respect to this measure. The optimal parameters for the RBF approximations are provided in Table \ref{table:RBF_rain}. The optimal parameters $k$ for the wQISA with $k$-NN weight function are: $k=1$ (00:30), $k=2$ (01:00), $k=5$ (01:30), $k=1$ (00:30), $k=5$ (02:00) and $k=7$ (02:30). 

\begin{table}[!h]
\centering
\caption{Optimal smoothing parameter $\alpha$ for the RBFs implementations.}
\label{table:RBF_rain}
\scalebox{0.9}{
\begin{tabular}{lccccc}
\toprule
\textbf{Method} & \textbf{00:30} & \textbf{01:00} & \textbf{01:30} & \textbf{02:00} & \textbf{02:30}\\
\midrule
RBF\textsuperscript{1} & $0.00$ & $1.00$ & $0.68$ & $1.00$ & $0.00$\\
RBF\textsuperscript{2} & $439.64$ & $1.00$ & $0.19$ & $0.99$ & $0.99$ \\
RBF\textsuperscript{3} & $0.00$ & $1.00$ & $0.00$ & $0.59$ & $0.00$\\
RBF\textsuperscript{4} & $104.84$ & $1.00$ & $0.98$ & $0.24$ & $0.23$\\
\bottomrule
\end{tabular}
}
\end{table}

\section{Concluding remarks}
\label{sec:conclusions}
The weighted quasi interpolant spline approximation (wQISA) is a simple and robust procedure to obtain a spline approximation of point clouds. 
This paper presents a data-driven implementation of wQISA \cite{Raffo2019}, inspired to the supervised learning paradigm. 
The method has been tested on several real-world data from different application domains and different levels of noise.  Experiments have shown that wQISA is able to accurately generalize to data different from those used to define it, i.e. it can effectively describe geometric shapes and measurements not only at those points where it was defined.
In other words, wQISA can be successfully used also in prediction tasks, as a linear regression model in supervised learning.

When dealing with data affected by low level of noise, the wQISA method offers approximations that are comparable with other well-known methods.
However, the experiments show a better adherence of the wQISA surfaces to point clouds, as reflected by the smaller Hausdorff distance.
On these relatively error-free models, we have experimentally noticed a convergence with less iterations than the other methods considered. On the one hand, the wQISA method has linear convergence. On the other hand, this becomes an advantage when dealing with strongly perturbed data, as it allows to adapt better to the underlying distribution.
Indeed, on very noisy data wQISA shows a very good stability with respect to standard deviation and mean squared error. In the simulation of precipitation events, wQISA outperforms meshless methods (RBF and Kriging), which in \cite{Patane2017} showed the best approximation capability.
Experiments confirm that wQISA is able to successfully handle hundreds of thousands of points.

To sum up, two conclusions can be drawn from the comparative analysis. First: it is difficult and possibly ill-advised to search for an \emph{absolute} approximation method, because different approximation techniques can provide better results in different contexts and using different evaluation metrics. Second: for the studied point clouds, the experimental results show that data-driven wQISA exhibits a good prediction capability when it comes to strongly perturbed data. 

As a further development of the method, we think it is possible to consider adaptive refinement schemes, such as in the case of LR B-splines or THB-splines \cite{Dokken2013,BRACCO2018}.
This is particularly relevant because these locally refining schemes naturally deal with isogeometric computations and simulation and offers the valuable perspective to practically adopt this work for Computer Aided Design and Manufacturing (CAD/CAM), Finite Element Analysis and IsoGeometric Analysis \cite{Johannessen2014,GIANNELLI2016}. Finally, we intend to extend the comparative analysis to include other real-world applications, methods, and high dimensional data.

\section*{Acknowledgments} 

This project has received funding from the European Union’s Horizon 2020 research and innovation programme under the Marie Sk\l{}odowska-Curie grant agreement No 675789. The work has been partially developed in the CNR-IMATI activities DIT.AD021.080.001 and DIT.AD009.091.001. The authors also thanks: Dr. Bianca Falcidieno and Dr. Michela Spagnuolo for the fruitful discussions; Dr. Oliver J. D. Barrowclough, Dr. Tor Dokken and Dr. Georg Muntingh for their concern as supervisors; Dr. Vibeke Skytt for her help in collecting the terrain data sets; Laboratoire Domaines Ocaniques – CNRS UMR 6538 for the submarine sand dunes data set; Dr. Simone Pittaluga for the rainfall field data sets; Eng. Simone Cammarasana for his help in generating pictures with Delaunay triangulations.

\bibliographystyle{plain}
\bibliography{main.bbl}

\begin{thebibliography}{10}

\bibitem{STARC}
{STARC} repository.
\newblock \url{http://public.cyi.ac.cy/starcRepo/}.

\bibitem{Ackermann2015}
J.~Ackermann and M.~Goesele.
\newblock A survey of photometric stereo techniques.
\newblock {\em Found. Trends Comput. Graph. Vis.}, 9(3--4):149--254, 2015.

\bibitem{Alexa2001}
M.~Alexa, J.~Behr, D.~Cohen-Or, S.~Fleishman, D.~Levin, and C.T. Silva.
\newblock Point set surfaces.
\newblock In {\em {Proceedings of the conference on visualization 01, VIS 01,
  IEEE Computer Society}}, pages 21--28, 2001.

\bibitem{Berger:2017}
Matthew Berger, Andrea Tagliasacchi, Lee~M. Seversky, Pierre Alliez, Gaël
  Guennebaud, Joshua~A. Levine, Andrei Sharf, and Claudio~T. Silva.
\newblock A survey of surface reconstruction from point clouds.
\newblock {\em Computer Graphics Forum}, 36(1):301--329, 2017.

\bibitem{BRACCO2018}
Cesare Bracco, Carlotta Giannelli, David Großmann, and Alessandra Sestini.
\newblock Adaptive fitting with thb-splines: Error analysis and industrial
  applications.
\newblock {\em Comput. Aided Geom. Des.}, 62:239 -- 252, 2018.

\bibitem{Carr2001}
J.~C. Carr, R.~K. Beatson, J.~B. Cherrie, T.~J. Mitchell, W.~R. Fright, B.~C.
  McCallum, and T.~R. Evans.
\newblock {Reconstruction and representation of 3D objects with radial basis
  functions}.
\newblock In {\em {Proc. ACM SIGGRAPH}}, pages 67--76, 2001.

\bibitem{Cheney1995}
E.~W. Cheney.
\newblock {\em {Approximation Theory, Wavelets and Applications}}, volume 454,
  chapter Quasi-interpolation, pages 37--45.
\newblock {NATO Science Series}, 1995.

\bibitem{chui2016introduction}
Charles~K Chui.
\newblock {\em An introduction to wavelets}.
\newblock Elsevier, 2016.

\bibitem{Daehlen1993}
M.~D{\AE}hlen and M.~Floater.
\newblock Iterative polynomial interpolation and data compression.
\newblock {\em Numerical Algorithms}, 5(3):165--177, Mar 1993.

\bibitem{deboor73}
C~de~Boor and G.J Fix.
\newblock Spline approximation by quasi-interpolants.
\newblock {\em Journal of Approximation Theory}, 8(1):19 -- 45, 1973.

\bibitem{Dokken2013}
T.~Dokken, T.~Lyche, and K.~F. Pettersen.
\newblock Polynomial splines over locally refined box-partitions.
\newblock {\em {Comput. Aided Geom. Des.}}, 30(3):331--356, 2013.

\bibitem{Farin1996}
Gerald~E. Farin.
\newblock {\em Curves and Surfaces for Computer-Aided Geometric Design: A
  Practical Code}.
\newblock Academic Press, Inc., USA, 4th edition, 1996.

\bibitem{Feng2014}
W.~Feng, Z.~Yang, and J.~Deng.
\newblock Moving multiple curves/surfaces approximation of mixed point clouds.
\newblock {\em Commun Math Stat}, 2(1):107--124, 2014.

\bibitem{Fleishman2005}
Shachar Fleishman, Daniel Cohen-Or, and Cl\'{a}udio~T. Silva.
\newblock Robust moving least-squares fitting with sharp features.
\newblock {\em ACM Trans. Graph.}, 24(3):544–552, July 2005.

\bibitem{FLOATER1996}
Michael~S. Floater and Armin Iske.
\newblock Multistep scattered data interpolation using compactly supported
  radial basis functions.
\newblock {\em Journal of Computational and Applied Mathematics}, 73(1):65 --
  78, 1996.

\bibitem{Forsey1995}
David~R. Forsey and Richard~H. Bartels.
\newblock Surface fitting with hierarchical splines.
\newblock {\em ACM Trans. Graph.}, 14(2):134–161, April 1995.

\bibitem{Franke}
Richard Franke.
\newblock Scattered data interpolation: Tests of some method.
\newblock {\em Mathematics of Computation}, 38(157):181--200, 1982.

\bibitem{FRANZETTI201317}
Marcaurelio Franzetti, Pascal~Le Roy, Christophe Delacourt, Thierry Garlan,
  Romain Cancouët, Alexey Sukhovich, and Anne Deschamps.
\newblock Giant dune morphologies and dynamics in a deep continental shelf
  environment: {Example of the banc du four (Western Brittany, France)}.
\newblock {\em Marine Geology}, 346:17 -- 30, 2013.

\bibitem{Friedman:1977}
Jerome~H. Friedman, Jon~Louis Bentley, and Raphael~Ari Finkel.
\newblock An algorithm for finding best matches in logarithmic expected time.
\newblock {\em ACM Trans. Math. Softw.}, 3(3):209--226, September 1977.

\bibitem{garnero2013comparisons}
Gabriele Garnero.
\newblock Comparisons between different interpolation techniques.
\newblock {\em International Archives of the Photogrammetry, Remote Sensing and
  Spatial Information Sciences}, 5:W3, 2013.

\bibitem{Georgopoulos2017}
A.~Georgopoulos.
\newblock {\em Heritage and Archaeology in the Digital Age}.
\newblock Springer, 2017.

\bibitem{GIANNELLI2016}
C.~Giannelli, B.~J\:{u}ttler, S.~K. Kleiss, A.~Mantzaflaris, B.~Simeon, and
  J.~\v{S}peh.
\newblock {THB}-splines: An effective mathematical technology for adaptive
  refinement in geometric design and isogeometric analysis.
\newblock {\em Comput Methods Appl Mech Eng}, 299:337 -- 365, 2016.

\bibitem{Giraudot:2013}
Simon Giraudot, David Cohen-Steiner, and Pierre Alliez.
\newblock Noise-adaptive shape reconstruction from raw point sets.
\newblock {\em Computer Graphics Forum}, 32(5):229--238, 2013.

\bibitem{Greiner1997}
Gnther Greiner and Kai Hormann.
\newblock Interpolating and approximating scattered {3D} data with hierarchical
  tensor product {B}-splines, surface fitting and multiresolution methods,
  1997.

\bibitem{Harrison2000}
D~L Harrison, S~J Driscoll, and M~Kitchen.
\newblock Improving precipitation estimates from weather radar using quality
  control and correction techniques.
\newblock {\em Meteorological Applications}, 7(2):135--144, 2000.

\bibitem{Hastie2009}
T.~Hastie, R.~Tibshirani, and J.~Friedman.
\newblock {\em The Elements of Statistical Learning}.
\newblock Springer, second edition, 2009.

\bibitem{Johannessen2014}
K.~A. Johannessen, T.~Kvamsdal, and T.~Dokken.
\newblock Isogeometric analysis using {LR B}-splines.
\newblock {\em Comput Methods Appl Mech Eng}, 2014.

\bibitem{Kazhdan2006}
M.~Kazhdan, M.~Bolitho, and H.~Hoppe.
\newblock Poisson surface reconstruction.
\newblock In {\em {$4^{th}$ EG Symp. on Geometry Processing}}, pages 61--70,
  2006.

\bibitem{Kiss2014}
Gábor Kiss, Carlotta Giannelli, Urška Zore, Bert Jüttler, David Großmann,
  and Johannes Barner.
\newblock Adaptive {CAD} model (re-)construction with {THB}-splines.
\newblock {\em Graphical Models}, 76(5):273 -- 288, 2014.

\bibitem{Lee1997}
S.~Lee, G.~Wolberg, and S.~Y. Shin.
\newblock Scattered data interpolation with multilevel b-splines.
\newblock {\em {IEEE T Vis Comput Gr}}, 3(3):228--244, 1997.

\bibitem{Levin2003}
D.~Levin.
\newblock {\em Geometric modeling for scientific visualization}, chapter
  Mesh-independent surface interpolation, pages 37--49.
\newblock Springer-Verlag, 2003.

\bibitem{Lyche2011}
T.~Lyche and K.~M{\o}rken.
\newblock {\em Spline Methods Draft}, chapter Tensor Product Spline Surfaces,
  pages 149--166.
\newblock University of Oslo, April 2011.

\bibitem{Ohtake2003}
Yutaka Ohtake, Alexander Belyaev, Marc Alexa, Greg Turk, and Hans-Peter Seidel.
\newblock Multi-level partition of unity implicits.
\newblock {\em ACM Trans. Graph.}, 22(3):463–470, July 2003.

\bibitem{Oliver1990}
M.~A. Oliver and R.~Webster.
\newblock Kriging: a method of interpolation for geographical information
  systems.
\newblock {\em International Journal of Geographical Information Systems},
  4(3):313--332, 1990.

\bibitem{Patane2017}
G.~Patan{\'e}, A.~Cerri, V.~Skytt, S.~Pittaluga, S.~Biasotti, D.~Sobrero,
  T.~Dokken, and M.~Spagnuolo.
\newblock Comparing methods for the approximation of rainfall fields in
  environmental applications.
\newblock {\em {ISPRS J. Photogramm. and Remote Sens.}}, 127:57--72, 2017.

\bibitem{PATANE2012387}
Giuseppe Patanè and Michela Spagnuolo.
\newblock Local approximation of scalar functions on {3D} shapes and volumetric
  data.
\newblock {\em Comput Graph-UK}, 36(5):387 -- 397, 2012.

\bibitem{Raffo2019}
A.~Raffo and S.~Biasotti.
\newblock {Weighted Quasi Interpolant Spline Approximations: Properties and
  Applications}, 2019.

\bibitem{Reitermanova2010}
Z~Reitermanova.
\newblock Data splitting.
\newblock In {\em WDS}, volume~10, pages 31--36, 2010.

\bibitem{Sablonniere05}
Paul Sablonni{\`e}re.
\newblock Recent progress on univariate and multivariate polynomial and spline
  quasi-interpolants.
\newblock In Detlef~H. Mache, J{\'o}zsef Szabados, and Marcel~G. de~Bruin,
  editors, {\em Trends and Applications in Constructive Approximation}, pages
  229--245. Birkh{\"a}user Basel, 2005.

\bibitem{Savchenko1995}
V.~V. Savchenko, A.~A. Pasko, O.~G. Okunev, and T.~L. Kunii.
\newblock Function representation of solids reconstructed from scattered
  surface points and contours.
\newblock {\em {Computer Graphics Forum}}, 14(4):181--188, 1995.

\bibitem{Schumaker2007}
L.~L. Schumaker.
\newblock {\em Spline Functions: Basic Theory}.
\newblock Cambridge University Press, 2007.

\bibitem{Sederberg2003}
T.~W. Sederberg, J.~Zheng, A.~Bakenov, and A.~Nasri.
\newblock {T-Splines and T-NURCCS}.
\newblock {\em {ACM Trans. Graph.}}, 22(3):477--484, 2003.

\bibitem{Shen2004}
C.~Shen, J.F. O{\textquotesingle}Brien, and J.R. Shewchuk.
\newblock Interpolating and approximating implicit surfaces from polygon soup.
\newblock {\em ACM Trans. Graph.}, 23(3):896--904, August 2004.

\bibitem{Skytt2015}
V.~Skytt, O.~Barrowclough, and T.~Dokken.
\newblock Locally refined spline surfaces for representation of terrain data.
\newblock {\em {Computer \& Graphics}}, 49:58--68, 2015.

\bibitem{Sorgente2018}
T.~Sorgente, S.~Biasotti, M.~Livesu, and M.~Spagnuolo.
\newblock Topology-driven shape chartification.
\newblock {\em Comput. Aided Geom. Des.}, 65:13 -- 28, 2018.

\end{thebibliography}

\appendix
\section{Additional material\label{appendix}}

When dealing with particularly small data sets, wQISA can still be used for single polynomial approximation. As an example, we consider the 2017 Particulate Matter (PM10) report for the city of Oslo (Norway), which is freely available at the municipal data-bank\footnote{http://statistikkbanken.oslo.kommune.no/webview/}. Particulates, also known as suspended particulate matter, are microscopic particles of solid or liquid matter suspended in the air. They are classified by size, such as for the particles smaller than 10 $\mu{}m$, called PM10, and the particles less than $2.5\mu{}m$, called PM2.5. PM10 contains particles that originate from combustion (e.g., heating) and road dust (e.g., car tyres and brakes). As particulates are among the most harmful form of air pollution, their approximation is of great interest.

\begin{figure}[!h]
    \centering
    \includegraphics[scale=0.5]{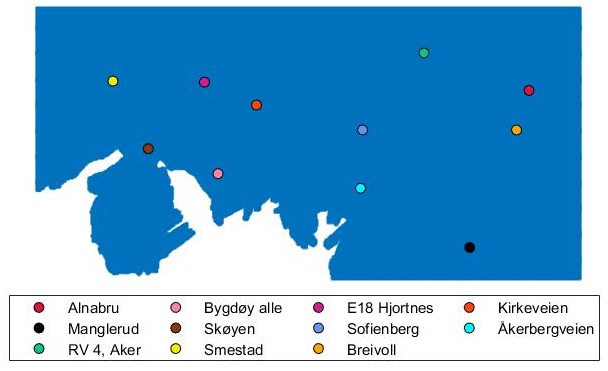}
    \caption{The 11 stations monitoring the PM10 concentration, owned by \emph{Oslo kommune} and \emph{Statens vegvesen}.}
    \label{figure:air_pollution_map}
\end{figure}

The data set contains only $11$ samples; see in Figure \ref{figure:air_pollution_map} an Oslo map with their position. The evaluation results give insight into the approximation performance in real condition of sparsity and provide an example of a situation where there is insufficient data to split into three parts. 

Given the small size of the data set, a Leave-One-Out cross-validation was performed to fix the free parameters and compute the statistics for the punctual error. The optimal parameters were computed by minimizing the mean squared error:
\begin{itemize}
    \item \emph{wQISA with k-NN weight}. Given the low number of points, we fix the restrict our attention to a degenerate tensor mesh containing just one element. The optimal number of neighbors $k$ is $k=5$, where we have checked for any $k=1,\ldots,10$.
    \item \emph{RBF}. The optimal smoothing parameter $\alpha$ for RBF approximations is computed by minimizing the mean squared error: $\alpha=679.69{\cdot}10^{-3}$ for Gaussian kernel (RBF\textsuperscript{1}) in $53$ iterations; $\alpha=190.06{\cdot}10^{-3}$ for RBF with multiquadric kernel (RBF\textsuperscript{2}) in $51$ iterations; $\alpha=25.00{\cdot}10^{-6}$ for RBF with inverse kernel (RBF\textsuperscript{3}) in $7$ iterations; $\alpha=983.44{\cdot}10^{-3}$ for RBF with modified Gaussian kernel (RBF\textsuperscript{4}) in $51$ iterations.
    \item \emph{MBA}. The mean squared error start increasing after the first iteration, thus we adopt a tensor mesh containing just one element as for wQISA.
\end{itemize}

 The results are presented in Table \ref{table:air_pollution}. The three best performances with respect to the MSE are given by MBA, wQISA and RBF with Gaussian kernel.

\begin{table}[!h]
\begin{center}
\caption{Air pollution, statistics for the error distribution for the LOO cross-validation. \label{table:air_pollution}}
\scalebox{0.9}{
\begin{tabular}{m{0.5in} m{0.5in} m{0.5in} m{0.5in} m{0.5in} m{0.5in}}
\midrule
Method & \vtop{\hbox{\strut Mean}\hbox{\strut [mm]}} & 
\vtop{\hbox{\strut Std}\hbox{\strut [mm]}} & \vtop{\hbox{\strut MSE}\hbox{\strut [mm\textsuperscript{2}]}} & \vtop{\hbox{\strut $L^\infty$}\hbox{\strut [mm]}}  \\[-2ex]
\midrule
wQISA                   & 3.005 & \textbf{1.984} & 12.967 & 7.245\\ 
Kriging                  & 3.365 & 2.020 & 15.401 & 7.068 \\
RBF\textsuperscript{1}  & 7.626 & 3.346 & 77.042 & 16.525\\
RBF\textsuperscript{2}  & 3.320 & 2.260 & 16.130 & 7.921 \\
RBF\textsuperscript{3}  & 3.501 & 2.531 & 18.660 & 7.704 \\
RBF\textsuperscript{4} & 14.731 & 3.714 & 230.801 & 21.193\\
MBA                  & \textbf{2.576} & 2.378 & \textbf{11.775} & \textbf{6.784} \\

\midrule
\end{tabular}
}
\end{center}
\end{table}

\end{document}